\documentclass[preprint,12pt]{elsarticle}
\usepackage{float,epsfig,epstopdf,bm,floatflt,here,relsize}
\usepackage{rotating}
\usepackage{amssymb}
\usepackage{amsmath}
\usepackage{mathtools}
\usepackage{tikz}
\usepackage{multirow}
\usepackage{booktabs}
\usepackage{array}
\usepackage{graphicx}\newcolumntype{H}{>{\setbox0=\hbox\bgroup}c<{\egroup}@{}}
\allowdisplaybreaks

\journal{Journal of Computational Physics}

\begin{document}

\begin{frontmatter}
\title{\textbf{Phase error analysis of implicit Runge-Kutta methods: Introducing new classes of minimal dissipation low dispersion high order schemes}}

\author{Subhajit Giri}
\ead{subhajit@tezu.ernet.in}
\author{Shuvam Sen\corref{cor1}}
\ead{shuvam@tezu.ernet.in}
\cortext[cor1]{Corresponding author}
\address{Department of Mathematical Sciences, Tezpur University, Tezpur 784028, INDIA}

\begin{abstract}
In current research, we analyse dissipation and dispersion characteristics of most accurate two and three stage Gauss-Legendre implicit Runge-Kutta (R-K) methods. These methods, known for their $A$-stability and immense accuracy, are observed to carry minimum dissipation error along with highest possible dispersive order in their respective classes. We investigate to reveal that these schemes are inherently optimized to carry low phase error only at small wavenumber. As larger temporal step size is imperative in conjunction with implicit R-K methods for physical problems, we interpret to derive a class of minimum dissipation and optimally low dispersion implicit R-K schemes. Schemes thus obtained by cutting down amplification error and maximum reduction of weighted phase error, suggest better accuracy for relatively bigger CFL number. Significantly, we are able to outline an algorithm that can be used to design stable implicit R-K methods for suitable time step with better accuracy virtues. The algorithm is potentially generalizable for implicit R-K class of methods. As we focus on two and three stage schemes a comprehensive comparison is carried out using numerical test cases.
\end{abstract}

\end{frontmatter}

\section{Introduction}
The general Runge-Kutta (R-K) methods widely known as implicit R-K methods are more challanging vis-a-vis explicit R-K methods. This dilemma stems from the necessity of iterative computations in implicit schemes as opposed to stage-by-stage implementation possible for explicit methods. However, as pointed out by Butcher \cite{but_08} there are compelling reasons to study them from theoretical and practical point of view. One of the reasons for interest in implicit R-K schemes lies in their weak stability characteristics, which are superior to those of explicit schemes. From practical view point, efficient solution of stiff problems require implicit R-K methods. Alexander \cite{ale_77} has noted in his work that for stiff problems only $A$-stable implicit R-K methods are useful. Explicit methods often suffer from stability limitations and as such small temporal step size become imperative. In numerical acoustics small time step lead to excessive computational cost. An $A$ stable implicit methods allows one to compute with bigger step size thereby somewhat compensating for the additional time spent on each step. Moreover even with lower stage number an implicit method can acquire higher order of accuracy compared to explicit methods. Maximum attainable order for various classes of $R$-stage implicit methods can be found in the works of Butcher \cite{but_08} and Alexander \cite{ale_77}. In this context Gauss-Legendre methods based on Legendre polynomials are known for their prolific accuracy. A $R$-stage implicit Gauss-Legendre method admit an order of accuracy as high as $2R$. Nevertheless, it must be said that ability of compute stiff problems that too with a relatively  bigger time step is the main motivation for implicit R-K methods.

It is well known that for unsteady flow problems, a convergent high order scheme does not guarantee good quality numerical solution \cite{tam_web_93}. Numerical schemes must be able to resolve all scales present in the flow. Further physical propagation speed of the respective scale should be appropriately matched by the numerical propagation speed of individually resolved scale. Hence for various problems in CFD, especially in computational acoustic, there has been historical interest for low-dissipation low-dispersion schemes. Altough a plethora of efficient low-dissipation and low-dispersion spatial discretization schemes are available in the literature, studies towards dissipation and dispersion relation preserving temporal integration procedures have been mostly confined to explicit Runge-Kutta algorithms. In this context works of Simos \cite{sim_93}, Hu et al. \cite{hu_hus_man_96}, Calvo et al. \cite{cal_fra_ran_03}, Bogey and Bailly \cite{bog_bai_04}, Berland et al. \cite{ber_bog_bai_06}, Anastassi and Simos \cite{ana_sim_05} and Tselios and Simos \cite{tse_sim_05} deserve special mention. For implicit R-K methods low-dissipation low-dispersion schemes have been advocated starting with the pioneering work of Franco et al. \cite{fra_gom_ran_97}. This has been followed more recently by the works of Najafi-Yazdi and Mongeau \cite{naj_mon_13} and Nazari et al. \cite{naz_moh_cha_14, naz_moh_cha_15}. But these works have been limited to diagonally implicit schemes. Although good stability characteristics can be found in diagonally implicit R-K methods \cite{fer_spi_08}, but it is clear from the above studies that the phase-lag virtues of diagonally implicit schemes are a compromise between explicit and fully implicit methods. To the best of our knowledge systematic study towards development of low-dissipation low-dispersion fully implicit R-K methods is not available in the literarture. Here it must be said that Bhaumik et al. \cite{bha_sen_sen_13} in their work have documented numerical properties of two stage implicit Gauss-Legendre method in conjunction with various spatial discretization procedures.

One of the aim of this work is to explore dispersion properties of maximally accurate two and three stage Gauss-Legendre implicit R-K schemes. We look to notify inherent wave resolving properties of these schemes known for their high order of accuracy. In the process a new algorithm which minimizes amplification error and optimally reduces weighted phase error is suggested. This algorithm used for implicit R-K class of methods, render a scheme $A$-stable and provide better numerical accuracy even for relatively larger time steps. This easy to implement algorithm is used to propose a class of schemes with special emphasis on certain regions of CFL number thereby making low-dissipation low-dispersion computation possible even with relatively bigger temporal step size.

We organize remainder of this paper into six sections. In Section 2, implicit R-K schemes are briefly introduced and its association to prototype equations are inspected. Analysis and derivation of new two and three stage minimal dissipation low dispersion schemes are done in Section 3 and Section 4 respectively. The algorithm is summarized in Section 5. Numerical examples are given in Section 6 and finally in Section 7 concluding remarks are offered.

\section{Implicit Runge-Kutta schemes}
Consider initial value problem (IVP) 
\begin{eqnarray}\label{1}
\frac{du}{dt}=f(t,u),\;\;\;u(t_0)=u_0.
\end{eqnarray}
A general $R$-stage R-K method can be defined as
\begin{eqnarray}\label{2}
u^{n+1}=u^n+\Delta t\sum_{r=1}^{R}b_rF_r
\end{eqnarray}
where
\begin{eqnarray}
F_r=f\bigg(t^n+\Delta t c_r, u^n+\Delta t\sum_{s=1}^Ra_{rs}F_s\bigg),\;\;\;r=1,2,...,R,\label{3}\\
c_r=\sum_{s=1}^Ra_{rs},\;\;\;\;\;\;\;\;\;\;\;\;\;\;\;\;\;\;\;\;\;\;\;\;\;\;\;\;\;\;\;\;\;\;\;\;\;\;\;\;\;\;\;r=1,2,...,R.\label{4}
\end{eqnarray}
Using Butcher tableau \cite{but_08} the above methods can be represented as
\begin{eqnarray}\label{5}
\begin{tabular}{c|c}
$\bm c$& $\bm A$\\
  \hline
       & $\bm b^T$ \\
\end{tabular}.
\end{eqnarray}
where $\bm A=(a_{rs})_{R\times R}$, $\bm b=(b_r)^T_R$, $\bm c=(c_r)^T_R$. For explicit R-K scheme $\bm A$ is strictly lower triangular and for diagonally implicit schemes $\bm A$ is lower triangular with non-zero diagonal entries. At times a diagonally implicit scheme is further categorized as singly diagonally implicit if its all diagonal entries are equal.

For R-K class of methods accuracy conditions upto fourth order can be represented as shown below.
\begin{eqnarray}
O(\Delta t):\;\;\; &\begin{tikzpicture}
     				\draw  node[fill,circle,scale=0.3]{} (0,0);
  				    \end{tikzpicture}&							\;\;\sum_{r=1}^{R}b_r=1,\label{6}\\
O(\Delta t^2):\;\;\;&\begin{tikzpicture}
  					 [scale=.3,auto=left,every node/.style={circle,scale=0.4,fill=black!100}]
  					 \node (n1) at (0,0) {};
  					 \node (n2) at (0,1)  {};
  					 \foreach \from/\to in {n1/n2}
   					 \draw (\from) -- (\to);
					\end{tikzpicture}&  \;\;\sum_{r=1}^{R}b_rc_r=\frac{1}{2},\label{8}\\
O(\Delta t^3):\;\;\;&\begin{tikzpicture}
  					 [scale=.3,auto=left,every node/.style={circle,scale=0.4,fill=black!100}]
  					 \node (n1) at (0,0) {};
  					 \node (n2) at (-0.5,1)  {};
  					 \node (n3) at (0.5,1)  {};
  					 \foreach \from/\to in {n1/n2,n1/n3}
   					 \draw (\from) -- (\to);
					\end{tikzpicture}&	\;\;\sum_{r=1}^{R}b_rc_r^2=\frac{1}{3},\label{9}\\
					&\begin{tikzpicture}
  					 [scale=.3,auto=left,every node/.style={circle,scale=0.4,fill=black!100}]
  					 \node (n1) at (0,0) {};
  					 \node (n2) at (0,1)  {};
  					 \node (n3) at (0,2)  {};
  					 \foreach \from/\to in {n1/n2,n2/n3}
   					 \draw (\from) -- (\to);
					\end{tikzpicture}& 	\;\;\sum_{r=1}^{R}\sum_{s=1}^{R}b_ra_{rs}c_s=\frac{1}{6},\label{10}\\
O(\Delta t^4):\;\;\;&\begin{tikzpicture}
  					 [scale=.3,auto=left,every node/.style={circle,scale=0.4,fill=black!100}]
  					 \node (n1) at (0,0) {};
  					 \node (n2) at (-1,1)  {};
  					 \node (n3) at (0,1)  {};
  					 \node (n4) at (1,1)  {};
  					 \foreach \from/\to in {n1/n2,n1/n3,n1/n4}
   					 \draw (\from) -- (\to);
					\end{tikzpicture}& \;\;\sum_{r=1}^{R}b_rc_r^3=\frac{1}{4},\label{11} \\
					&\begin{tikzpicture}
  					 [scale=.3,auto=left,every node/.style={circle,scale=0.4,fill=black!100}]
  					 \node (n1) at (0,0) {};
  					 \node (n2) at (-0.5,1)  {};
  					 \node (n3) at (0.5,1)  {};
  					 \node (n4) at (0.5,2)  {};
  					 \foreach \from/\to in {n1/n2,n1/n3,n3/n4}
   					 \draw (\from) -- (\to);
					\end{tikzpicture}& \;\;\sum_{r=1}^{R}\sum_{s=1}^{R}b_rc_ra_{rs}c_s=\frac{1}{8},\label{12}\\
					&\begin{tikzpicture}
  					 [scale=.3,auto=left,every node/.style={circle,scale=0.4,fill=black!100}]
  					 \node (n1) at (0,0) {};
  					 \node (n2) at (0,1)  {};
  					 \node (n3) at (-0.5,2)  {};
  					 \node (n4) at (0.5,2)  {};
  					 \foreach \from/\to in {n1/n2,n2/n3,n2/n4}
   					 \draw (\from) -- (\to);
					\end{tikzpicture}& \;\;\sum_{r=1}^{R}\sum_{s=1}^{R}b_ra_{rs}c_s^2=\frac{1}{12},\label{13}\\
					&\begin{tikzpicture}
  					 [scale=.3,auto=left,every node/.style={circle,scale=0.4,fill=black!100}]
  					 \node (n1) at (0,0) {};
  					 \node (n2) at (0,0.75)  {};
  					 \node (n3) at (0,1.5)  {};
  					 \node (n4) at (0,2.25)  {};
  					 \foreach \from/\to in {n1/n2,n2/n3,n3/n4}
   					 \draw (\from) -- (\to);
					\end{tikzpicture}& \;\;\sum_{r=1}^{R}\sum_{s=1}^{R}\sum_{l=1}^{R}b_ra_{rs}a_{sl}c_l=\frac{1}{24}.\label{14}
\end{eqnarray}
Rooted trees shown against Eqs. (\ref{6}), (\ref{8})-(\ref{14}) are pictorial representations of various order conditions \cite{but_08}. All of which may not be independent. In the above only two of the three Eqs. (\ref{12})-(\ref{14}) are independent. Fifth and sixth order accuracy require satisfaction of additional nine and twenty conditions respectively. Using rooted trees these are inscribed as		$$\begin{tikzpicture}
  					 [scale=.15,auto=left,every node/.style={circle,scale=0.2,fill=black!100}]
  					 \node (n1) at (0,0) {};
  					 \node (n2) at (-1.2,1)  {};
  					 \node (n3) at (-0.4,1)  {};
  					 \node (n4) at (0.4,1)  {};
  					 \node (n5) at (1.2,1)  {};
  					 \foreach \from/\to in {n1/n2,n1/n3,n1/n4,n1/n5}
   					 \draw (\from) -- (\to);
					\end{tikzpicture},
					\begin{tikzpicture}
  					 [scale=.15,auto=left,every node/.style={circle,scale=0.2,fill=black!100}]
  					 \node (n1) at (0,0) {};
  					 \node (n2) at (-1,1)  {};
  					 \node (n3) at (0,1)  {};
  					 \node (n4) at (1,1)  {};
  					 \node (n5) at (1,2)  {};
  					 \foreach \from/\to in {n1/n2,n1/n3,n1/n4,n4/n5}
   					 \draw (\from) -- (\to);
					\end{tikzpicture},
					\begin{tikzpicture}
  					 [scale=.15,auto=left,every node/.style={circle,scale=0.2,fill=black!100}]
  					 \node (n1) at (0,0) {};
  					 \node (n2) at (-0.5,1)  {};
  					 \node (n3) at (0.5,1)  {};
  					 \node (n4) at (0,2)  {};
  					 \node (n5) at (1,2)  {};
  					 \foreach \from/\to in {n1/n2,n1/n3,n3/n4,n3/n5}
   					 \draw (\from) -- (\to);
					\end{tikzpicture},
					\begin{tikzpicture}
  					 [scale=.15,auto=left,every node/.style={circle,scale=0.2,fill=black!100}]
  					 \node (n1) at (0,0) {};
  					  \node (n2) at (-0.5,1)  {};
  					 \node (n3) at (0.5,1)  {};
  					 \node (n4) at (0.5,2)  {};
  					 \node (n5) at (0.5,3)  {};
  					 \foreach \from/\to in {n1/n2,n1/n3,n3/n4,n4/n5}
   					 \draw (\from) -- (\to);
					\end{tikzpicture},
					\begin{tikzpicture}
  					 [scale=.15,auto=left,every node/.style={circle,scale=0.2,fill=black!100}]
  					 \node (n1) at (0,0) {};
  					 \node (n2) at (-0.5,1)  {};
  					 \node (n3) at (0.5,1)  {};
  					 \node (n4) at (-0.5,2)  {};
  					 \node (n5) at (0.5,2)  {};
  					 \foreach \from/\to in {n1/n2,n1/n3,n2/n4,n3/n5}
   					 \draw (\from) -- (\to);
					\end{tikzpicture},
					\begin{tikzpicture}
  					 [scale=.15,auto=left,every node/.style={circle,scale=0.2,fill=black!100}]
  					 \node (n1) at (0,0) {};
  					 \node (n2) at (0,1)  {};
  					 \node (n3) at (0,2)  {};
  					 \node (n4) at (-1,2)  {};
  					 \node (n5) at (1,2)  {};
  					 \foreach \from/\to in {n1/n2,n2/n3,n2/n4,n2/n5}
   					 \draw (\from) -- (\to);
					\end{tikzpicture},
					\begin{tikzpicture}
  					 [scale=.15,auto=left,every node/.style={circle,scale=0.2,fill=black!100}]
  					 \node (n1) at (0,0) {};
  					 \node (n2) at (0,1)  {};
  					 \node (n3) at (-0.5,2)  {};
  					 \node (n4) at (0.5,2)  {};
  					 \node (n5) at (0.5,3)  {};
  					 \foreach \from/\to in {n1/n2,n2/n3,n2/n4,n4/n5}
   					 \draw (\from) -- (\to);
					\end{tikzpicture},
					\begin{tikzpicture}
  					 [scale=.15,auto=left,every node/.style={circle,scale=0.2,fill=black!100}]
  					 \node (n1) at (0,0) {};
  					 \node (n2) at (0,1)  {};
  					 \node (n3) at (0,2)  {};
  					 \node (n4) at (-0.5,3)  {};
  					 \node (n5) at (0.5,3)  {};
  					 \foreach \from/\to in {n1/n2,n2/n3,n3/n4,n3/n5}
   					 \draw (\from) -- (\to);
					\end{tikzpicture},
					\begin{tikzpicture}
  					 [scale=.15,auto=left,every node/.style={circle,scale=0.2,fill=black!100}]
  					 \node (n1) at (0,0) {};
  					 \node (n2) at (0,1)  {};
  					 \node (n3) at (0,2)  {};
  					 \node (n4) at (0,3)  {};
  					 \node (n5) at (0,4)  {};
  					 \foreach \from/\to in {n1/n2,n2/n3,n3/n4,n4/n5}
   					 \draw (\from) -- (\to);
					\end{tikzpicture}$$
					 and
					$$\begin{tikzpicture}
  					 [scale=.15,auto=left,every node/.style={circle,scale=0.2,fill=black!100}]
  					 \node (n1) at (0,0)   {};
  					 \node (n2) at (-2,1)  {};
  					 \node (n3) at (-1,1)  {};
  					 \node (n4) at (0,1)   {};
  					 \node (n5) at (1,1)   {};
  					 \node (n6) at (2,1)   {};
  					 \foreach \from/\to in {n1/n2,n1/n3,n1/n4,n1/n5,n1/n6}
   					 \draw (\from) -- (\to);
					\end{tikzpicture},
					\begin{tikzpicture}
  					 [scale=.15,auto=left,every node/.style={circle,scale=0.2,fill=black!100}]
  					 \node (n1) at (0,0)     {};
  					 \node (n2) at (-1.2,1)  {};
  					 \node (n3) at (-0.4,1)  {};
  					 \node (n4) at (0.4,1)   {};
  					 \node (n5) at (1.2,1)   {};
  					 \node (n6) at (1.2,2)   {};
  					 \foreach \from/\to in {n1/n2,n1/n3,n1/n4,n1/n5,n5/n6}
   					 \draw (\from) -- (\to);
					\end{tikzpicture},
					\begin{tikzpicture}
  					 [scale=.15,auto=left,every node/.style={circle,scale=0.2,fill=black!100}]
  					 \node (n1) at (0,0)     {};
  					 \node (n2) at (-0.8,1)  {};
  					 \node (n3) at (-0.8,2)  {};
  					 \node (n4) at (0,1)     {};
  					 \node (n5) at (0.8,1)   {};
  					 \node (n6) at (0.8,2)   {};
  					 \foreach \from/\to in {n1/n2,n2/n3,n1/n4,n1/n5,n5/n6}
   					 \draw (\from) -- (\to);
					\end{tikzpicture},
					\begin{tikzpicture}
  					 [scale=.15,auto=left,every node/.style={circle,scale=0.2,fill=black!100}]
  					 \node (n1) at (0,0)     {};
  					 \node (n2) at (-0.8,1)  {};
  					 \node (n3) at (0,1)     {};
  					 \node (n4) at (0.8,1)   {};
  					 \node (n5) at (0.5,2)   {};
  					 \node (n6) at (1.5,2)   {};
  					 \foreach \from/\to in {n1/n2,n1/n3,n1/n4,n4/n5,n4/n6}
   					 \draw (\from) -- (\to);
					\end{tikzpicture},
					\begin{tikzpicture}
  					 [scale=.15,auto=left,every node/.style={circle,scale=0.2,fill=black!100}]
  					 \node (n1) at (0,0)     {};
  					 \node (n2) at (-0.8,1)  {};
  					 \node (n3) at (0,1)     {};
  					 \node (n4) at (0.8,1)   {};
  					 \node (n5) at (0.8,2)   {};
  					 \node (n6) at (0.8,3)   {};
  					 \foreach \from/\to in {n1/n2,n1/n3,n1/n4,n4/n5,n5/n6}
   					 \draw (\from) -- (\to);
					\end{tikzpicture},
					\begin{tikzpicture}
  					 [scale=.15,auto=left,every node/.style={circle,scale=0.2,fill=black!100}]
  					 \node (n1) at (0,0)     {};
  					 \node (n2) at (-0.8,1)  {};
  					 \node (n3) at (0.8,1)   {};
  					 \node (n4) at (0,2)     {};
  					 \node (n5) at (0.8,2)   {};
  					 \node (n6) at (1.6,2)   {};
  					 \foreach \from/\to in {n1/n2,n1/n3,n3/n4,n3/n5,n3/n6}
   					 \draw (\from) -- (\to);
					\end{tikzpicture},
					\begin{tikzpicture}
  					 [scale=.15,auto=left,every node/.style={circle,scale=0.2,fill=black!100}]
  					 \node (n1) at (0,0)     {};
  					 \node (n2) at (-0.8,1)  {};
  					 \node (n3) at (0.8,1)   {};
  					 \node (n4) at (0,2)     {};
  					 \node (n5) at (1.6,2)   {};
  					 \node (n6) at (1.6,3)   {};
  					 \foreach \from/\to in {n1/n2,n1/n3,n3/n4,n3/n5,n5/n6}
   					 \draw (\from) -- (\to);
					\end{tikzpicture},
					\begin{tikzpicture}
  					 [scale=.15,auto=left,every node/.style={circle,scale=0.2,fill=black!100}]
  					 \node (n1) at (0,0)     {};
  					 \node (n2) at (-0.8,1)  {};
  					 \node (n3) at (0.8,1)   {};
  					 \node (n4) at (0.8,2)   {};
  					 \node (n5) at (0,3)     {};
  					 \node (n6) at (1.6,3)   {};
  					 \foreach \from/\to in {n1/n2,n1/n3,n3/n4,n4/n5,n4/n6}
   					 \draw (\from) -- (\to);
					\end{tikzpicture},
					\begin{tikzpicture}
  					 [scale=.15,auto=left,every node/.style={circle,scale=0.2,fill=black!100}]
  					 \node (n1) at (0,0)      {};
  					 \node (n2) at (-0.8,1)   {};
  					 \node (n3) at (-0.8,2)   {};
  					 \node (n4) at (0.8,1)    {};
  					 \node (n5) at (0,2)      {};
  					 \node (n6) at (1.6,2)    {};
  					 \foreach \from/\to in {n1/n2,n2/n3,n1/n4,n4/n5,n4/n6}
   					 \draw (\from) -- (\to);
					\end{tikzpicture},
					\begin{tikzpicture}
  					 [scale=.15,auto=left,every node/.style={circle,scale=0.2,fill=black!100}]
  					 \node (n1) at (0,0)     {};
  					 \node (n2) at (-0.8,1)  {};
  					 \node (n3) at (-0.8,2)  {};
  					 \node (n4) at (0.8,1)   {};
  					 \node (n5) at (0.8,2)   {};
  					 \node (n6) at (0.8,3)   {};
  					 \foreach \from/\to in {n1/n2,n2/n3,n1/n4,n4/n5,n5/n6}
   					 \draw (\from) -- (\to);
					\end{tikzpicture},
					\begin{tikzpicture}
  					 [scale=.15,auto=left,every node/.style={circle,scale=0.2,fill=black!100}]
  					 \node (n1) at (0,0)     {};
  					 \node (n2) at (-0.8,1)  {};
  					 \node (n3) at (0.8,1)   {};
  					 \node (n4) at (0.8,2)   {};
  					 \node (n5) at (0.8,3)   {};
  					 \node (n6) at (0.8,4)   {};
  					 \foreach \from/\to in {n1/n2,n1/n3,n3/n4,n4/n5,n5/n6}
   					 \draw (\from) -- (\to);
					\end{tikzpicture},
					\begin{tikzpicture}
  					 [scale=.15,auto=left,every node/.style={circle,scale=0.2,fill=black!100}]
  					 \node (n1) at (0,0)     {};
  					 \node (n2) at (0,1)     {};
  					 \node (n3) at (-0.8,2)  {};
  					 \node (n4) at (-0.8,3)  {};
  					 \node (n5) at (0.8,2)   {};
  					 \node (n6) at (0.8,3)   {};
  					 \foreach \from/\to in {n1/n2,n2/n3,n3/n4,n2/n5,n5/n6}
   					 \draw (\from) -- (\to);
					\end{tikzpicture},
					\begin{tikzpicture}
  					 [scale=.15,auto=left,every node/.style={circle,scale=0.2,fill=black!100}]
  					 \node (n1) at (0,0)     {};
  					 \node (n2) at (0,1)     {};
  					 \node (n3) at (-0.8,2)  {};
  					 \node (n4) at (0.8,2)   {};
  					 \node (n5) at (0,3)     {};
  					 \node (n6) at (1.6,3)   {};
  					 \foreach \from/\to in {n1/n2,n2/n3,n2/n4,n4/n5,n4/n6}
   					 \draw (\from) -- (\to);
					\end{tikzpicture},
					\begin{tikzpicture}
  					 [scale=.15,auto=left,every node/.style={circle,scale=0.2,fill=black!100}]
  					 \node (n1) at (0,0)     {};
  					 \node (n2) at (0,1)     {};
  					 \node (n3) at (-0.8,2)  {};
  					 \node (n4) at (0.8,2)   {};
  					 \node (n5) at (0.8,3)   {};
  					 \node (n6) at (0.8,4)   {};
  					 \foreach \from/\to in {n1/n2,n2/n3,n2/n4,n4/n5,n5/n6}
   					 \draw (\from) -- (\to);
					\end{tikzpicture},
					\begin{tikzpicture}
  					 [scale=.15,auto=left,every node/.style={circle,scale=0.2,fill=black!100}]
  					 \node (n1) at (0,0)     {};
  					 \node (n2) at (0,1)     {};
  					 \node (n3) at (-1.2,2)  {};
  					 \node (n4) at (-0.4,2)  {};
  					 \node (n5) at (0.4,2)   {};
  					 \node (n6) at (1.2,2)   {};
  					 \foreach \from/\to in {n1/n2,n2/n3,n2/n4,n2/n5,n2/n6}
   					 \draw (\from) -- (\to);
					\end{tikzpicture},
					\begin{tikzpicture}
  					 [scale=.15,auto=left,every node/.style={circle,scale=0.2,fill=black!100}]
  					 \node (n1) at (0,0)     {};
  					 \node (n2) at (0,1)     {};
  					 \node (n3) at (-0.8,2)  {};
  					 \node (n4) at (0,2)     {};
  					 \node (n5) at (0.8,2)   {};
  					 \node (n6) at (0.8,3)   {};
  					 \foreach \from/\to in {n1/n2,n2/n3,n2/n4,n2/n5,n5/n6}
   					 \draw (\from) -- (\to);
					\end{tikzpicture},
					\begin{tikzpicture}
  					 [scale=.15,auto=left,every node/.style={circle,scale=0.2,fill=black!100}]
  					 \node (n1) at (0,0)      {};
  					 \node (n2) at (0,1)      {};
  					 \node (n3) at (0,2)      {};
  					 \node (n4) at (-0.8,3)   {};
  					 \node (n5) at (0,3)      {};
  					 \node (n6) at (0.8,3)    {};
  					 \foreach \from/\to in {n1/n2,n2/n3,n3/n4,n3/n5,n3/n6}
   					 \draw (\from) -- (\to);
					\end{tikzpicture},
					\begin{tikzpicture}
  					 [scale=.15,auto=left,every node/.style={circle,scale=0.2,fill=black!100}]
  					 \node (n1) at (0,0)     {};
  					 \node (n2) at (0,1)     {};
  					 \node (n3) at (0,2)     {};
  					 \node (n4) at (-0.8,3)  {};
  					 \node (n5) at (0.8,3)   {};
  					 \node (n6) at (0.8,4)   {};
  					 \foreach \from/\to in {n1/n2,n2/n3,n3/n4,n3/n5,n5/n6}
   					 \draw (\from) -- (\to);
					\end{tikzpicture},
					\begin{tikzpicture}
  					 [scale=.15,auto=left,every node/.style={circle,scale=0.2,fill=black!100}]
  					 \node (n1) at (0,0)    {};
  					 \node (n2) at (0,1)    {};
  					 \node (n3) at (0,2)    {};
  					 \node (n4) at (0,3)    {};
  					 \node (n5) at (-0.8,4) {};
  					 \node (n6) at (0.8,4)  {};
  					 \foreach \from/\to in {n1/n2,n2/n3,n3/n4,n4/n5,n4/n6}
   					 \draw (\from) -- (\to);
					\end{tikzpicture},
					\begin{tikzpicture}
  					 [scale=.15,auto=left,every node/.style={circle,scale=0.2,fill=black!100}]
  					 \node (n1) at (0,0)   {};
  					 \node (n2) at (0,1)   {};
  					 \node (n3) at (0,2)   {};
  					 \node (n4) at (0,3)   {};
  					 \node (n5) at (0,4)   {};
  					 \node (n6) at (0,5)   {};
  					 \foreach \from/\to in {n1/n2,n2/n3,n3/n4,n4/n5,n5/n6}
   					 \draw (\from) -- (\to);
					\end{tikzpicture}$$
					respectively. 
					
\subsection{Analysis of implicit Runge-Kutta method}					
\subsubsection{Model test equation}
To analyse implicit R-K method we consider linear first order ODE
\begin{eqnarray}\label{15}
\dot{u}=I\lambda u
\end{eqnarray}
where $I=\sqrt{-1}$ \cite{naj_mon_13}.\\
For the general $R$-stage R-K method we write
\begin{eqnarray}\label{16}
\bm F^{T}=[F_1, F_2, ..., F_R]
\end{eqnarray}
to denote the stages of solutions for $(n+1)$-th time step.\\
Thus for the test equation,
\begin{eqnarray}\label{17}
\bm F=I\lambda(\bm 1 u^n+\Delta t \bm A\bm F)
\end{eqnarray}
with $\bm 1=(1, 1, ..., 1)^T$ being a vector of length $R$.\\
Hence
\begin{eqnarray}\label{18}
\bm F=I\lambda(\bm I_R-I\sigma \bm A)^{-1}\bm 1 u^n
\end{eqnarray}
where $\bm I_R$ is the identity matrix of order $R$ and $\sigma=\lambda\Delta t$.\\
Therefore using Eq. (\ref{2}), the solution at $(n+1)$-th time step is given by
\begin{eqnarray}\label{19}
u^{n+1}=\big(1+I\sigma\bm b^T(\bm I_R-I\sigma \bm A)^{-1}\bm 1\big) u^n.
\end{eqnarray}
The numerical amplification can thus be represented as
\begin{eqnarray}\label{20}
G_N(\sigma)= 1+I\sigma\bm b^T(\bm I_R-I\sigma \bm A)^{-1}\bm 1.
\end{eqnarray}
Comparing with the exact amplification
\begin{eqnarray}\label{21}
G_E(\sigma)= e^{I\sigma},
\end{eqnarray}
we see that for a numerically stable R-K scheme $|G_N(\sigma)|\le 1$, $\forall \sigma$.

Following Simos \cite{sim_93} the amplification (dissipation) and phase (dispersion) errors can be  represented by the quantities
\begin{eqnarray}\label{21.0}
a(\sigma)=1-|G_N(\sigma)|\;\;\;\;\text{and}\;\;\;\;\phi(\sigma)=\sigma-\bm {arg}(G_N(\sigma))
\end{eqnarray}
respectively. If
\begin{eqnarray}\label{21.01}
a(\sigma)=O(\sigma^{p+1})\;\;\;\;\text{and}\;\;\;\;\phi(\sigma)=O(\sigma^{q+1})
\end{eqnarray}
the method is often said to possess dissipative order $p$ and dispersive order $q$ respectively.

Low-dissipation low-dispersion scheme is one where amplification error and phase error remain low even for relative bigger values of $\sigma$. For complete range of values of $\sigma$ values it may be appropriate to define dissipation and dispersion error as
\begin{eqnarray}\label{21.1}
a_{[0, \pi]}=\left[\int_{0}^{\pi}\left|1-|G_N(\sigma)|\right|^2d\sigma\right]^{1/2},
\end{eqnarray}
\begin{eqnarray}\label{21.2}
\phi_{[0, \pi]}=\left[\int_{0}^{\pi}\left|\sigma-\bm {arg}(G_N(\sigma))\right|^2d\sigma\right]^{1/2}.
\end{eqnarray}
\subsubsection{1D convection equation}
Another methodical procedure \cite{bha_sen_sen_13} to study temporal discretization in conjunction to a suitable spatial approximation is via standard one-dimensional convection equation
\begin{eqnarray}\label{37a}
\frac{\partial u}{\partial t}+c\frac{\partial u}{\partial x}=0.
\end{eqnarray}
For a plane wave solution
\begin{eqnarray}\label{37b}
u(x,t)=e^{I(\lambda t-k x)}, \;\;\;k\in \mathbb{R}
\end{eqnarray}
with wave number $k$ and frequency $\lambda$ the Eq. (\ref{37a}) admit linear dispersion relation $
\lambda=ck$. As such solution propagates with phase velocity $v_p=\frac{\lambda}{k}$ and group velocity $v_g=\frac{d\lambda}{dk}$ identically equal to convection velocity $c$. Thus for  convection equation scaled phase velocity and group velocity are both unity. Dispersion relation preserving character of any discretization of Eq. (\ref{37a}) can thus be quantified by computing numerical phase velocity $v_{pN}$ and numerical group velocity $v_{gN}$.

As the evolution of a wave packet containing several wave numbers is rather complicated we express $u(x,t)$ in terms of its Fourier components
\begin{eqnarray}\label{37c}
u(x,t)=\int U(k,t)e^{-I(k x-\lambda t)}dk
\end{eqnarray}
where the integral is performed from $-k_m$ to $k_m$, defined by the Nyquist limit of $k_m=\pi/h$.
Hence at the $j$th node
\begin{eqnarray}\label{37d}
u_j(t)=\int U(k,t)e^{-I(k x_j-\lambda t)}dk
\end{eqnarray}
where $h$ is the uniform grid size. The exact derivative is
\begin{eqnarray}\label{37e}
\frac{\partial u}{\partial x}(x,t)=-\int IkU(k,t)e^{-I(k x-\lambda t)}dk.
\end{eqnarray}
From which we see that
\begin{eqnarray}\label{37f}
\left(\frac{\partial u}{\partial x}\right)_j(t)=-\int IkU(k,t)e^{-I(k x_j-\lambda t)}dk.
\end{eqnarray}
Assuming numerical approximation of spatial derivative is carried out using
\begin{eqnarray}\label{37g}
&&{\bm M_1}\left[\left(\frac{\partial u}{\partial x}\right)_j\right]=\frac{1}{h}{\bm M_2}[u_j]\nonumber\\
&\Rightarrow &\left[\left(\frac{\partial u}{\partial x}\right)_j\right]=\frac{1}{h}{\bm C}[u_j]
\end{eqnarray}
where ${\bm C}={\bm M_1^{-1}}{\bm M_2}=[C_{ij}]$ we see that
\begin{eqnarray}\label{37h}
\left[\left(\frac{\partial u}{\partial x}\right)_j\right]&=&\int\frac{1}{h}\sum_{l=1}^NC_{jl}e^{-Ik(x_l-x_j)}U(k,t)e^{-I(k x_j-\lambda t)}dk\nonumber\\
&=&-\int I[k_{eq}]_jU(k,t)e^{-I(k x_j-\lambda t)}dk
\end{eqnarray}
with
\begin{eqnarray}\label{37i}
[k_{eq}]_j=\frac{I}{h}\sum_{l=1}^NC_{jl}e^{-Ik(x_l-x_j)}.
\end{eqnarray}
Thus semi-discretized form of Eq. (\ref{37a}) can be written as
\begin{eqnarray}\label{37j}
\left(\frac{\partial u}{\partial t}\right)_j=Ic\int [k_{eq}]_jU(k,t)e^{-I(k x_j-\lambda t)}dk.
\end{eqnarray}
Assuming as before that temporal integration lead to numerical amplification factor $G_N$ we see that initial distribution
\begin{eqnarray}\label{37k}
u_j^0=u_j(0)=\int U(k,0)e^{-Ik x_j}dk
\end{eqnarray}
lead to
\begin{eqnarray}\label{37l}
u_j^n&=&\int U(k,0)|G_{Nj}|^ne^{In\beta_j}e^{-Ik x_j}dk\nonumber\\
&=&\int U(k,0)|G_{Nj}|^ne^{-I(k x_j-\frac{\beta_j}{\Delta t}n\Delta t)}dk
\end{eqnarray}
where $G_{Nj}$ is the nodal amplification factor with amplitude $|G_{Nj}|$ and argument $\beta_j$. Eq. (\ref{37l}) on comparison with Eq. (\ref{37c}) reveals numerical circular frequency
\begin{eqnarray}\label{37m}
\lambda_{Nj}=\frac{\beta_j}{\Delta t}.
\end{eqnarray}
Finally comparing Eq. (\ref{15}) with Eq. (\ref{37j})
\begin{eqnarray}\label{37n}
G_{Nj}(\sigma_j)= 1+I\sigma_j\bm b^T(\bm I_R-I\sigma_j \bm A)^{-1}\bm 1.
\end{eqnarray}
where
\begin{eqnarray}\label{37o}
\sigma_j= c[k_{eq}]_j\Delta t=N_ch[k_{eq}]_j.
\end{eqnarray}
$N_c$ being the CFL number. Thus expressions for scaled numerical phase velocity and group velocity at $j$th node can be obtained as
\begin{eqnarray}\label{37p}
\left(\frac{v_{pN}}{c}\right)_j=\frac{\lambda_{Nj}}{kc}=\frac{1}{N_c}\frac{\beta_j}{kh}
\end{eqnarray}
and
\begin{eqnarray}\label{37q}
\left(\frac{v_{gN}}{c}\right)_j=\frac{1}{c}\frac{\partial\lambda_{Nj}}{\partial k}=\frac{1}{N_c}\frac{\partial \beta_j}{\partial (kh)}
\end{eqnarray}
respectively.
\section{Wave analysis of two stage schemes}
In two stage implicit schemes $\bm A=(a_{rs})_{2\times 2}$. Thus numerical amplification factor $G_{N,2}$ can be written as
\begin{eqnarray}\label{22}
G_{N,2}(\sigma)&=&1+I\sigma \left(
                                                                   \begin{array}{c}
                                                                     b_1 \\
                                                                     b_2  \\
                                                                   \end{array}
                                                                 \right)^T\left(
                                   \begin{array}{cc}
                                    1-I\sigma a_{11} &   -I\sigma a_{12}\\
                                     -I\sigma a_{21} &  1-I\sigma a_{22}\\
                                   \end{array}
                                 \right)^{-1}\left(
                                                                   \begin{array}{c}
                                                                     1 \\
                                                                     1  \\
                                                                   \end{array}
                                                                 \right)\nonumber\\
&=&\frac{\mathfrak{Num} G_{N,2}(\sigma)}{\mathfrak{Den} G_{N,2}(\sigma)}
\end{eqnarray}
with
\begin{eqnarray}
\mathfrak{Num} G_{N,2}(\sigma)&=&1+I\sigma(b_1+b_2-a_{11}-a_{22})\nonumber\\
&&+\sigma^2(a_{12}a_{21}-a_{11}a_{22}-b_1(a_{12}-a_{22})+b_2(a_{11}-a_{21})),\label{23}\\
\mathfrak{Den} G_{N,2}(\sigma)&=&1-I\sigma(a_{11}+a_{22})+\sigma^2(a_{12}a_{21}-a_{11}a_{22}).\label{24}
\end{eqnarray}
If the two stage implicit scheme possesses atleast second order accuracy then
\begin{eqnarray}
&&b_1+b_2=1, \label{25}\\
&&b_1(a_{11}+a_{12})+b_2(a_{21}+a_{22})=\frac{1}{2} \label{26}
\end{eqnarray}
and the expressions given in Eqs. (\ref{23}) and (\ref{24}) reduces to
\begin{eqnarray}\label{27}
\mathfrak{Num} G_{N,2}(\sigma)&=&1+I\sigma(1-a_{11}-a_{22})\nonumber\\
&&+\sigma^2(a_{11}+a_{22}-\frac{1}{2}+a_{12}a_{21}-a_{11}a_{22})
\end{eqnarray}
and
\begin{eqnarray}\label{28}
\mathfrak{Den} G_{N,2}(\sigma)&=&1-I\sigma(a_{11}+a_{22})+\sigma^2(a_{12}a_{21}-a_{11}a_{22}).
\end{eqnarray}
For complete minimization of dissipation i.e. infinite dissipative order can be realized if $|G_{N,2}(\sigma)|=1$. The same is achieved by forcing
\begin{eqnarray}\label{29}
a_{11}+a_{22}=\frac{1}{2}.
\end{eqnarray}
Thus in terms of dissipation error we have infinite order of accuracy. This further implies
\begin{eqnarray}\label{29.1}
\bm {arg}(G_{N,2}(\sigma))=2\tan^{-1}\left[{\frac{\sigma/2}{1+\sigma^2Y}}\right],
\end{eqnarray}
with $Y=a_{12}a_{21}-a_{11}a_{22}$.  Eqs. (\ref{25}), (\ref{26}) and (\ref{29}) provide us with a set of three equations in six coefficients $b_1$, $b_2$, $a_{11}$, $a_{12}$, $a_{21}$, $a_{22}$. With more unknowns a copious opportunity arises for minimization of dispersion error. For this purpose a weighted phase error in $L^2$-norm is defined over the entire wave range $\sigma\in[0,\pi]$ as shown below:
\begin{eqnarray}\label{30}
\|PE(Y)\|_{L^2[0,\pi]}=\left[\int_{0}^{\pi}\left|\sigma-2\tan^{-1}\left({\frac{\sigma/2}{1+\sigma^2Y}}\right)\right|^2 \left|e^{-\alpha\sigma^2}\right|^2d\sigma\right]^{1/2}.
\end{eqnarray}
Emphasis then will be on systematic reduction of the above defined error. The idea of weighted minimization of phase error, probably used for the first time in temporal discretization, has been previously used for spatial approxiation by Haras and Ta'asan \cite{har_taa_94}. Authors in \cite{har_taa_94}, made a detailed analysis of the effect of weight function and opined that presence of weight function in $L^2$-norm minimization help better resolve frequencies occurring in solution. In our context it can been seen that the one-dimensional Gaussian weight function $e^{-\alpha\sigma^2}$ with $\alpha\ge0$ give more emphasis to phase error corresponding to smaller values of $\sigma$ at the expense of higher ones. For a fixed $\alpha$ weightage is unity at $\sigma=0$ and decreases with $\sigma$ increasing to $\pi$. This decrease in weightage for higher values of $\sigma$ is accelerated with $\alpha$ increasing, ultimately leading to a situation where weightage is concentrated in a small neighbourhood of $\sigma=0$ for very high $\alpha$ values.

\begin{figure}[!h]
\begin{minipage}[b]{.6\linewidth}\hspace{-1cm}
\centering\psfig{file=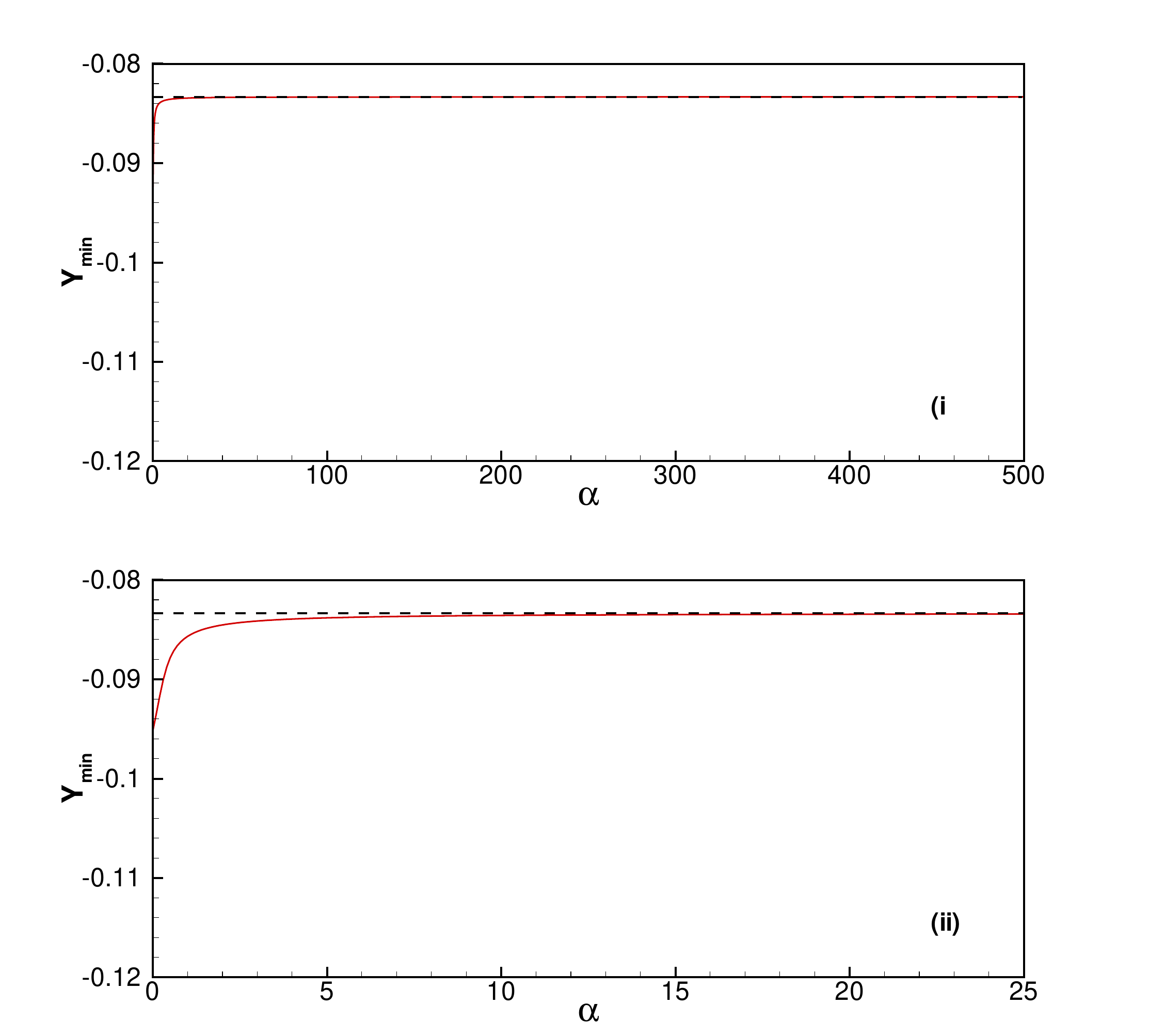,width=0.9\linewidth}\\(a)
\end{minipage}
\begin{minipage}[b]{.6\linewidth}\hspace{-1cm}
\centering\includegraphics[width=75mm]{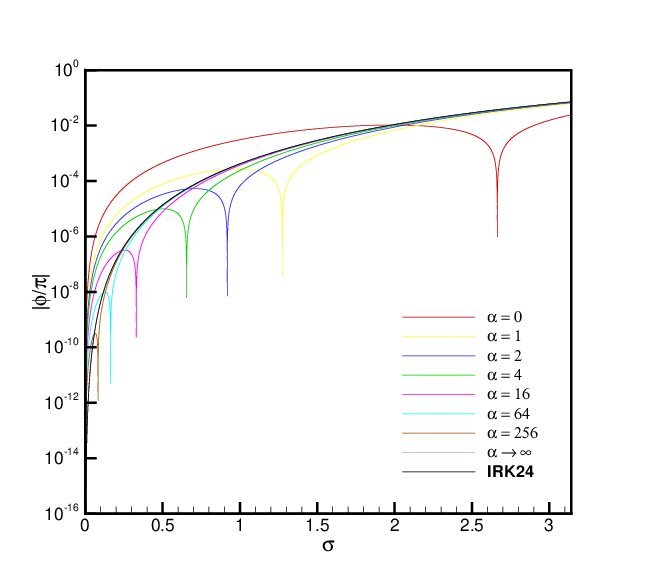}\\(b)
\end{minipage}
\begin{center}
\caption{{\sl Two stage implicit R-K scheme: (a) Variation of points of minima of $L^2$-norm phase error (i) full view, (ii) zoomed view. (b) Comparison of phase error for different values of $\alpha$.} }
    \label{fig:2S_phase_com}
\end{center}
\end{figure}

Eq. (\ref{30}) reveals variation of $\|PE\|_{{L^2[0,\pi]}}$ with parameter $\alpha$ apart from $Y$. With $\alpha$ changing, minimum $L^2$-norm weighted phase error is obtained corresponding to distinct values of $Y$. Plot of point of minima ($Y_{min}$) with $\alpha$ changing is shown in Figure \ref{fig:2S_phase_com}(a)(i).
The plot is found to asymptotically approach $Y=-1/12$. A close-up view in Figure \ref{fig:2S_phase_com}(a)(ii) reveals significant variation of $Y_{min}$ at small values of $\alpha$. To understand its effect on dispersion error we plot normalized magnitude of $\phi=\sigma-2\tan^{-1}\left[{\frac{\sigma/2}{1+\sigma^2Y_{min}}}\right]$ in Figure \ref{fig:2S_phase_com}(b) corresponding to different values of $\alpha$. This figure confirms that with $\alpha$ increasing more emphasis lay on minimizing phase error for smaller angular frequency at the expense of higher ones. In asymptotic case $\alpha\rightarrow\infty$ priority is clearly narrowed down to a small neighbourhood of $\sigma=0$. Mathematically for sufficiently small values of $\sigma$, the asymptotic approach to $Y_{min}=-1/12$, as seen in figure \ref{fig:2S_phase_com}(a) can be explained by using inverse trigonometric expansion in Eq. (\ref{29.1}). Dispersion order quantification then reveals
\begin{eqnarray}\label{30.0}
\phi(\sigma)=\left(Y+\frac{1}{12}\right)\sigma^3+\left(Y^2+\frac{Y}{4}+\frac{1}{80}\right)\sigma^5+O(\sigma^7).
\end{eqnarray}
Hence if one restricts to small $\sigma$ values minimization of Eq. (\ref{30}) should not differ from the minimization of Eq. (\ref{30.0}).

From the figure \ref{fig:2S_phase_com}(b) it is clear that with $\alpha$ values decreasing regions of very low dissipation error can be found even for relatively bigger $\sigma$ values. These regions gradually shifts from higher to lower $\sigma$ values reaching a limiting case for $\alpha\rightarrow\infty$. Thus for $\alpha$ appropriately chosen there lies the possibility to design schemes which will lead to lesser error even at higher $\sigma$ values. Implying in turn better accuracy even at bigger but fittingly chosen CFL number. To understand clearly we work with a set of four $\alpha$ values $\alpha=0, 4, 16$ and $\alpha\rightarrow\infty$ in this work.
\subsection{$\alpha=0$}
Here corresponding minimum value $\|PE\|_{{L^2[0,\pi]}_{\min}}$ is obtained at
\begin{eqnarray}\label{30.1}
Y=Y_{min}=-0.0952154410.
\end{eqnarray}
The set of four equations (\ref{25}), (\ref{26}), (\ref{29}) and (\ref{30.1}) can be solved together. A double infinite family of solution exists all having dissipation error identically equal to zero and an overall dispersion error as given by equation (\ref{21.2}) fixed at $\phi_{[0, \pi]}=4.238151 \times 10^{-2}$. Three such representative schemes are presented in the table \ref{table 1_0}. The first set termed as S2A1 is obtained by putting additional conditions $b_1=b_2$ and $a_{11}=a_{22}$. The second set, S2A2 is reached using $b_1=b_2$ and $a_{12}=2a_{22}$ whereas for the third set S2A3, $b_1=2b_2$ and $a_{12}=2a_{22}$ is assumed.
\begin{table}[h!]
\caption{Second order two stage low-dissipation low-dispersion implicit R-K (S2A) schemes obtained with weight parameter $\alpha=0$.}
\centering
\begin{tabular}{c c c c H}
 \hline
 			&\multicolumn{3}{c}{Schemes}\\ \cline{2-5}
 Parameter 	& S2A1 			& S2A2 			& S2A3          & S2A4\\
 \hline
 $b_1$ 		& 0.5000000000 	&0.5000000000 	&0.6666666667	&0.5000000000      \\
 $b_2$ 		& 0.5000000000 	&0.5000000000 	&0.3333333333	&0.5000000000      \\
 $a_{11}$ 	& 0.2500000000 	&0.2199869148 	&0.3333333333	&0.0952154410      \\
 $a_{12}$ 	& -0.0585699937 &0.5600261703 	&0.3848586017	&-0.0952154410     \\
 $a_{21}$ 	& 0.5585699937 	&-0.0600261703 	&-0.1030505367	&0.5952154410     \\
 $a_{22}$ 	& 0.2500000000 	&0.2800130852 	&0.1666666667	&0.4047845590      \\
 \hline
\end{tabular}
\label{table 1_0}
\end{table}

\subsection{$\alpha=4$}
On minimizing  $\|PE(Y)\|_{{L^2[0,\pi]}}$ given by Eq. (\ref{30}) with $\alpha=4$ we get
\begin{eqnarray}\label{30.2}
Y=Y_{min}=-0.0839362135.
\end{eqnarray}
Again a doubly infinite solution set is obtained on solving Eq. (\ref{30.2}) in conjunction with equations (\ref{25}), (\ref{26}) and (\ref{29}) three of which can be found in table \ref{table 1_4}. Compared to the case of $\alpha=0$, here the set of schemes have overall more dissipation error $\phi_{[0, \pi]}=1.274510\times 10^{-1}$. Additional conditions used to arrive at S2B1, S2B2 and S2B3 are on the same lines as those considered in section 3.1.
\begin{table}[h!]
\caption{Second order two stage low-dissipation low-dispersion implicit R-K (S2B) schemes obtained with weight parameter $\alpha=4$.}
\centering
\begin{tabular}{c c c c H}
 \hline
 			&\multicolumn{3}{c}{Schemes}\\ \cline{2-5}
 Parameter 	& S2B1 			& S2B2 			& S2B3           &S2B4\\
 \hline
 $b_1$ 		& 0.5000000000 	&0.5000000000 	&0.6666666667	&0.5000000000       \\
 $b_2$ 		& 0.5000000000 	&0.5000000000 	&0.3333333333	&0.5000000000       \\
 $a_{11}$ 	& 0.2500000000 	&0.2297892142 	&0.3333333333	&0.0839362135       \\
 $a_{12}$ 	& -0.0397174719 &0.5404215718 	&0.3715278556	&-0.0839362135       \\
 $a_{21}$ 	& 0.5397174719 	&-0.0404215718 	&-0.0763890446	&0.5839362135       \\
 $a_{22}$ 	& 0.2500000000 	&0.2702105718 	&0.1666666667	&0.4160637865         \\
 \hline
\end{tabular}
\label{table 1_4}
\end{table}
\subsection{$\alpha=16$}
With $\alpha$ increased to 16 minimum value of $\|PE(Y)\|_{{L^2[0,\pi]}}$ is obtained at
\begin{eqnarray}\label{30.3}
Y=Y_{min}=-0.0834849563.
\end{eqnarray}
Again we present three out of a system of doubly infinite solution all with similar overall phase error $\phi_{[0, \pi]}=1.319268\times 10^{-1}$ in table \ref{table 1_16}.
\begin{table}[h!]
\caption{Second order two stage low-dissipation low-dispersion implicit R-K (S2C) schemes obtained with weight parameter $\alpha=16$.}
\centering
\begin{tabular}{c c c c H}
 \hline
 			&\multicolumn{3}{c}{Schemes}\\ \cline{2-5}
 Parameter 	& S2C1 			& S2C2 			& S2C3          & S2C4\\
 \hline
 $b_1$ 		& 0.5000000000 	&0.5000000000 	&0.6666666667	&0.5000000000        \\
 $b_2$ 		& 0.5000000000 	&0.5000000000 	&0.3333333333	&0.5000000000        \\
 $a_{11}$ 	& 0.2500000000 	&0.2301921022 	&0.3333333333	&0.0834849563        \\
 $a_{12}$ 	& -0.0389376339 &0.5396157957 	&0.3709764270	&-0.0834849563       \\
 $a_{21}$ 	& 0.5389376339 	&-0.0396157957 	&-0.0752861872	&0.5834849563       \\
 $a_{22}$ 	& 0.2500000000 	&0.2698078978 	&0.1666666667	&0.4165150437        \\
 \hline
\end{tabular}
\label{table 1_16}
\end{table}
\subsection{$\alpha\rightarrow\infty$}
With $\alpha\rightarrow\infty$, $Y_{min}\rightarrow-1/12$. In this context we shall like to point out that a two stage second order implicit Runge-Kutta method attains temporally third order of accuracy if additional conditions
\begin{eqnarray} \label{35}
b_1(a_{11}+a_{12})^2+b_2(a_{21}+a_{22})^2=\frac{1}{3},
\end{eqnarray}
\begin{eqnarray} \label{36}
b_1a_{11}(a_{11}+a_{12})+b_1a_{12}(a_{21}+a_{22})+b_2a_{21}(a_{11}+a_{12})+b_2a_{22}(a_{21}+a_{22})=\frac{1}{6}
\end{eqnarray}
are satisfied. It has been found that the above two equations (\ref{35}) and (\ref{36}) together with the eqs. (\ref{25}), (\ref{26}) and (\ref{29}) implies
\begin{eqnarray}\label{37}
Y=Y_{min}=-\frac{1}{12}.
\end{eqnarray}
Hence all third and higher order two stage $A$ stable schemes enjoys same overall dispersion error of $\phi_{[0, \pi]}=1.334335\times 10^{-1}$ without further scope of minimization. The above choice of $Y_{min}$ also reveals that highest possible fourth order dispersion accuracy. Eqs. (\ref{25}), (\ref{26}), (\ref{29}), (\ref{35}) and (\ref{36}) being five equations in six unknowns admit one parameter family of infinite solutions. We note down three such solutions in the table \ref{table 1_INF}.
\begin{table}[h!]
\caption{Third order two stage low-dissipation low-dispersion implicit R-K (S2D) schemes. They correspond to weight parameter $\alpha\rightarrow\infty$.}
\centering
\begin{tabular}{c c c c}
 \hline
  			&\multicolumn{3}{c}{Schemes}\\ \cline{2-4}
 Parameter 	& S2D1 		   & S2D2 		 & S2D3 \\
 \hline
 $b_1$ 		&0.5000000000  & 0.8367053706 &0.6666666667 \\
 $b_2$ 		&0.5000000000  & 0.1632946294 &0.3333333333 \\
 $a_{11}$ 	&0.2500000000  & 0.4183526852 &0.3333333333 \\
 $a_{12}$ 	&-0.0386751346 & 0.2091763426 &0.3707908119 \\
 $a_{21}$ 	&0.5386751346  &-0.2350933158 &-0.0749149571\\
 $a_{22}$ 	&0.2500000000  & 0.0816473148 &0.1666666667 \\
 \hline
\end{tabular}
\label{table 1_INF}
\end{table}
S2D1 require additional condition $b_1=b_2$. S2D2 and S3D3 are arrived with $a_{11}=2a_{12}$ and $a_{11}=2a_{22}$ respectively. Significantly the first set of values (S2D1) in Table \ref{table 1_INF} dovetails Gauss-Legendre two stage fourth order (IRK24) method which is further discussed in the next subsection.
\subsection{Gauss-Legendre two stage method (IRK24)}
This scheme can be essentially found by taking zeros of Legendre polynomials of degree two viz. $6x^2-6x+1$ as $c_1$ and $c_2$ and striving for highest possible accuracy \cite{but_08}. This method possesses optimal fourth order of accuracy at two stage with Butcher tableau representation
\begin{eqnarray}\nonumber
\begin{tabular}{c|c}
\begin{tabular}{c}
  $\frac{1}{2}-\frac{\sqrt{3}}{6}$ \\
  \\
  $\frac{1}{2}+\frac{\sqrt{3}}{6}$ \\
\end{tabular}
& \begin{tabular}{c c}
            $\frac{1}{4}$ & $\frac{1}{4}-\frac{\sqrt{3}}{6}$  \\
            \\
            $\frac{1}{4}+\frac{\sqrt{3}}{6}$ & $\frac{1}{4}$  \\
          \end{tabular}
\\
  \hline
       & \begin{tabular}{c c}
           $\frac{1}{2}\;\;\;\;\;\;\;\;\;\;$ & $\frac{1}{2}$  \\
         \end{tabular}
\end{tabular}.
\end{eqnarray}
From the previous section it is clear that IRK24 automatically possesses highest possible dispersion accuracy of order four in addition to little dissipation error. That the dispersion characteristics of asymptotic case, $\alpha\rightarrow\infty$, correspond to Gauss-Legendre two stage method can be clearly observed in figure \ref{fig:2S_phase_com}(b) where the respective plots completely overlap each other. Hence if the choice is restricted to a small neighbourhood of $\sigma=0$, IRK24 should be the preferred procedure for time integration. It is well known in the literature that higher order of accuracy does not guarantee lesser phase error. As an example RK44 can be highlighted. Our above analysis indicate that fourth order of accuracy of IRK24 may be rooted in the fourth order dispersive accuracy of the scheme along with zero dissipation error.

\subsection{Comparison of numerical characteristics}
\begin{figure}[!h]
\begin{minipage}[b]{.6\linewidth}\hspace{-1cm}
\centering\includegraphics[width=75mm]{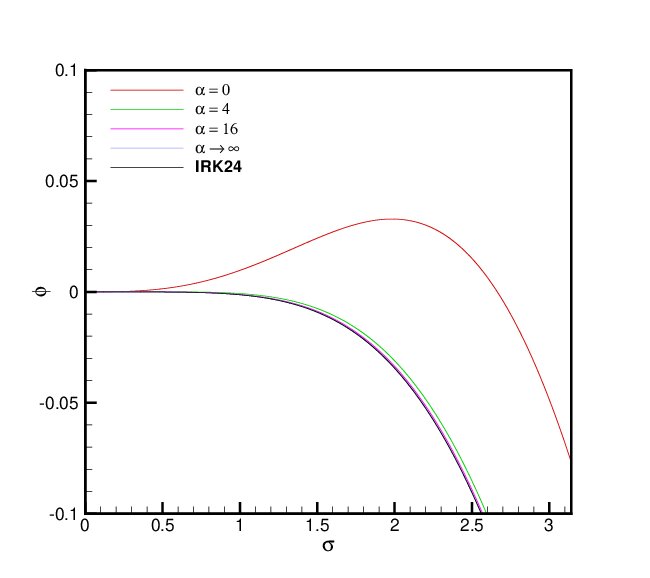}\\(a)
\end{minipage}
\begin{minipage}[b]{.6\linewidth}\hspace{-1cm}
\centering\includegraphics[width=75mm]{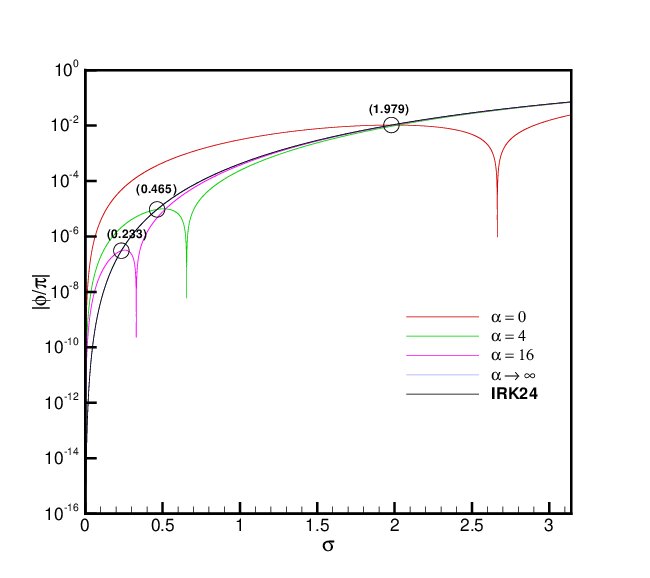}\\(b)
\end{minipage}
\begin{center}
\caption{{\sl (a) Phase difference and (b) dispersion error in logarithmic scale of various
schemes.} }
    \label{fig:2S_phase}
\end{center}
\end{figure}

Numerical phase difference of four different classes of schemes discussed above along with IRK24 is  compared in figure \ref{fig:2S_phase}(a). This figure indicates phase variation of new S2A class of schemes, derived by minimizing $\|PE\|$ with $\alpha=0$ and least overall dispersion error, is completely different from other classes of schemes. In figure  \ref{fig:2S_phase}(b) under logarithmic scale dispersion error of various types of two stage schemes discussed above become prominent. This figure amply demonstrated that there is no difference in the dispersive nature of the S2D schemes obtained using asymptotic value of $\alpha$ and those of IRK24. In the context of low-dissipation low-dispersion schemes IRK24 is indeed a member of a S2D family of schemes albeit with higher order of accuracy. Angular frequency at which phase variation graphs intersect are also shown in this figure. It is seen that S2A class of schemes holds advantage over S2D class of schemes only for $\sigma$ values more than 1.979. From phase error point of view, it can be said that, S2A schemes is better to S2D for waves with more than 3.175 time step per period since $T/\Delta t=2\pi/\sigma$. Signifying its advantage only at sufficiently big temporal step size.
As truncation error grows significantly with bigger step sizes advantages of such scheme may not get reflected in overall error reported. In terms of phase error new S2B and new S2C class of schemes are seen to be better than S2D and hence IRK24 for angular frequencies greater than 0.465 and 0.223 respectively. But overall quantum of advantage of S2C class of schemes vis-a-vis S2D class is quite restricted. Nevertheless there seems to be threshold values of $\sigma$ beyond which new S2B and new S2C class of schemes should report less dispersive solution compared to IRK24. For new S2A class this threshold is quite high.

\begin{figure}[p]
\vspace{-18 mm}
\begin{minipage}[b]{.6\linewidth}
\centering\psfig{file=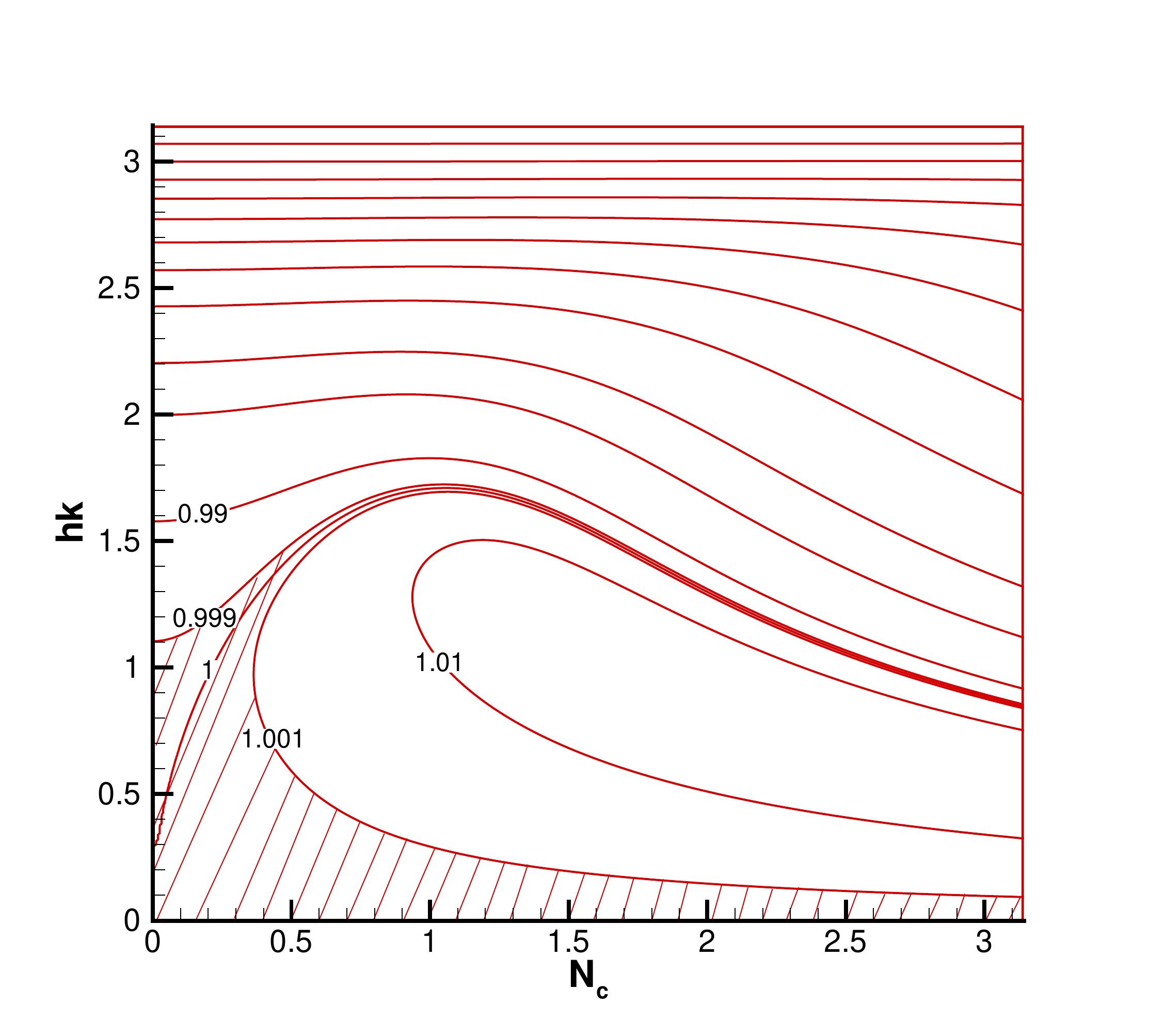,width=0.7\linewidth}
 (a)
\end{minipage}            \hspace{-2.5mm}
\begin{minipage}[b]{.6\linewidth}
\centering\psfig{file=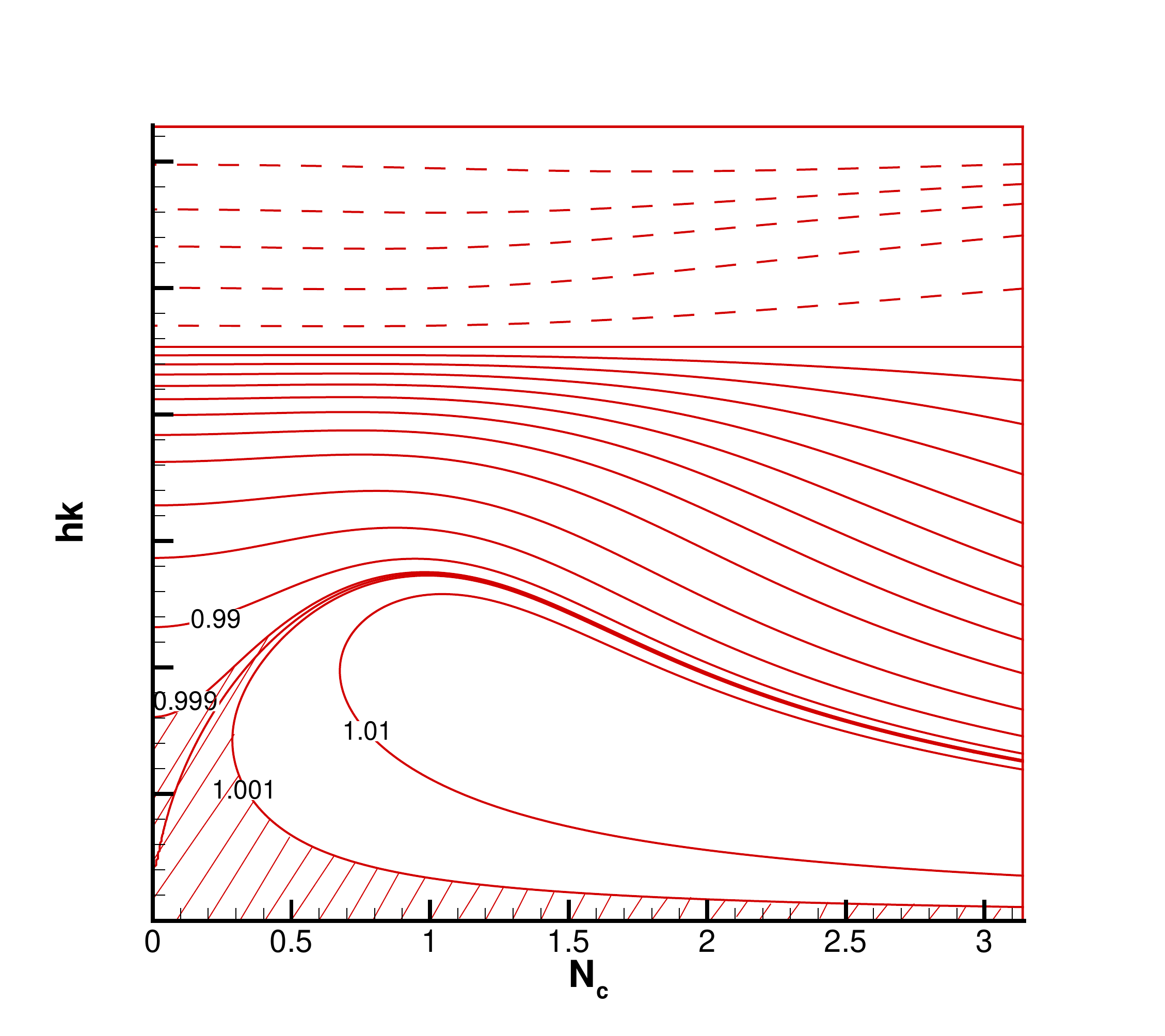,width=0.7\linewidth}
 (b)
\end{minipage}            \hspace{-2.5mm}
\begin{minipage}[b]{.6\linewidth}   \
\centering\psfig{file=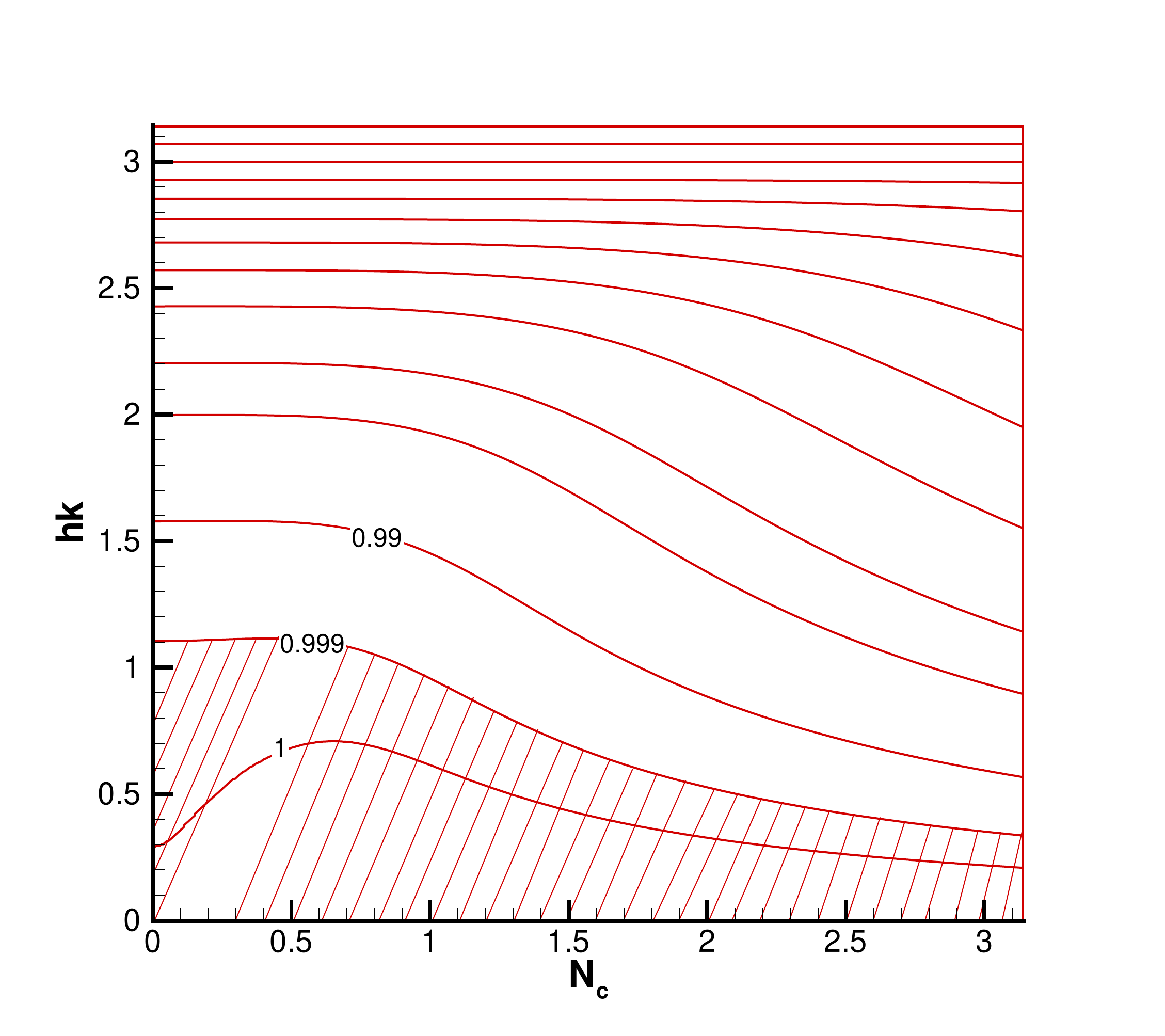,width=0.7\linewidth}
 (c)
\end{minipage}            \hspace{-2.5mm}
\begin{minipage}[b]{.6\linewidth}
\centering\psfig{file=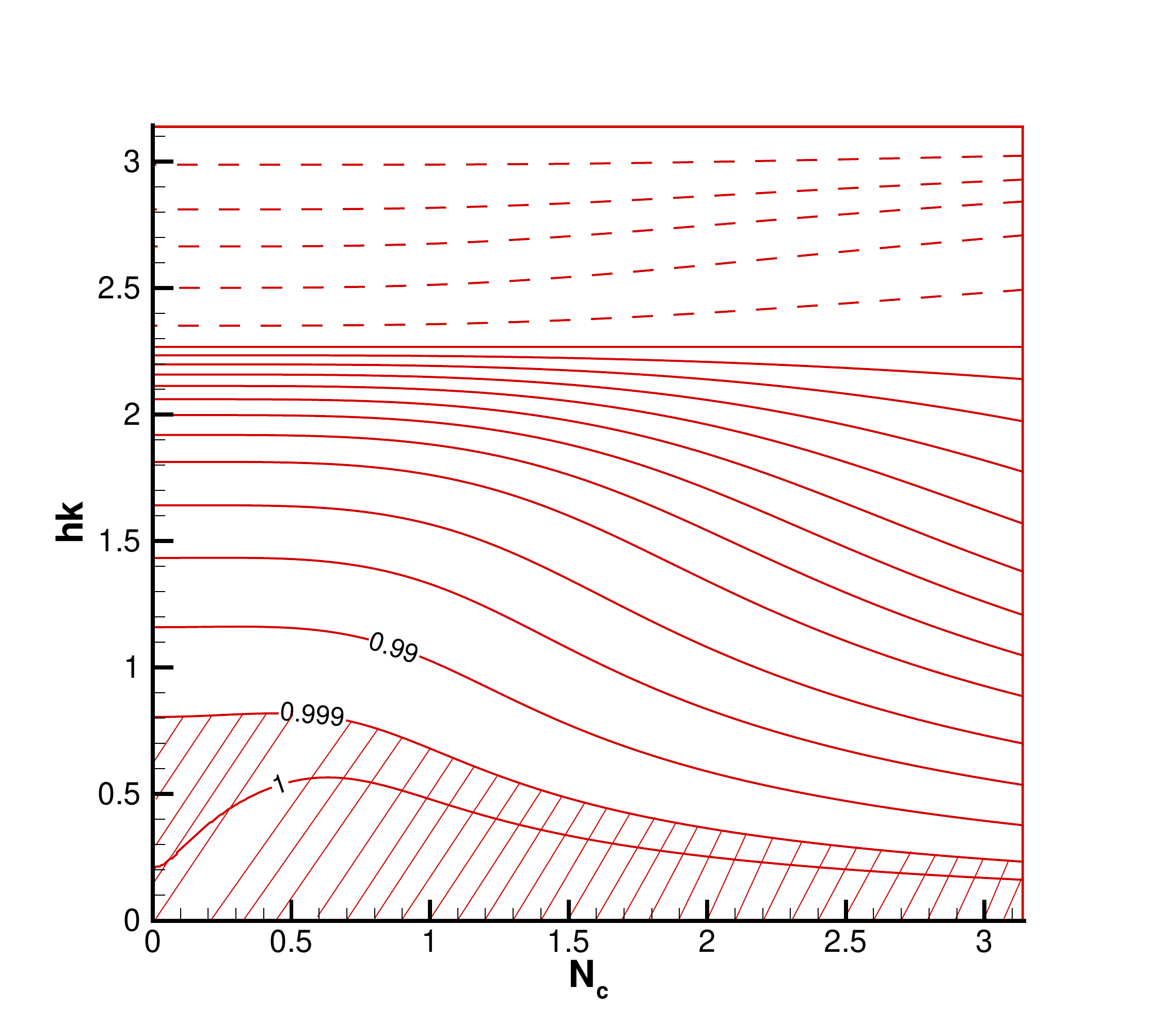,width=0.7\linewidth}
 (d)
\end{minipage}            \hspace{-2.5mm}
\begin{minipage}[b]{.6\linewidth}   \
\centering\psfig{file=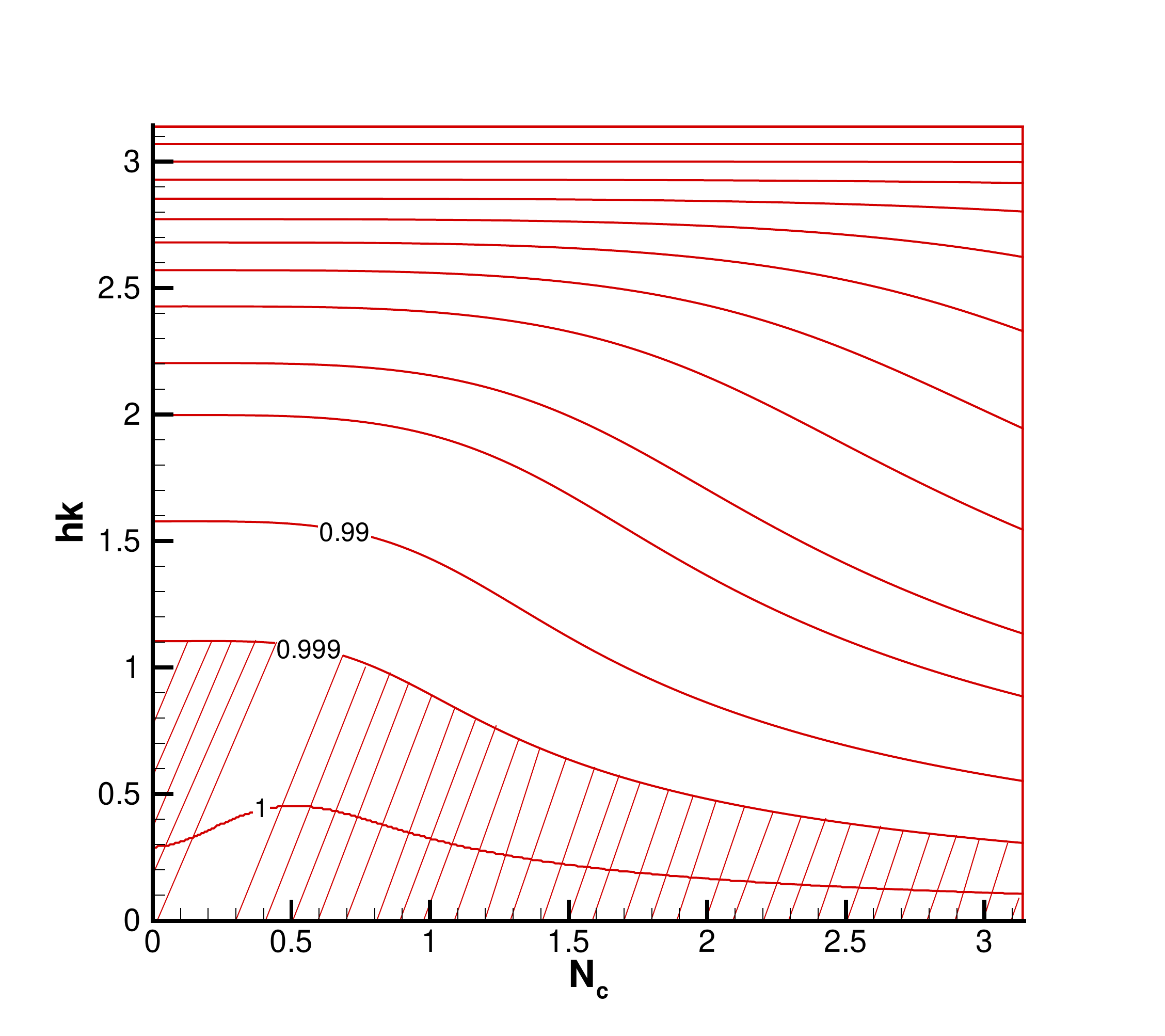,width=0.7\linewidth}
 (e)
\end{minipage}            \hspace{-2.5mm}
\begin{minipage}[b]{.6\linewidth}
\centering\psfig{file=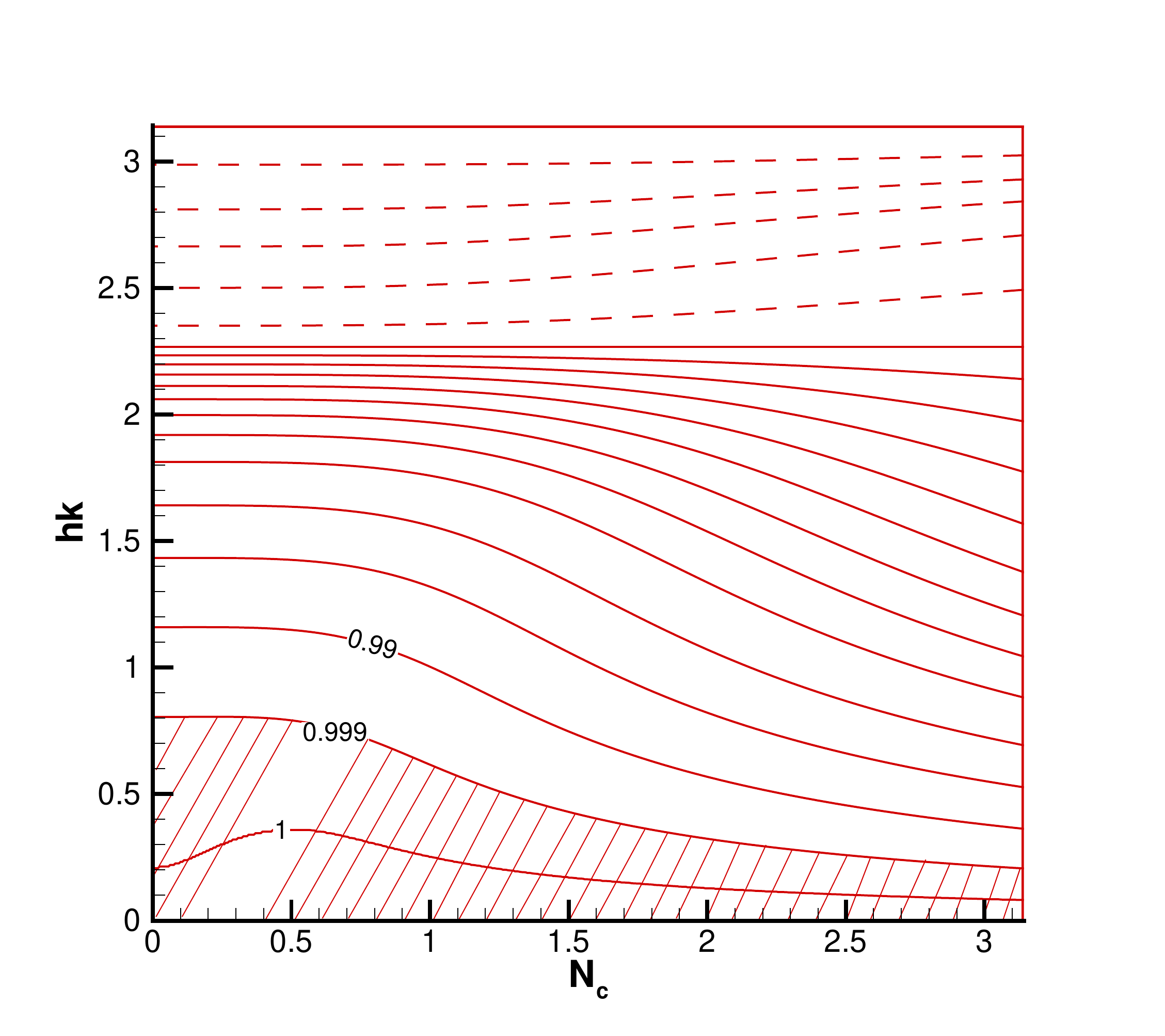,width=0.7\linewidth}
 (f)
\end{minipage}            \hspace{-2.5mm}
\begin{minipage}[b]{.6\linewidth}   \
\centering\psfig{file=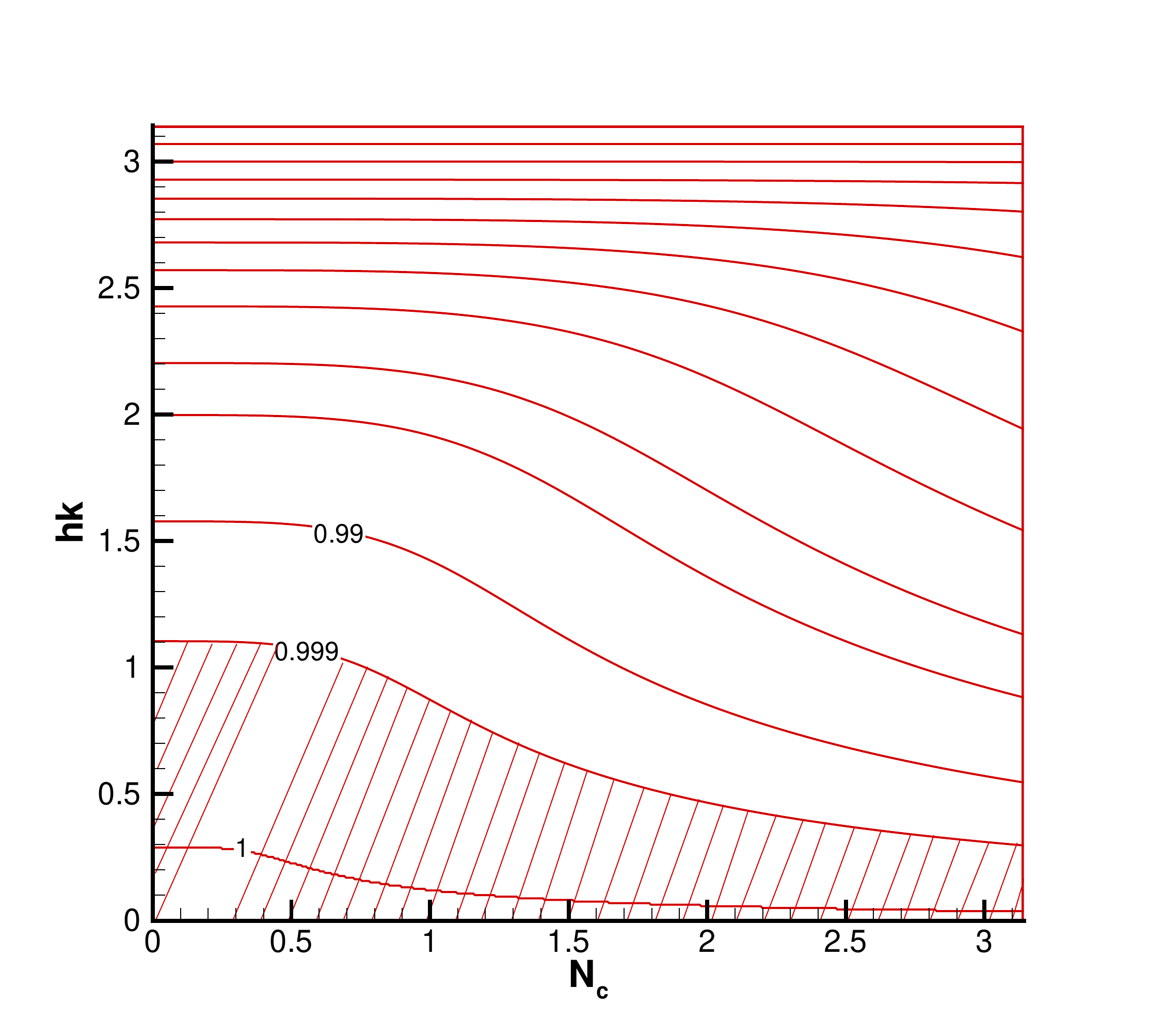,width=0.7\linewidth}
 (g)
\end{minipage}            \hspace{-2.5mm}
\begin{minipage}[b]{.6\linewidth}
\centering\psfig{file=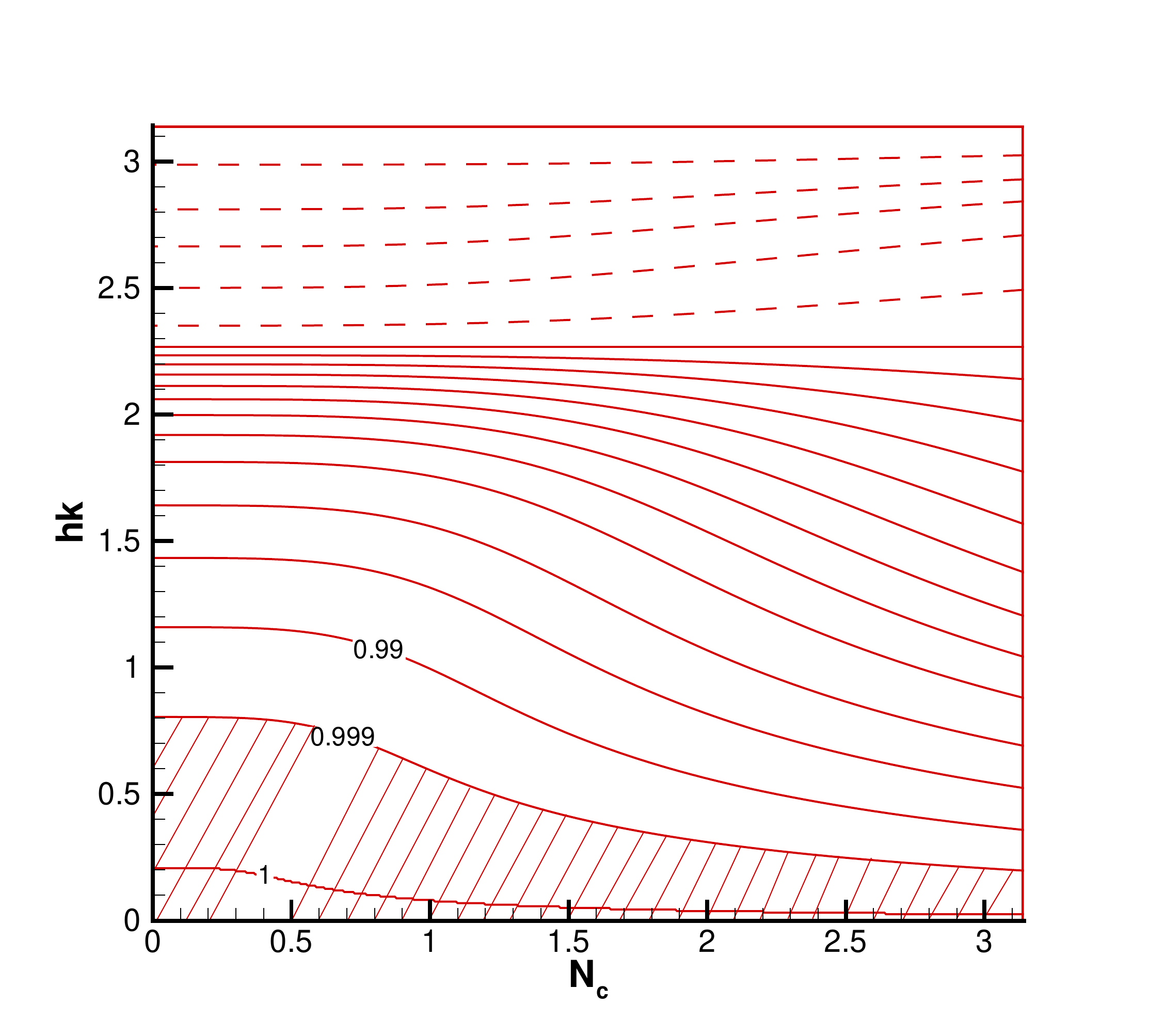,width=0.7\linewidth}
 (h)
\end{minipage}            \hspace{-2.5mm}
\begin{center}
\caption{{\sl Contours of normalized numerical phase velocity (left) and group velocity (right) for indicated schemes plotted in the $(N_c, hk)$ plane at mid-node when Lele scheme is used for spatial discretization: (a)-(b) S2A $(\alpha=0)$, (c)-(d) S2B $(\alpha=4)$, (e)-(f) S2C $(\alpha=16)$, (g)-(h) IRK24, S2D $(\alpha\rightarrow\infty)$.}}
\label{fig:2S_impl}
\end{center}
\end{figure}

\begin{figure}[p]
\vspace{-18 mm}
\begin{minipage}[b]{.6\linewidth}
\centering\psfig{file=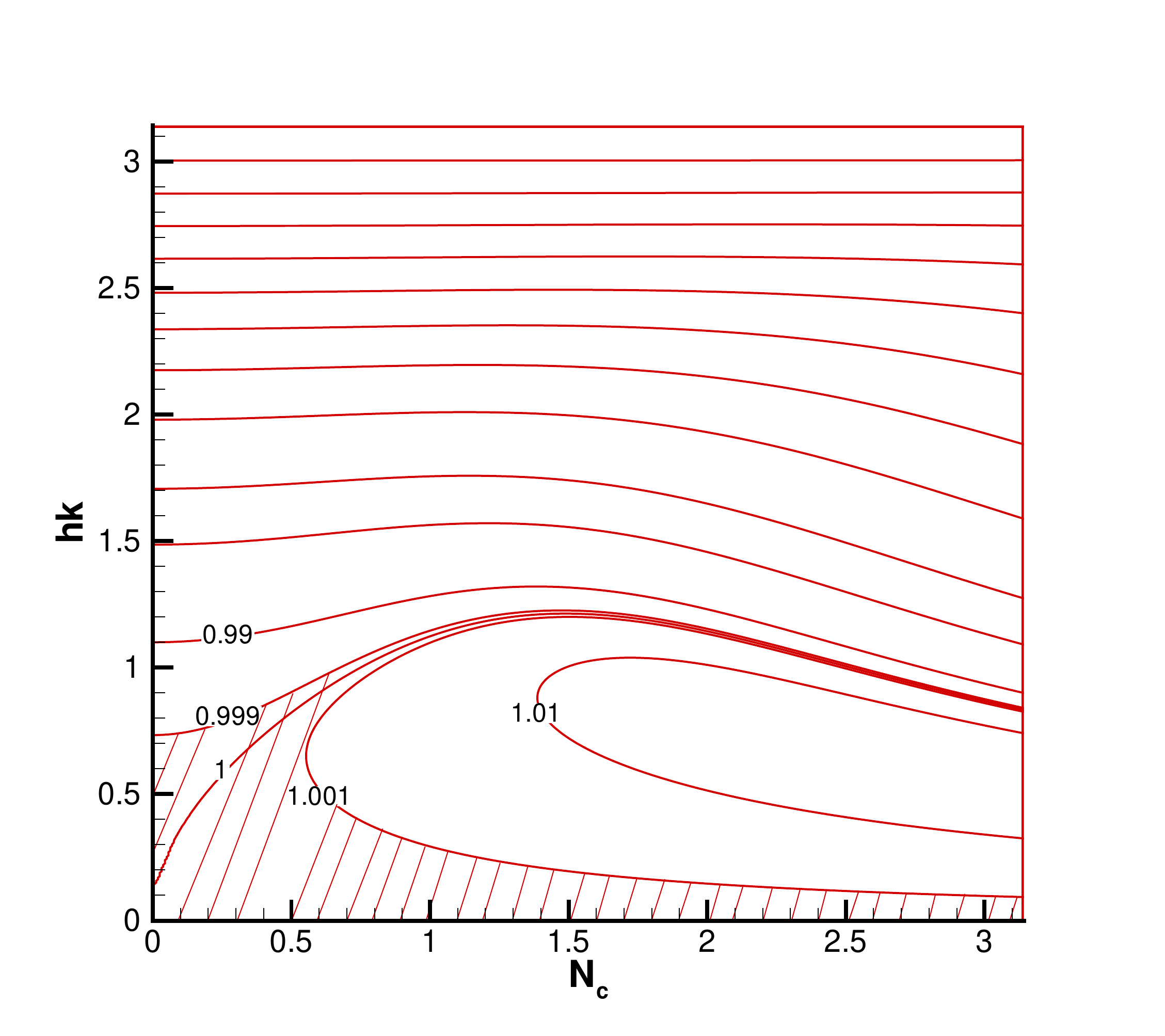,width=0.7\linewidth}
 (a)
\end{minipage}            \hspace{-2.5mm}
\begin{minipage}[b]{.6\linewidth}
\centering\psfig{file=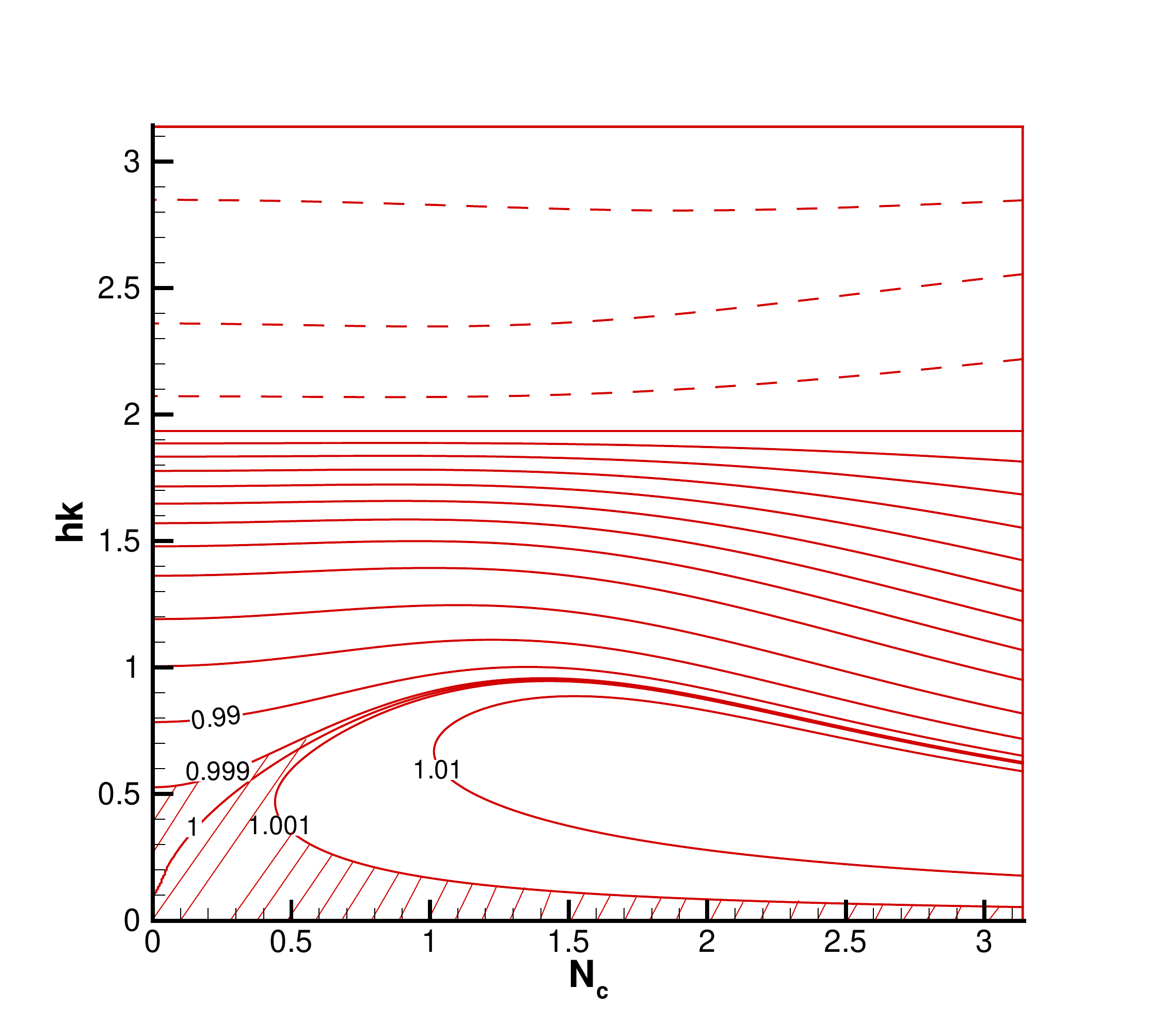,width=0.7\linewidth}
 (b)
\end{minipage}            \hspace{-2.5mm}
\begin{minipage}[b]{.6\linewidth}   \
\centering\psfig{file=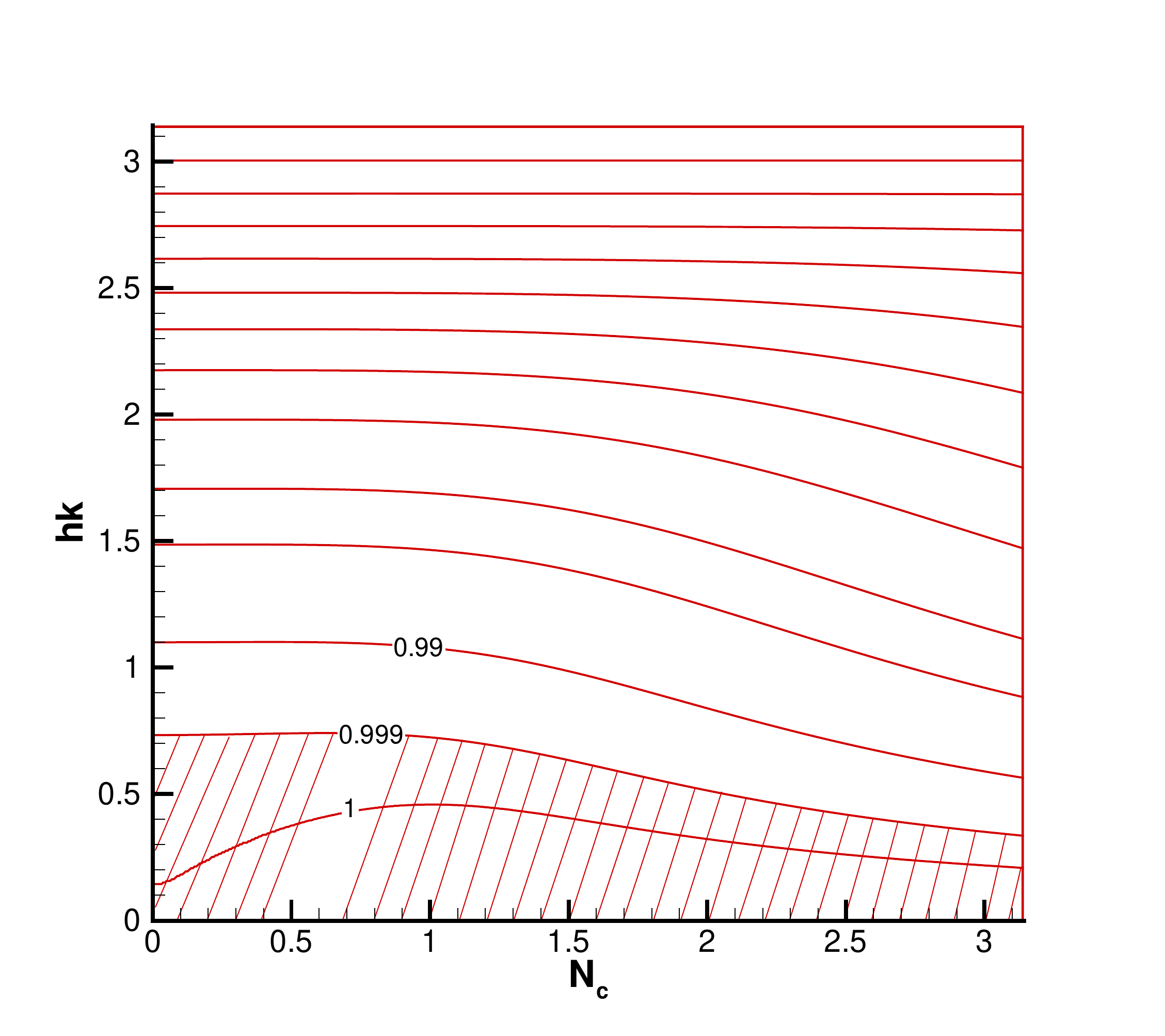,width=0.7\linewidth}
 (c)
\end{minipage}            \hspace{-2.5mm}
\begin{minipage}[b]{.6\linewidth}
\centering\psfig{file=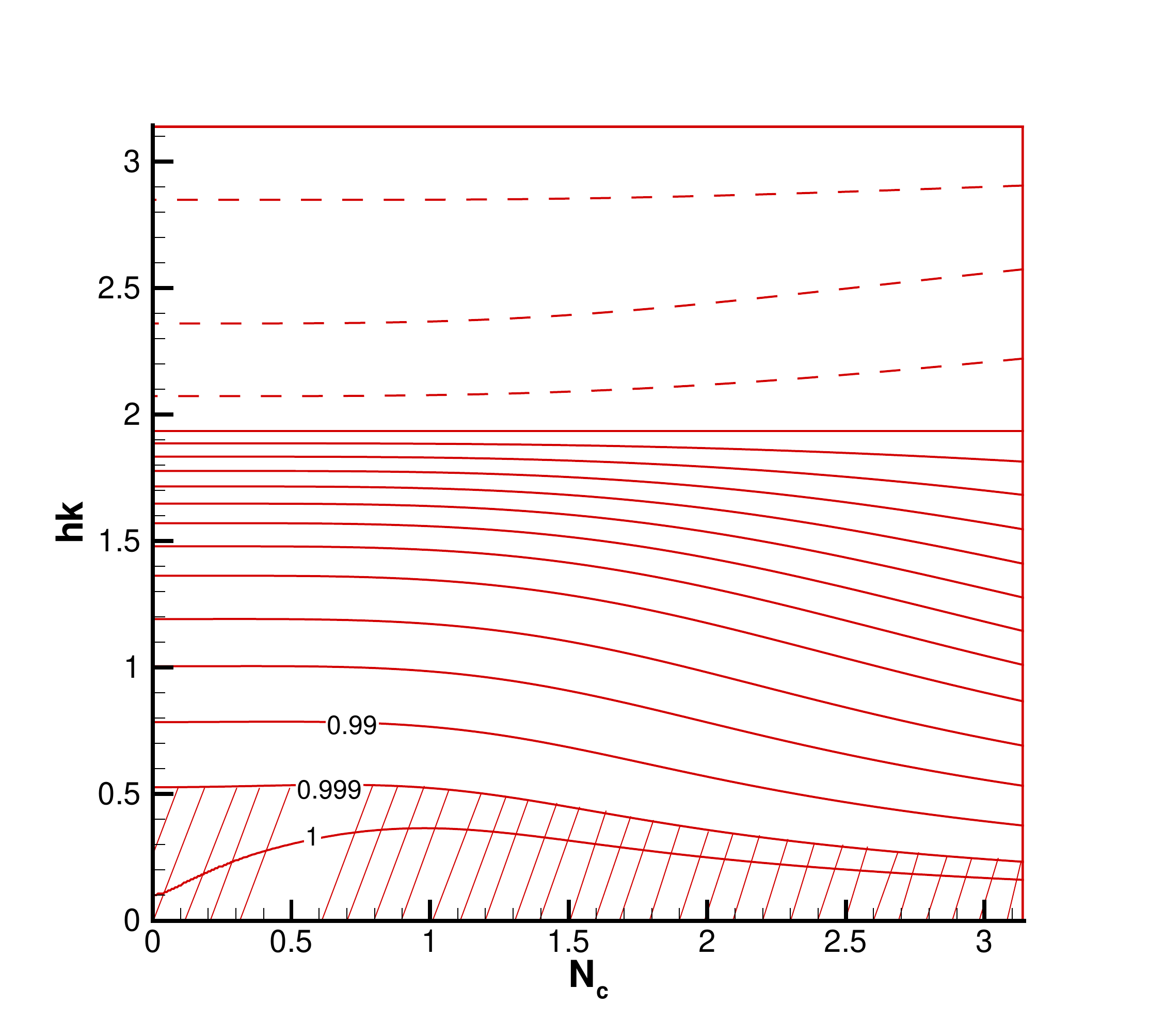,width=0.7\linewidth}
 (d)
\end{minipage}            \hspace{-2.5mm}
\begin{minipage}[b]{.6\linewidth}   \
\centering\psfig{file=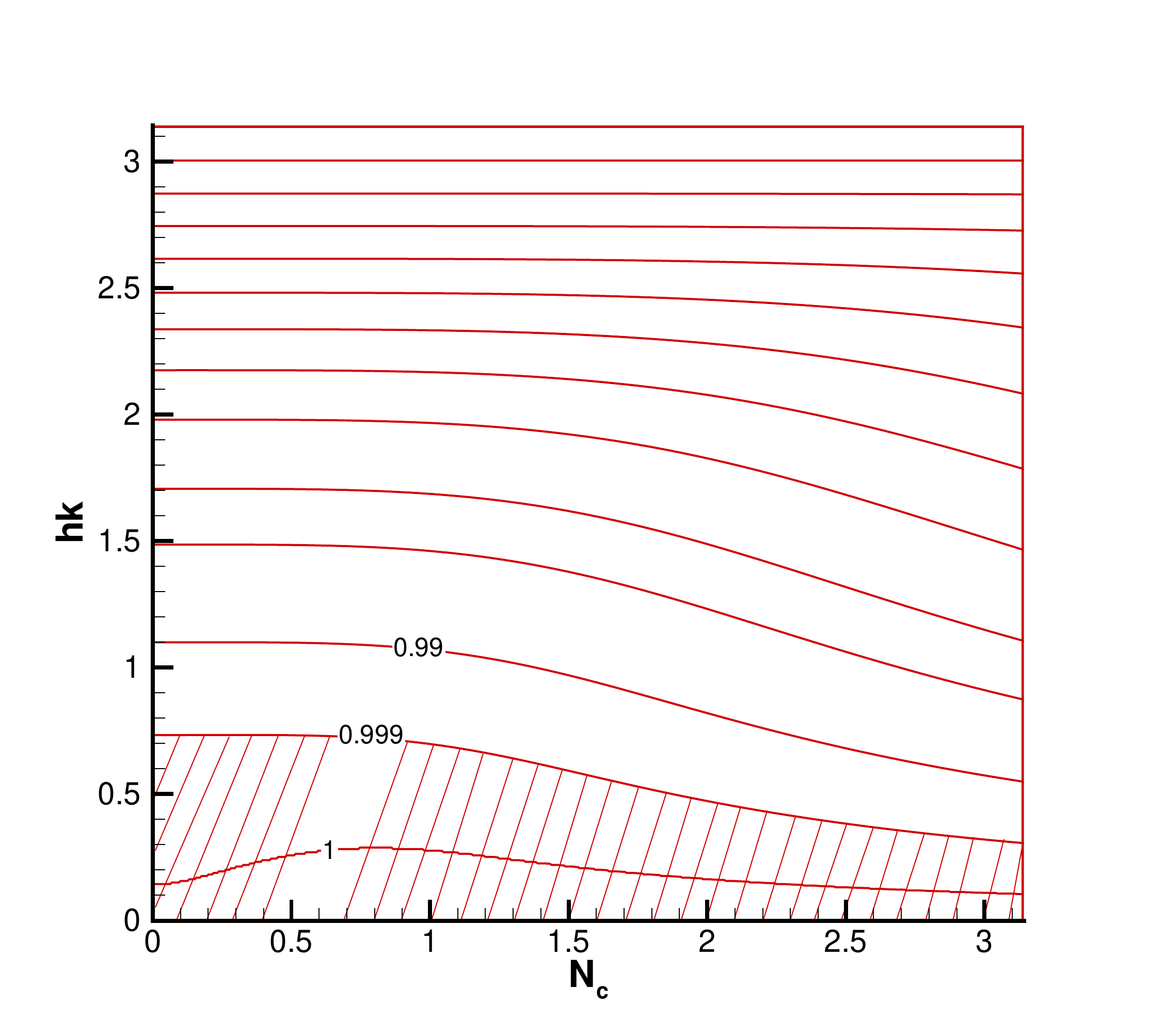,width=0.7\linewidth}
 (e)
\end{minipage}            \hspace{-2.5mm}
\begin{minipage}[b]{.6\linewidth}
\centering\psfig{file=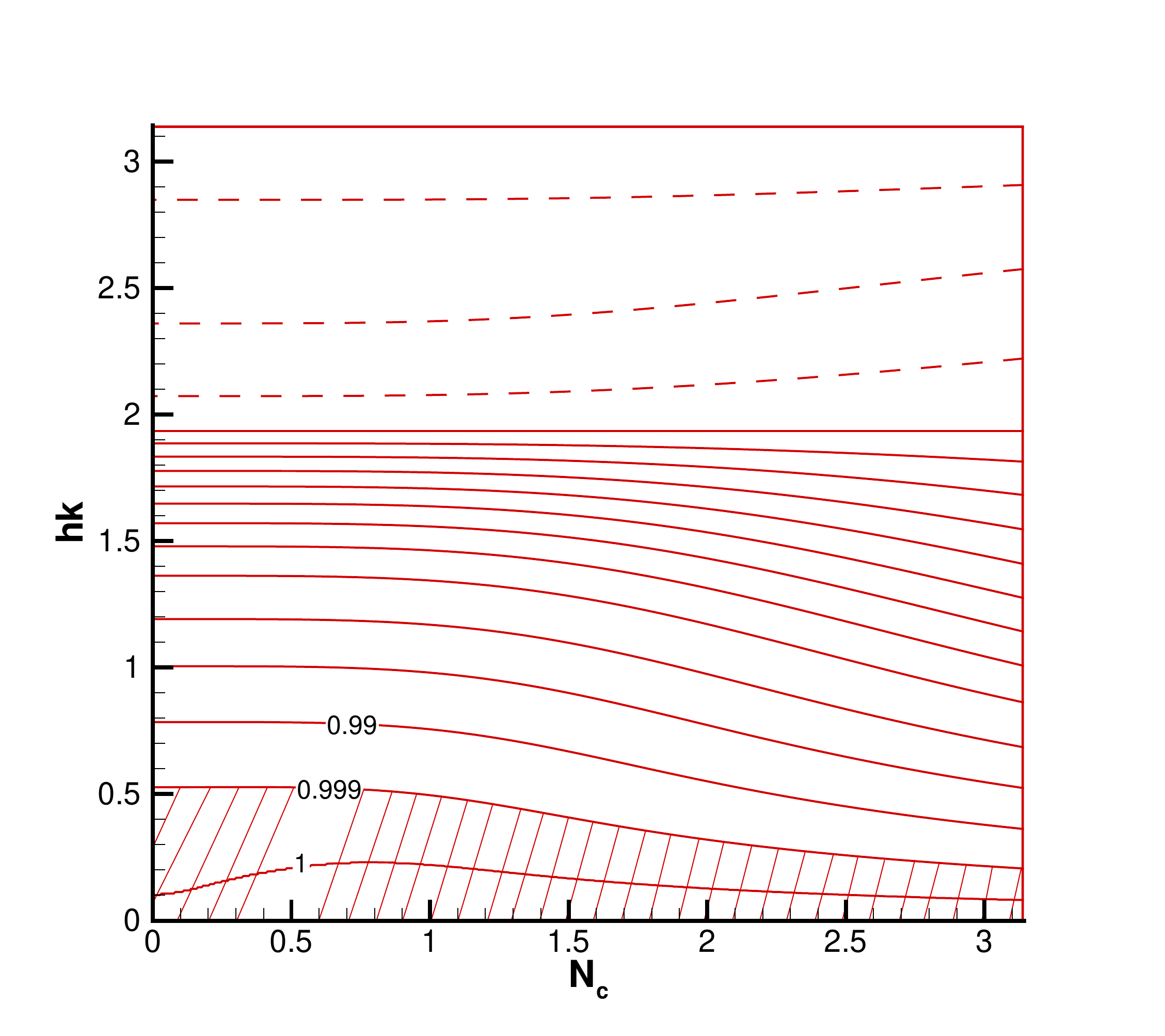,width=0.7\linewidth}
 (f)
\end{minipage}            \hspace{-2.5mm}
\begin{minipage}[b]{.6\linewidth}   \
\centering\psfig{file=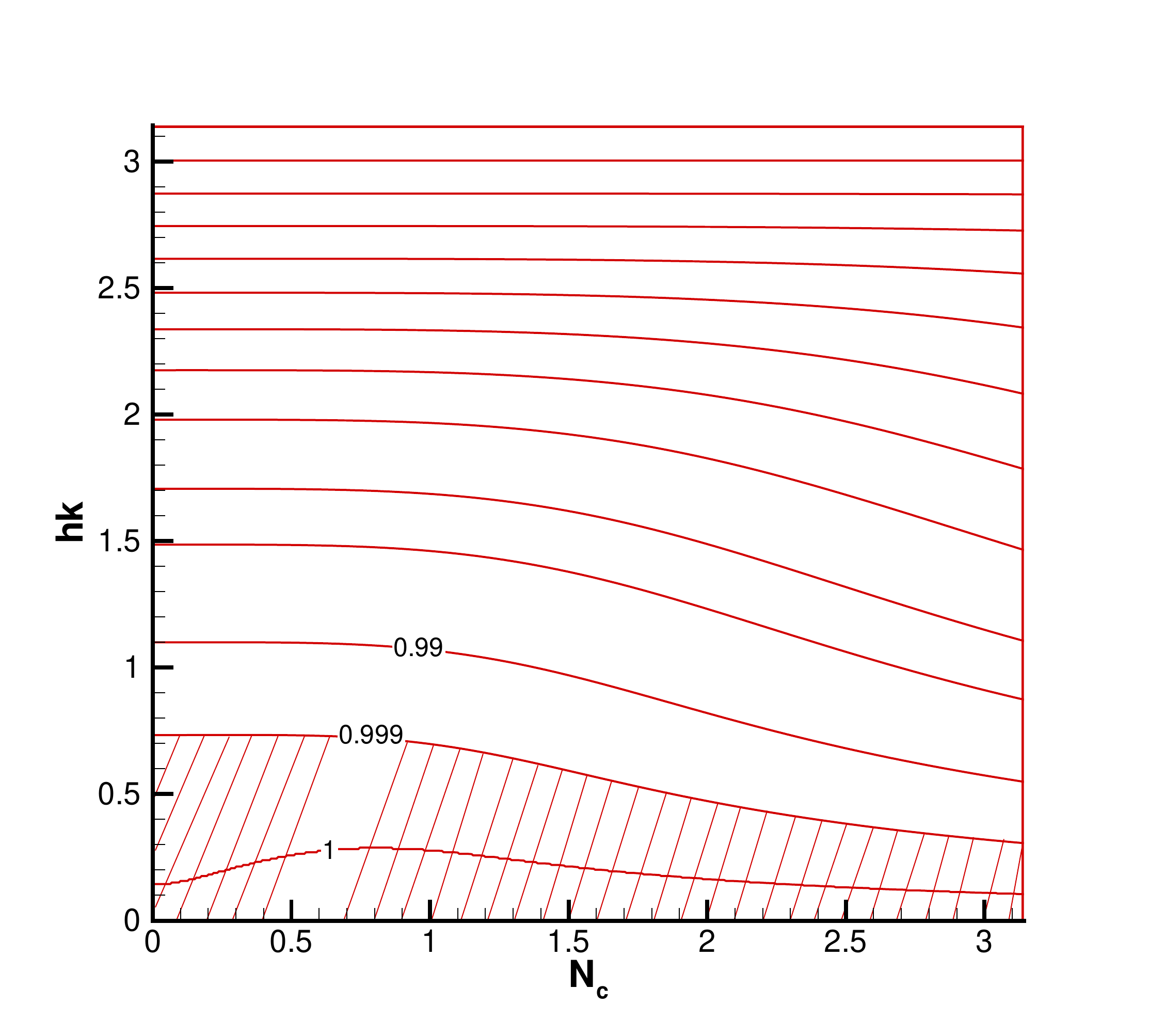,width=0.7\linewidth}
 (g)
\end{minipage}            \hspace{-2.5mm}
\begin{minipage}[b]{.6\linewidth}
\centering\psfig{file=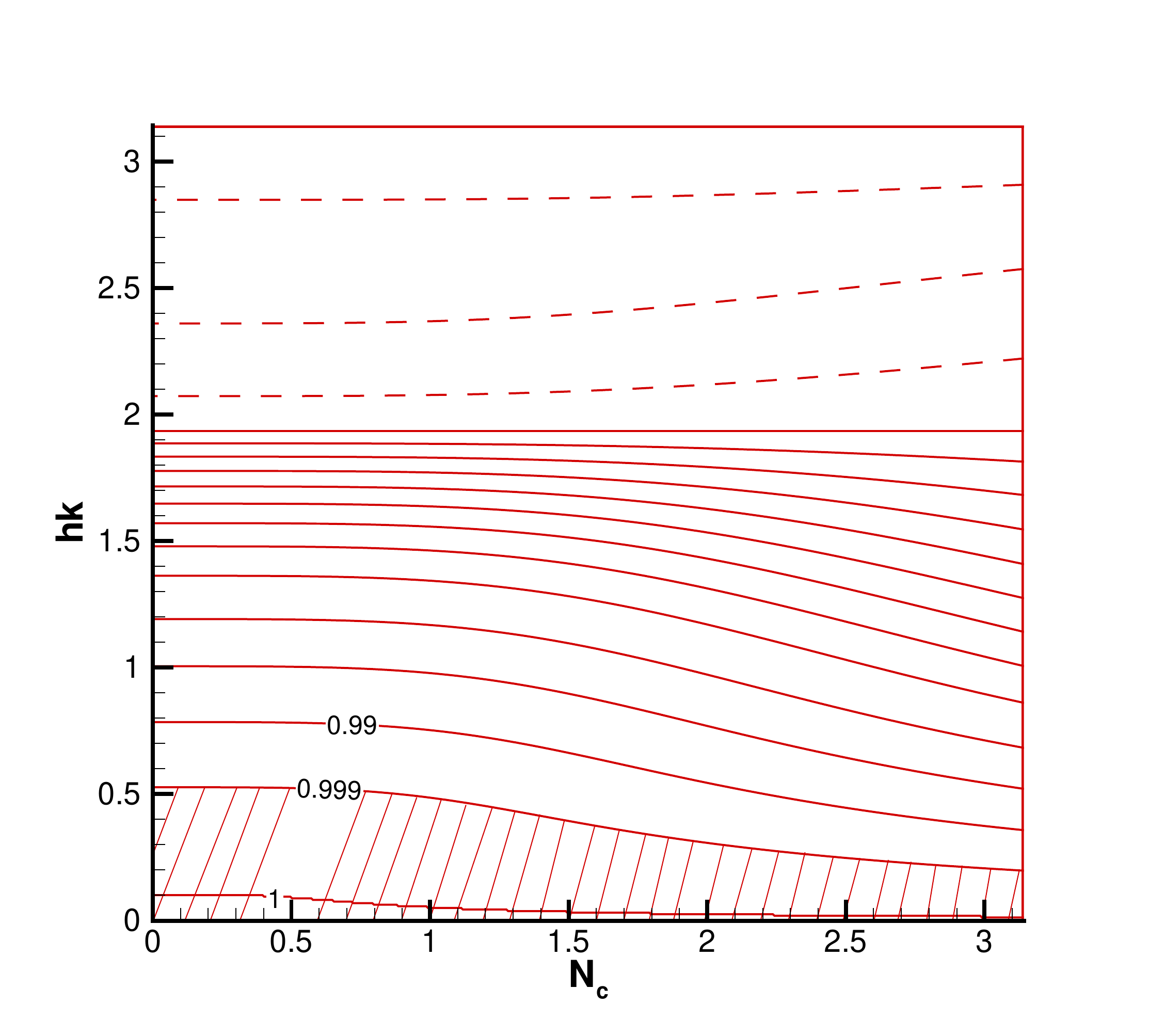,width=0.7\linewidth}
 (h)
\end{minipage}            \hspace{-2.5mm}
\begin{center}
\caption{{\sl Contours of normalized numerical phase velocity (left) and group velocity (right) for indicated schemes plotted in the $(N_c, hk)$ plane at mid-node when explicit CD6 scheme is used for spatial discretization: (a)-(b) S2A $(\alpha=0)$, (c)-(d) S2B $(\alpha=4)$, (e)-(f) S2C $(\alpha=16)$, (g)-(h) IRK24, S2D $(\alpha\rightarrow\infty)$.}}
\label{fig:2S_expl}
\end{center}
\end{figure}

To probe further we compare the dispersion property of the numerical schemes discussed here in terms of $v_{pN}/c$ and $v_{gN}/c$. It is well documented in the literature \cite{bha_sen_sen_13, sen_dip_sau_07, raj_sen_dut_10} that at all nodal points both $\left[1-\left(\frac{v_{pN}}{c}\right)_j\right]$ and $\left[1-\left(\frac{v_{gN}}{c}\right)_j\right]$ should be kept close to zero for dispersion error free computation. Here nodal spectral properties are evaluated using 501 grid point in conjunction with implicit Lele scheme \cite{lel_92} and explicit central CD6 scheme both of which are spatially sixth order accurate. In figure \ref{fig:2S_impl}, $(v_{pN}/c)_j$ (left) and $(v_{gN}/c)_j$ (right) contours are plotted for the central node for various two stage schemes discussed in this work when the implicit scheme of Lele \cite{lel_92} is used to discretize space derivative. The region of $\pm0.1\%$ tolerance in the scaled values of numerical phase and group velocity have been hatched to highlight regions of minimum dispersion error. For S2A schemes this region is significantly less than S2B, S2C and S2D invalidating possibility of long term computation with high $N_c$ values. Dispersion character of S2B, S2C and S2D are found to be similar and rules out $q$-wave formation for $hk<2.27$. Here it must be pointed out that contours of $(v_{pN}/c)_j$ and $(v_{gN}/c)_j$ are exactly same for IRK24 and other S2D set of schemes and are not shown separately. In combination with explicit CD6 spatial descritization the region of tolerance decreases as seen in figure \ref{fig:2S_expl}. But the overall pattern remain same pointing towards robustness of the schemes developed. $q$-wave formation in this case is ruled out for $hk<1.94$.

\section{Wave analysis of three stage schemes}
\begin{figure}[!h]
\begin{minipage}[b]{.6\linewidth}\hspace{-1cm}
\centering\psfig{file=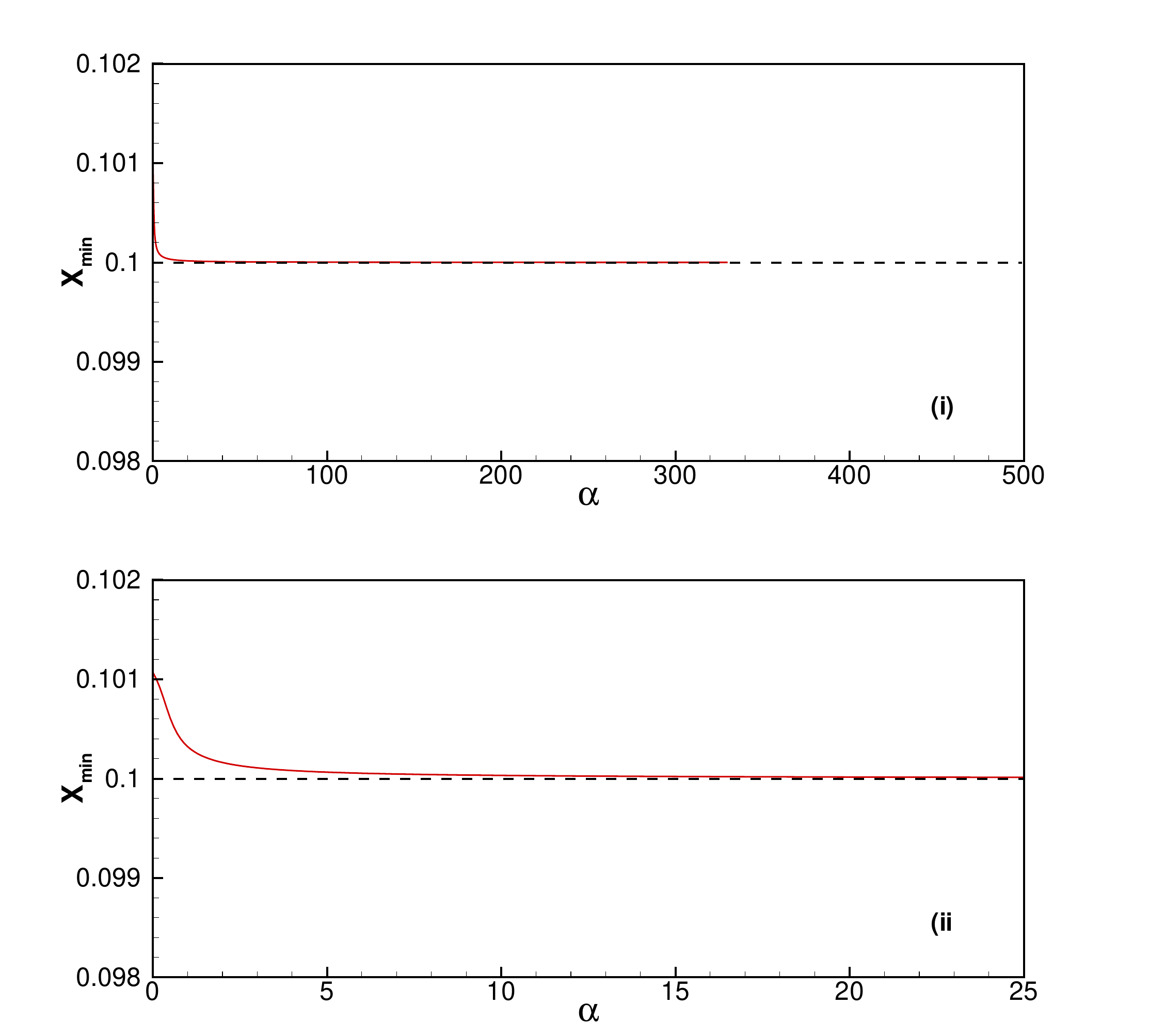,width=0.9\linewidth}\\(a)
\end{minipage}
\begin{minipage}[b]{.6\linewidth}\hspace{-1cm}
\centering\includegraphics[width=75mm]{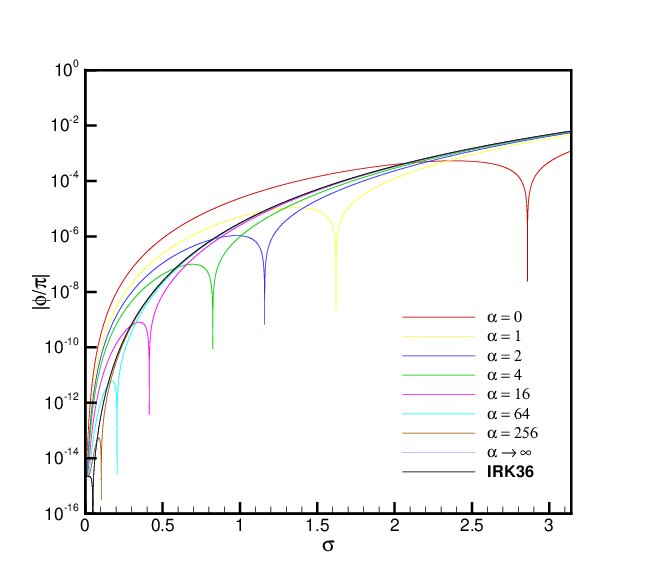}\\(b)
\end{minipage}
\begin{center}
\caption{{\sl Three stage implicit R-K scheme: (a) Variation of points of minima of $L^2$-norm phase error (i) full view, (ii) zoomed view. (b) Comparison of phase error for different values of $\alpha$.} }
    \label{fig:3S_phase_com}
\end{center}
\end{figure}

For three stage methods $\bm A=(a_{rs})_{3\times 3}$, $\bm b=(b_r)_3$. Hence we are required to find as many as twelve free parameters $\{b_r, a_{rs}: 1\le r, s \le3\}$. Numerical amplification factor $G_{N,3}(\sigma)$ for a three stage implicit scheme with not less than second order accuracy is
\begin{eqnarray}\label{38}
G_{N,3}(\sigma)= \frac{\mathfrak{Num} G_{N,3}(\sigma)}{\mathfrak{Den} G_{N,3}(\sigma)}
\end{eqnarray}
where
\begin{eqnarray}\label{39}
\begin{split}
\mathfrak{Num} G_{N,3}(\sigma)=&[\sigma(1-a_{11}-a_{22}-a_{33})+\sigma^3(a_{11}a_{22}a_{33}
+a_{12}a_{23}a_{31}+a_{21}a_{32}a_{13}\\
&-a_{11}a_{23}a_{32}-a_{22}a_{31}a_{13}-a_{33}a_{12}a_{21}\\
&-(a_{22}a_{33}+a_{12}a_{23}+a_{13}a_{32}-a_{23}a_{32}-a_{12}a_{33}-a_{13}a_{22})b_1 \\
&-(a_{33}a_{11}+a_{23}a_{31}+a_{21}a_{13}-a_{31}a_{13}-a_{23}a_{11}-a_{21}a_{33})b_2\\
&-(a_{11}a_{22}+a_{31}a_{12}+a_{32}a_{21}-a_{12}a_{21}-a_{31}a_{22}-a_{32}a_{11})b_3)]\\
&+I[-1+\sigma^2(1/2-a_{11}-a_{22}-a_{33}\\
&+a_{11}a_{22}+a_{22}a_{33}+a_{33}a_{11}-a_{12}a_{21}-a_{23}a_{32}-a_{31}a_{13})],
\end{split}
\end{eqnarray}
\begin{eqnarray}\label{40}
\begin{split}
\mathfrak{Den} G_{N,3}(\sigma)=&-[\sigma(a_{11}+a_{22}+a_{33})+\sigma^3(-a_{11}a_{22}a_{33}-a_{12}a_{23}a_{31}-a_{21}a_{32}a_{13}\\
&+a_{11}a_{23}a_{32}+a_{22}a_{31}a_{13}+a_{33}a_{12}a_{21})]
\\&+I[-1+\sigma^2(a_{11}a_{22}+a_{22}a_{33}+a_{33}a_{11}-a_{12}a_{21}-a_{23}a_{32}-a_{31}a_{13})].
\end{split}
\end{eqnarray}
Resulting scheme can be endowed with minimal dissipation error if we can compute the coefficients in a fashion such that $|G_{N,3}(\sigma)|=1$. Careful inspection of eqs. (\ref{39}) and (\ref{40}) reveals that it is imperative to impose 
\begin{eqnarray}\label{41}
a_{11}+a_{22}+a_{33}=\frac{1}{2}.
\end{eqnarray}
for dissipation error free computation to be made possible. This leaves us with ample opportunity to aim for fourth order accuracy which will require satisfaction of additional five conditions given in Eqs. (\ref{9})-(\ref{13}), Eq. (\ref{14}) being dependent on others. Analysis and further simplification of nine equations (\ref{6}), (\ref{8})-(\ref{14}), (\ref{41}) divulge that
\begin{eqnarray}\label{42}
\begin{split}
&2(a_{11}a_{22}a_{33}+ a_{12}a_{23}a_{31}+a_{21}a_{32}a_{13}-a_{11}a_{23}a_{32}-a_{22}a_{31}a_{13}-a_{33}a_{12}a_{21})\\
&-(a_{22}a_{33}+a_{12}a_{23}+a_{13}a_{32}-a_{23}a_{32}-a_{12}a_{33}-a_{13}a_{22})b_1\\
&-(a_{33}a_{11}+a_{23}a_{31}+a_{21}a_{13}-a_{31}a_{13}-a_{23}a_{11}-a_{21}a_{33})b_2\\
&-(a_{11}a_{22}+a_{31}a_{12}+a_{32}a_{21}-a_{12}a_{21}-a_{31}a_{22}-a_{32}a_{11})b_3=0
\end{split}
\end{eqnarray}
and
\begin{eqnarray}\label{43}
\begin{split}
&-2(a_{11}a_{22}a_{33}+a_{12}a_{23}a_{31}+a_{21}a_{32}a_{13}-a_{11}a_{23}a_{32}-a_{22}a_{31}a_{13}-a_{33}a_{12}a_{21})\\
&+a_{11}a_{22}+a_{33}a_{11}+a_{22}a_{33}-a_{12}a_{21}-a_{23}a_{32}-a_{31}a_{13}=\frac{1}{12}.
\end{split}
\end{eqnarray}
Eqs. (\ref{41}) and (\ref{42}) together implies $|G_{N,3}(\sigma)|=1$. Thus all schemes satisfying these conditions will automatically lead to minimum dissipation error for entire range of $\sigma$ values. Further
\begin{eqnarray}\label{44}
\bm {arg}(G_{N,3}(\sigma))=2\tan^{-1}\left(
    \frac{
    \splitfrac{
      \sigma^2(a_{11}a_{22}+a_{22}a_{33}+a_{33}a_{11}
    }
    {
      -a_{12}a_{21}-a_{23}a_{32}-a_{31}a_{13})-1
    }
    }
    {
    \splitfrac{
      {\sigma/2}-\sigma^3 (a_{11}a_{22}a_{33}+a_{12}a_{23}a_{31}+a_{21}a_{32}a_{13}
      }{
        -a_{11}a_{23}a_{32}-a_{22}a_{31}a_{13}-a_{33}a_{12}a_{21})
      }
      }
  \right)+\pi.
\end{eqnarray}
Using Eq. (\ref{43}) we obtain
\begin{eqnarray}\label{44.1}
\bm {arg}(G_{N,3}(\sigma))=2\tan^{-1}\left({\frac{2 (\sigma^2X-1)}{\sigma-\sigma^3(X-\frac{1}{12})}}\right)+\pi
\end{eqnarray}
where
\begin{eqnarray}\label{45}
X=a_{11}a_{22}+a_{22}a_{33}+a_{33}a_{11}-a_{12}a_{21}-a_{23}a_{32}-a_{31}a_{13}.
\end{eqnarray}
We propose to formulate a new error function for all fourth order minimal dissipation three stage implicit method and is expressed as
\begin{eqnarray}\label{46}
\|PE(X)\|_{L^2[0,\pi]}=\left[\int_{0}^{\pi}\left|\sigma-2\tan^{-1}\left({\frac{2 (\sigma^2X-1)}{\sigma-\sigma^3(X-\frac{1}{12})}}\right)-\pi\right|^2\left|e^{-\alpha\sigma^2}\right|^2 d\sigma\right]^{1/2}.
\end{eqnarray}
The minimum value of $\|PE(X)\|_{L^2[0,\pi]}$ and the point of minima $X_{min}$ varies with $\alpha$. We plot $X_{min}$ against $\alpha$ in Figure \ref{fig:3S_phase_com}(a).
From the figure it is clear that point of minima $X_{min}$ increases with $\alpha$ only to asymptotically approach $X=1/10$ replicating behaviour of $Y_{min}$ observed in two stage schemes. Normalized absolute value of $\phi=\sigma-2\tan^{-1}\left({\frac{2 (\sigma^2X_{min}-1)}{\sigma-\sigma^3(X_{min}-\frac{1}{12})}}\right)-\pi$ is shown in Figure \ref{fig:3S_phase_com}(b). Again for the limiting case ($X_{min}=1/10$) phase error is found to be indistinguishable from Gauss-Legendre scheme (IRK36) \cite{but_08}, similar to two stage methods. Thus continuing from earlier section it seems to us that Gauss-Legendre methods, as they are derived, are inherently optimized to carry minimum dissipation error at small angular frequency. To the best of our knowledge, this singular character of optimally accurate IRK24 and IRK36 methods is not documented in the literature. A correct mathematically validation of the asymptotic approach of $X_{min}$ to $1/10$, as seen in figure \ref{fig:3S_phase_com}(a) can be found by using series expansion of
\begin{eqnarray}\label{46.0}
\phi(\sigma)=\left(\frac{1-10X}{120}\right)\sigma^5-\left(\frac{336X^2-84X+5}{4032}\right)\sigma^7+O(\sigma^9).
\end{eqnarray}
With emphasis zeroed down to very small values of $\sigma$, only leading term in the above expression remains relevant and is minimum at $X=1/10$ justifying the variation of phase error for the case $\alpha\rightarrow\infty$. Indeed at $X=1/10$ the three stage R-K methods carries sixth order dispersion error which is highest possible for a $A$-stable three stage implicit method. As in two stage here again we work with $\alpha=0$, $4$, $16$, and $\alpha\rightarrow\infty$ to understand effect of phase reduction.

\subsection{$\alpha=0$}
We begin our quest with $\alpha=0$. For this case minimum value of phase error $\|PE\|_{{L^2[0,\pi]}_{\min}}$ is attained at
\begin{eqnarray}\label{46.1}
X=X_{min}=0.1010711100.
\end{eqnarray}
Corresponding dispersion error is found to be $\phi_{[0, \pi]}=1.781038\times10^{-3}$. The nine equations (\ref{6}), (\ref{8})-(\ref{13}), (\ref{41}) and (\ref{46.1}) are to be satisfied by twelve coefficients for a three stage implicit scheme with overall minimum dissipation and dispersion error. Solving we get infinite set of solutions all with dimension three. A few such schemes have been shown in Table \ref{table 4_0}. The first scheme S3A1 can be realized using additional constraints $b_1=b_2$, $a_{12}=0=a_{13}$. For second set S3A2 we use $b_1=b_3$, $b_2=2a_{22}$, $a_{11}=a_{33}$ and $a_{13}+a_{31}=5/18$ which is partially motivated by the coefficient of IRK36. S3A3 and S3A4 are obtained with $b_2=b_3$, $a_{23}=0=a_{13}$ and $a_{11}=1/4$, $a_{11}=2a_{12}$, $a_{22}=2a_{21}$ respectively.
\begin{table}[h!]
\caption{Fourth order three stage low-dissipation low-dispersion implicit R-K (S3A) schemes obtained with weight parameter $\alpha=0$.}
\centering
\begin{tabular}{c c c c c}
 \hline
  			&\multicolumn{4}{c}{Schemes}\\ \cline{2-5}
Parameter 	& S3A1 		  & S3A2 		 	 &S3A3 		   & S3A4 \\
 \hline
 $b_1$ 		&0.4902164042 & 0.2777777778 &0.4990278482 & 0.6702568370 \\
 $b_2$ 		&0.4902164042 & 0.4444444444 &0.2504860759 & 1.5072733738 \\
 $b_3$ 		&0.0195671916 & 0.2777777778 &0.2504860759 &-1.1775302108 \\
 $a_{11}$ 	&0.2267610814 & 0.1388888889 &0.2548461218 & 0.2500000000 \\
 $a_{12}$ 	&0.0000000000 & -0.0386992007 &-0.0438954380& 0.1250000000 \\
 $a_{13}$ 	&0.0000000000 & 0.0125119772 &0.0000000000 &-0.1193016952 \\
 $a_{21}$ 	&0.5149632492 & 0.3019647782 &0.7842232807 & 0.5223474224 \\
 $a_{22}$ 	&0.2396583441 & 0.2222222222 &0.0183927967 & 1.0446948445 \\
 $a_{23}$ 	&0.0381882637 & -0.0241870005 &0.0000000000 &-0.8704412903 \\
 $a_{31}$ 	&0.7895342543 & 0.2652658006 &0.2800365570 & 0.3872607826 \\
 $a_{32}$ 	&-0.8134251058& 0.4831436452 &0.2664412801& 1.0200312043 \\
 $a_{33}$ 	&0.0335805745 & 0.1388888889 &0.2267610814 &-0.7946948449 \\
 \hline
\end{tabular}
\label{table 4_0}
\end{table}

\subsection{$\alpha=4$}
The choice of $\alpha=4$ can be seen to balance relative prominence of angular frequency. At $\alpha=4$ minimization is materialized at
\begin{eqnarray}\label{47}
X_{min}=0.1000815539
\end{eqnarray}
with dispersion error $\phi_{[0, \pi]}=8.878927\times10^{-3}$. Four such methods out of many possible, derived as discussed in section 4.1, have been tabulated below.
\begin{table}[h!]
\caption{Fourth order three stage low-dissipation low-dispersion implicit R-K (S3B) schemes obtained with weight parameter $\alpha=4$.}
\centering
\begin{tabular}{c c c c c}
 \hline
  			&\multicolumn{4}{c}{Schemes}\\ \cline{2-5}
Parameter 	& S3B1 		  & S3B2 		 	 &S3B3 		   & S3B4 \\
 \hline
 $b_1$ 		&0.4968572595 & 0.2777777778 &0.0062854810 & 0.6655575665 \\
 $b_2$ 		&0.4968572595 & 0.4444444444 &0.4968572595 & 1.4996524498 \\
 $b_3$ 		&0.0062854810 & 0.2777777778 &0.4968572595 &-1.1652100163 \\
 $a_{11}$ 	&0.2162822020 & 0.1388888889 &0.0441899847 & 0.2500000000 \\
 $a_{12}$ 	&0.0000000000 & -0.0361869818 &0.9526770898& 0.1250000000 \\
 $a_{13}$ 	&0.0000000000 & 0.0099997583 &0.0000000000 &-0.1202281088 \\
 $a_{21}$ 	&0.5384237709 & 0.3003946414 &-0.0295312165 & 0.5184708679 \\
 $a_{22}$ 	&0.2395278133 & 0.2222222222 &0.2395278133 & 1.0369417358 \\
 $a_{23}$ 	&0.0120518190 & -0.0226168636 &0.0000000000 &-0.8579660534 \\
 $a_{31}$ 	&2.2933385501 & 0.2677780195 &0.0290118251 & 0.3844142138 \\
 $a_{32}$ 	&-2.3343956093& 0.4806314263 &0.5384237709& 1.0165737263 \\
 $a_{33}$ 	&0.0441899847 & 0.1388888889 &0.2162822020 &-0.7869417358 \\
 \hline
\end{tabular}
\label{table 4}
\end{table}

\subsection{$\alpha=16$}
With $\alpha$ fixed at 16, we obtain
\begin{eqnarray}\label{47.1}
X_{min}=0.1000204444
\end{eqnarray}
and the corresponding overall dispersion error increases slightly to $\phi_{[0, \pi]}=9.400444\times10^{-3}$. Again the nine equations (\ref{6}), (\ref{8})-(\ref{13}), (\ref{41}) and (\ref{47.1}) admit infinite solutions, some of which have been represented in table \ref{table 4_16}.
\begin{table}[h!]
\caption{Fourth order three stage low-dissipation low-dispersion implicit R-K (S3C) schemes obtained with weight parameter $\alpha=16$.}
\centering
\begin{tabular}{c c c c c}
 \hline
  			&\multicolumn{4}{c}{Schemes}\\ \cline{2-5}
Parameter 	& S3C1 		  & S3C2 		 &S3C3 		   &S3C4 \\
 \hline
 $b_1$ 		&0.4973158852 & 0.2777777778 &0.0053682296 & 0.6652690242 \\
 $b_2$ 		&0.4973158852 & 0.4444444444 &0.4973158852 & 1.4991870878 \\
 $b_3$ 		&0.0053682296 & 0.2777777778 &0.4973158852 &-1.1644561120 \\
 $a_{11}$ 	&0.2155587380 & 0.1388888889 &0.0449040217 & 0.2500000000 \\
 $a_{12}$ 	&0.0000000000 & -0.0360294386 &0.9524190295& 0.1250000000 \\
 $a_{13}$ 	&0.0000000000 & 0.0098422150 &0.0000000000 &-0.1202853659 \\
 $a_{21}$ 	&0.5399915310 & 0.3002961769 &-0.0293468092 & 0.5182333238 \\
 $a_{22}$ 	&0.2395372403 & 0.2222222222 &0.2395372403 & 1.0364666476 \\
 $a_{23}$ 	&0.0102807976 & -0.0225183991 &0.0000000000 &-0.8572009786 \\
 $a_{31}$ 	&2.6764782769 & 0.2679355627 &0.0288909930 & 0.3842400678 \\
 $a_{32}$ 	&-2.7187053498& 0.4804738830 &0.5399915310& 1.0163632850 \\
 $a_{33}$ 	&0.0449040217 & 0.1388888889 &0.2155587380 &-0.7864666476 \\
 \hline
\end{tabular}
\label{table 4_16}
\end{table}

\subsection{$\alpha\rightarrow\infty$}
For three stage methods as $\alpha\rightarrow\infty$ we see that
\begin{eqnarray}\label{47.2}
X_{min}\rightarrow\frac{1}{10}
\end{eqnarray}
and the dispersion error become $\phi_{[0, \pi]}=9.575026\times10^{-3}$. This minimum value automatically reveal highest possible sixth order dispersion error via Eq. (\ref{46.0}) for three stage implicit methods. All these methods possesses identical dispersion error. Four schemes are reported in table \ref{table 4_INF}. The second set of values carry special significance since they correspond to Gauss-Legendre three stage methods and carry overall accuracy of order six.
\begin{table}[h!]
\caption{Fourth order three stage low-dissipation low-dispersion implicit R-K (S3D) schemes obtained with weight parameter $\alpha\rightarrow\infty$.}
\centering
\begin{tabular}{c c c c c}
 \hline
  			&\multicolumn{4}{c}{Schemes}\\ \cline{2-5}
Parameter 	& S3D1 		  & S3D2 		 &S3D3 		   &S3D4 \\
 \hline
 $b_1$ 		&0.4974707660 & 0.2777777778 &0.0050584680 & 0.6651725342 \\
 $b_2$ 		&0.4974707660 & 0.4444444444 &0.4974707660 & 1.4990315365 \\
 $b_3$ 		&0.0050584680 & 0.2777777778 &0.4974707660 &-1.1642040707 \\
 $a_{11}$ 	&0.2153144231 & 0.1388888889 &0.0451446417 & 0.2500000000 \\
 $a_{12}$ 	&0.0000000000 & -0.0359766675 &0.9523324891& 0.1250000000 \\
 $a_{13}$ 	&0.0000000000 & 0.0097894440 &0.0000000000 &-0.1203045227 \\
 $a_{21}$ 	&0.5405195031 & 0.3002631950 &-0.0292850447 & 0.5181539006 \\
 $a_{22}$ 	&0.2395409352 & 0.2222222222 &0.2395409352 & 1.0363078012 \\
 $a_{23}$ 	&0.0096836713 & -0.0224854172 &0.0000000000 &-0.8569451582 \\
 $a_{31}$ 	&2.8373912650 & 0.2679883338 &0.0288516507 & 0.3841818496 \\
 $a_{32}$ 	&-2.8800130375& 0.4804211120 &0.5405195031& 1.0162929609 \\
 $a_{33}$ 	&0.0451446417 & 0.1388888889 &0.2153144231 &-0.7863078012\\
 \hline
\end{tabular}
\label{table 4_INF}
\end{table}
\subsection{Gauss-Legendre three stage method (IRK36)}
Gauss-Legendre three stage scheme having Butcher tableau representation
\begin{eqnarray}\label{53}
\begin{tabular}{c|c}
\begin{tabular}{c}
  $\frac{1}{2}-\frac{\sqrt{15}}{10}$ \\
  \\
  $\frac{1}{2}$ \\
  \\
  $\frac{1}{2}+\frac{\sqrt{15}}{10}$ \\
\end{tabular}
& \begin{tabular}{c c c}
            $\frac{5}{36}$ & $\frac{2}{9}-\frac{\sqrt{15}}{15}$ & $\frac{5}{36}-\frac{\sqrt{15}}{30}$ \\
            \\
            $\frac{5}{36}+\frac{\sqrt{15}}{24}$ & $\frac{2}{9}$ & $\frac{5}{36}-\frac{\sqrt{15}}{24}$ \\
            \\
            $\frac{5}{36}+\frac{\sqrt{15}}{30}$ & $\frac{2}{9}+\frac{\sqrt{15}}{15}$ & $\frac{5}{36}$ \\
          \end{tabular}
\\
  \hline
       & \begin{tabular}{c c c}
           $\frac{5}{18}\;\;\;\;\;\;\;\;\;\;$ & $\frac{4}{9}$ & $\;\;\;\;\;\;\;\;\;\;\frac{5}{18}$ \\
         \end{tabular}
\end{tabular}.
\end{eqnarray}					
are derived by taking coefficients $c_s$, $s=1,2,3$ as root of Legendre polynomial $\mathbb P^*_R$ \cite{but_08}. This scheme attains highest possible sixth order of accuracy and is denoted here as IRK36. It is seen that the coefficients of this methods are such that $X=\frac{1}{10}$ emphasising its sixth order dispersion accuracy. Thus it is safe to conclude that IRK36 carries optimized dispersion error for small values of $\sigma$, apart from infinite dissipation accuracy and is as good as that of any other method presented in table \ref{table 4_INF} amply demonstrated by the complete overlap of phase error variation. But the key for the IRK36 is that given infinite dissipation accuracy and sixth order dispersion accuracy, it is able to attain highest possible sixth order overall accuracy and is the only scheme to do so.

\subsection{Comparison of numerical characteristics}
\begin{figure}[!h]
\begin{minipage}[b]{.6\linewidth}\hspace{-1cm}
\centering\includegraphics[width=75mm]{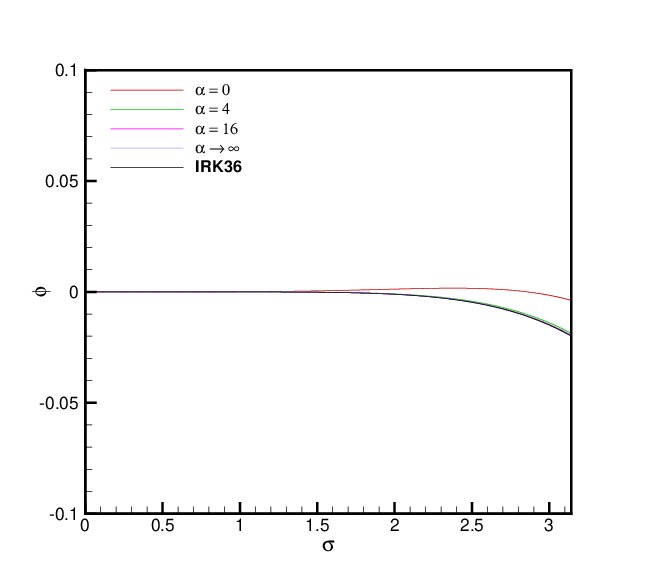}\\(a)
\end{minipage}
\begin{minipage}[b]{.6\linewidth}\hspace{-1cm}
\centering\includegraphics[width=75mm]{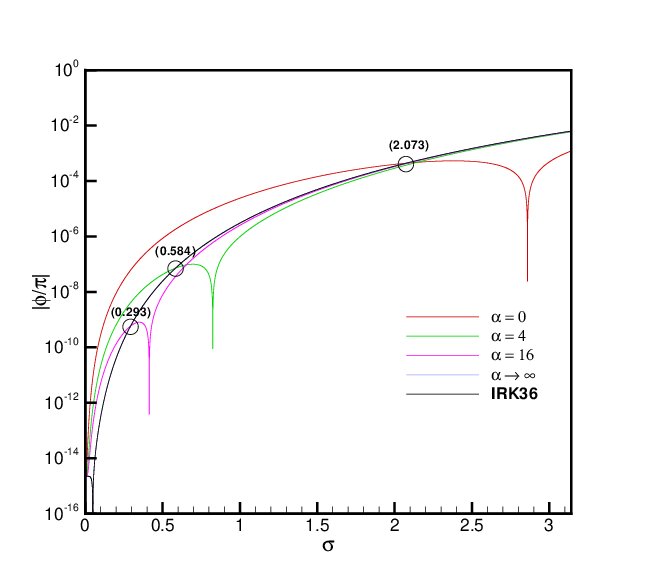}\\(b)
\end{minipage}
\begin{center}
\caption{{\sl (a) Phase difference and (b) dispersion error in logarithmic scale of various
schemes.} }
    \label{fig:3S_phase}
\end{center}
\end{figure}
We compare numerical characteristics of diverse groups of three stage schemes derived in this section. In figure \ref{fig:3S_phase}(a) we see that phase difference is minimum for S3A class of schemes, corresponding to the unity weight kernel. But a closer look in figure \ref{fig:3S_phase}(b) reveals that this set of schemes hold advantage over all other schemes only for $\sigma$ values higher than 2.073. In terms of wave number we can conclude that S3A holds advantage for waves having more than 3.031 time step per period. As the weight function is minimized by taking $\alpha=4$ we see that the new set of schemes S3B are less dispersive compared to other in the region $(0.584, 2.073)$. Correspondingly for S3C class of schemes dispersive error is minimum in the $\sigma$ range $(0.293, 0.584)$ vis-a-vis other groups of scheme.
\begin{figure}[p]
\vspace{-18 mm}
\begin{minipage}[b]{.6\linewidth}
\centering\psfig{file=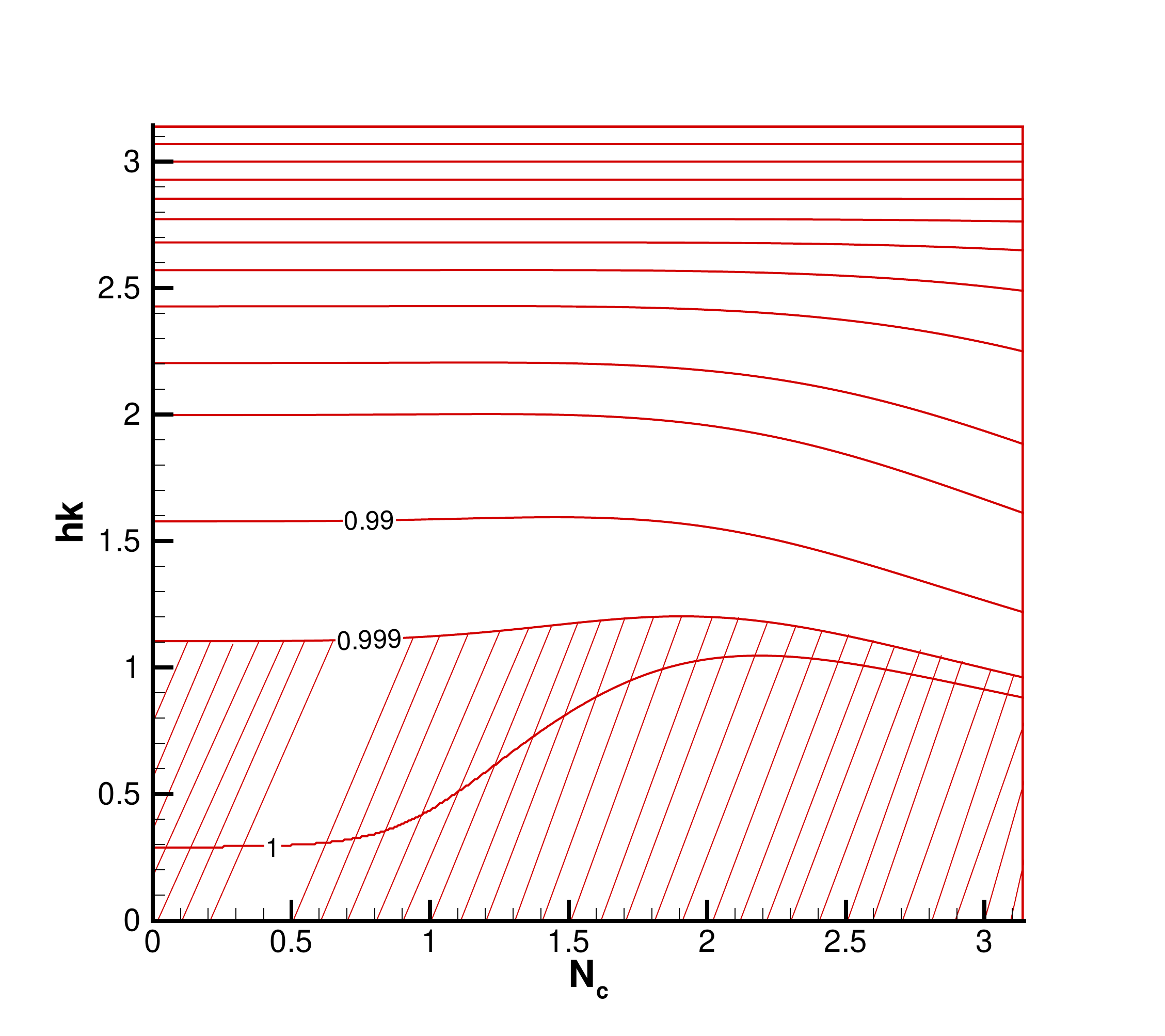,width=0.7\linewidth}
 (a)
\end{minipage}            \hspace{-2.5mm}
\begin{minipage}[b]{.6\linewidth}
\centering\psfig{file=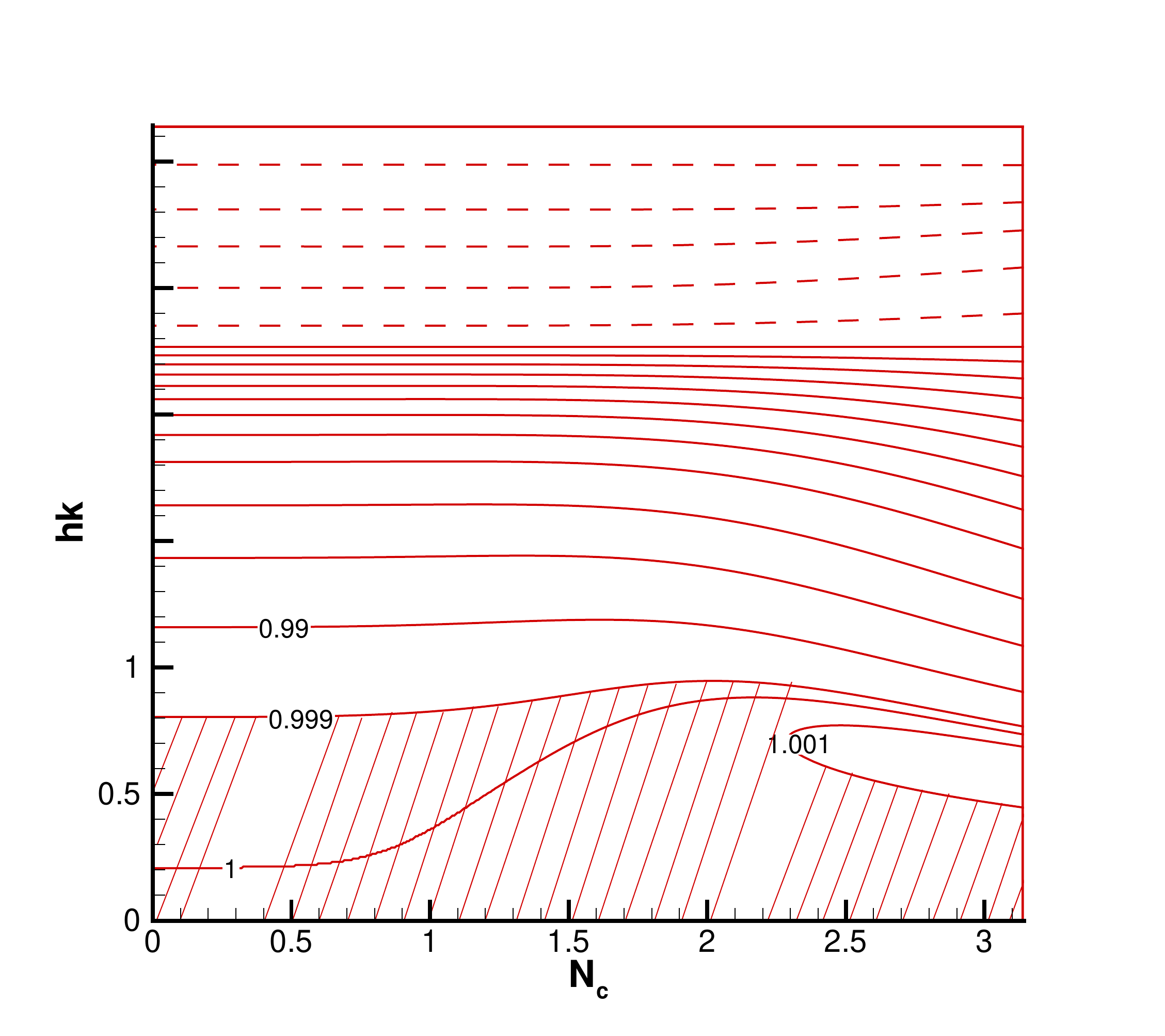,width=0.7\linewidth}
 (b)
\end{minipage}            \hspace{-2.5mm}
\begin{minipage}[b]{.6\linewidth}   \
\centering\psfig{file=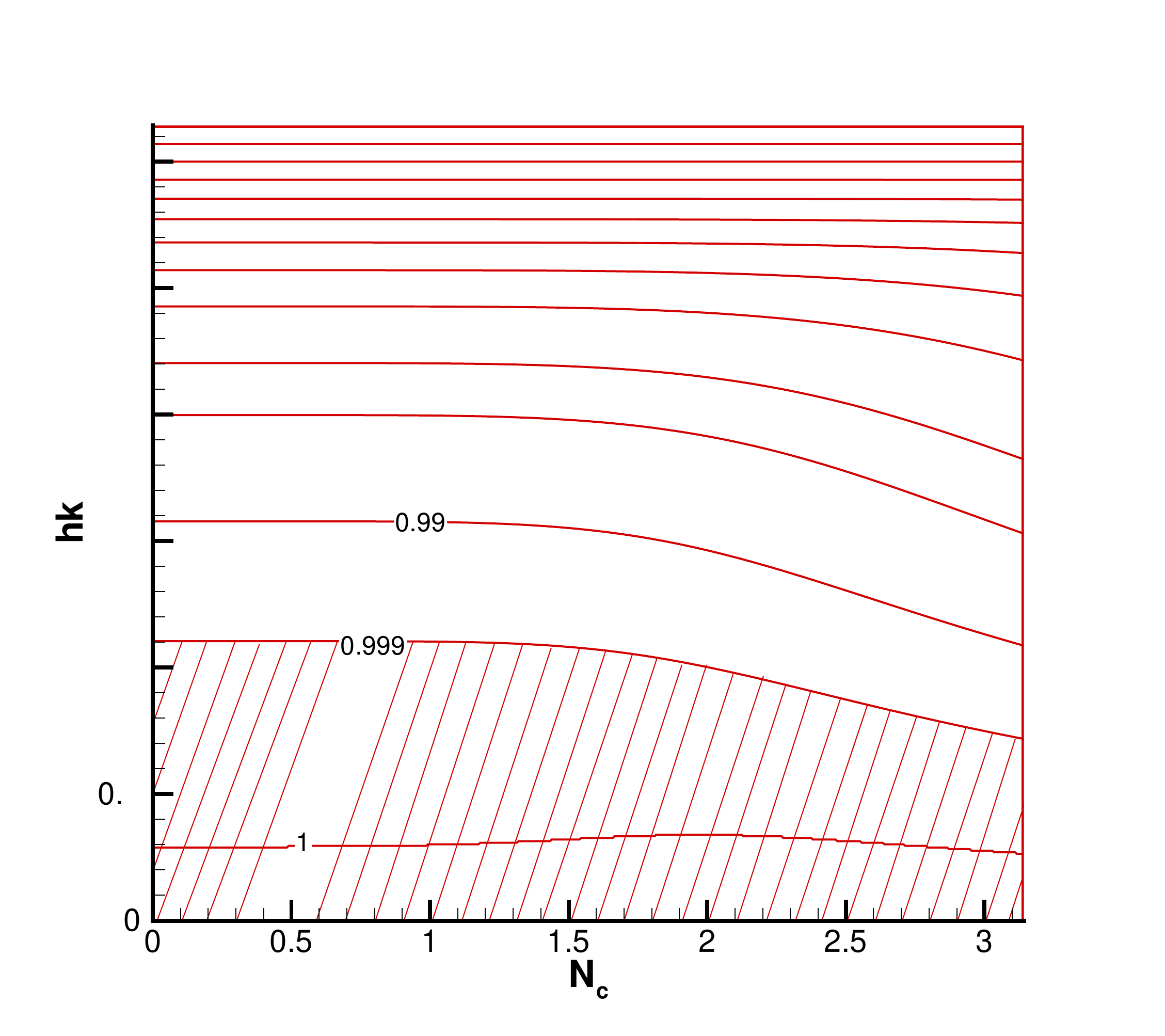,width=0.7\linewidth}
 (c)
\end{minipage}            \hspace{-2.5mm}
\begin{minipage}[b]{.6\linewidth}
\centering\psfig{file=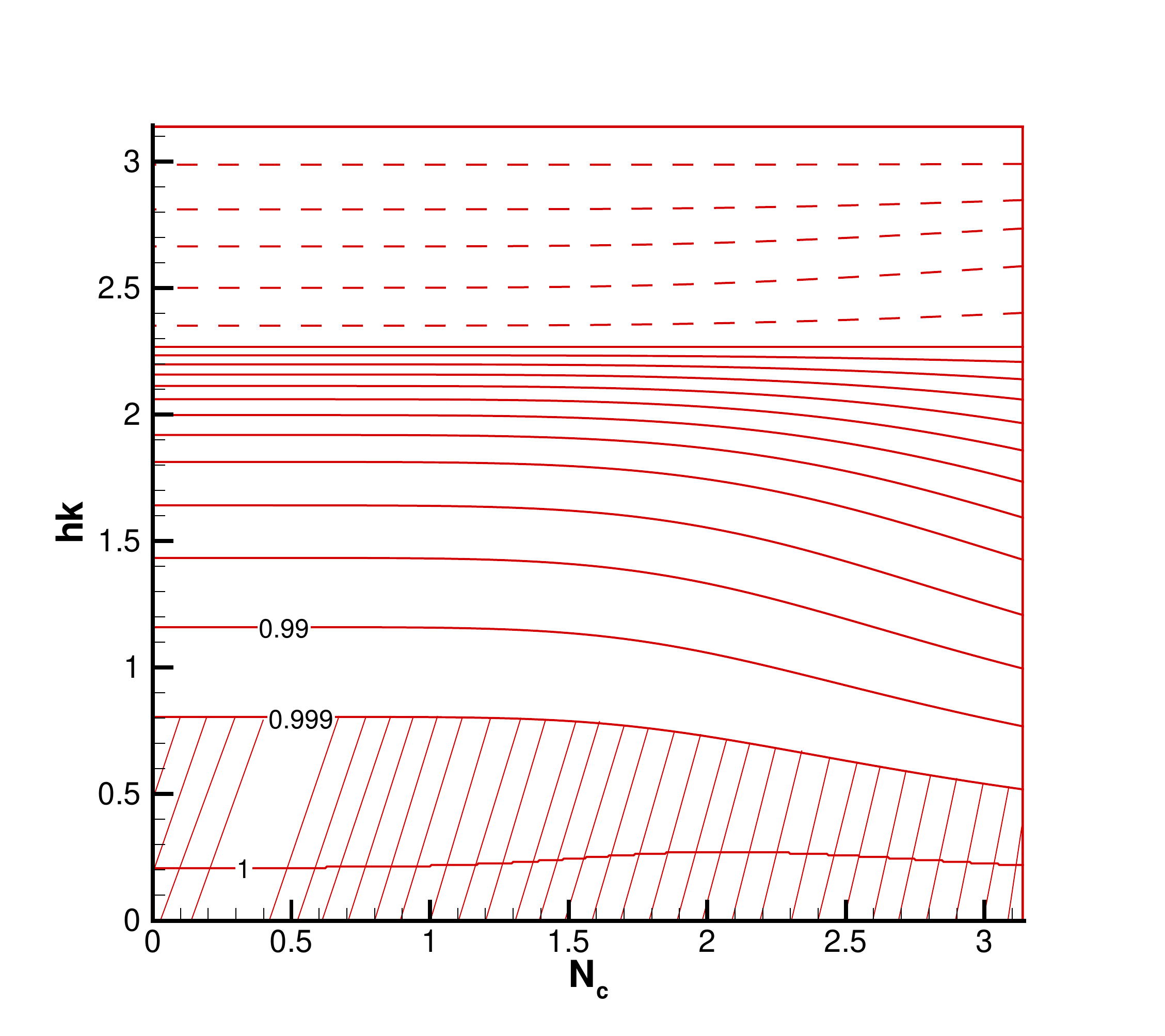,width=0.7\linewidth}
 (d)
\end{minipage}            \hspace{-2.5mm}
\begin{minipage}[b]{.6\linewidth}   \
\centering\psfig{file=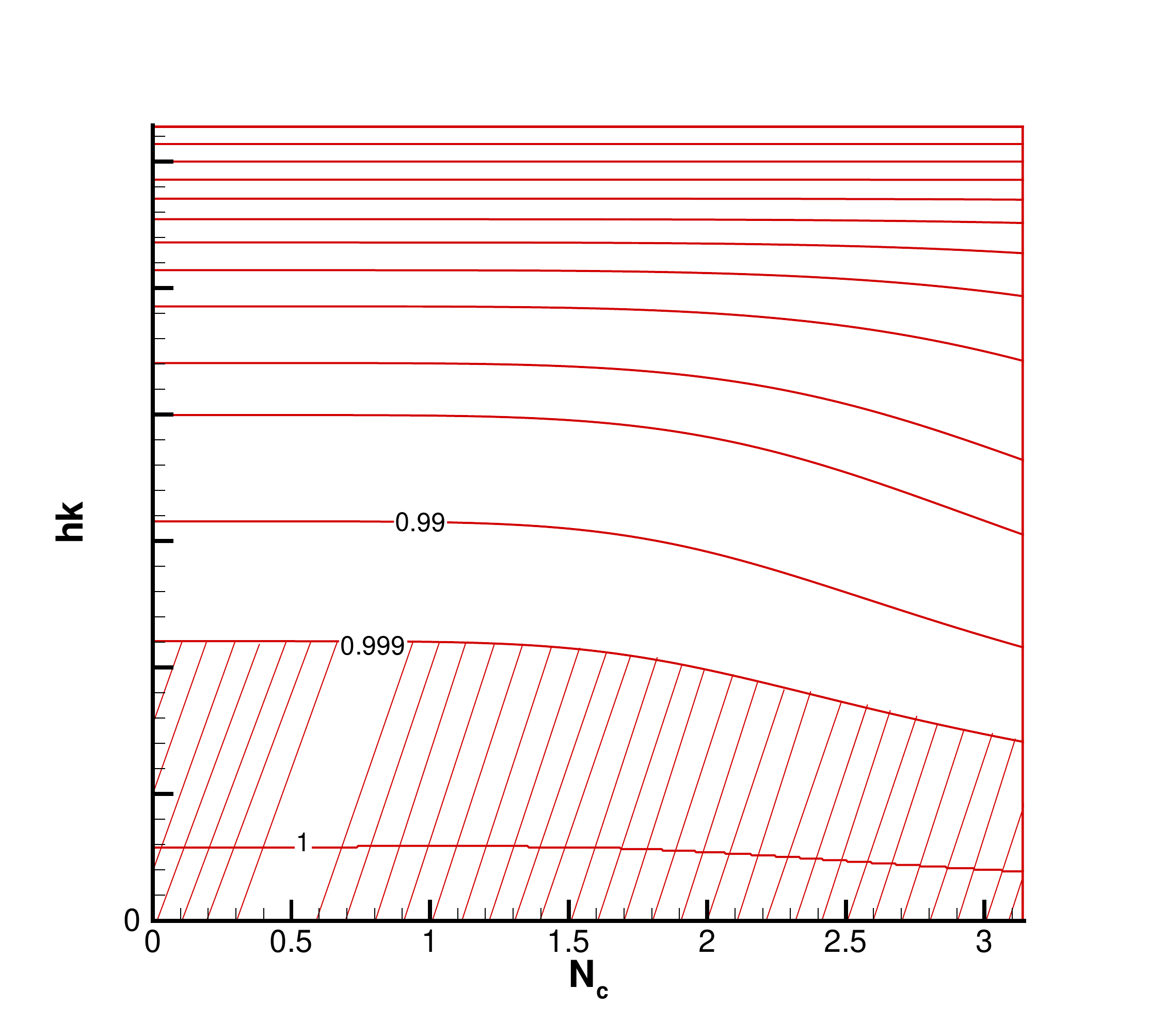,width=0.7\linewidth}
 (e)
\end{minipage}            \hspace{-2.5mm}
\begin{minipage}[b]{.6\linewidth}
\centering\psfig{file=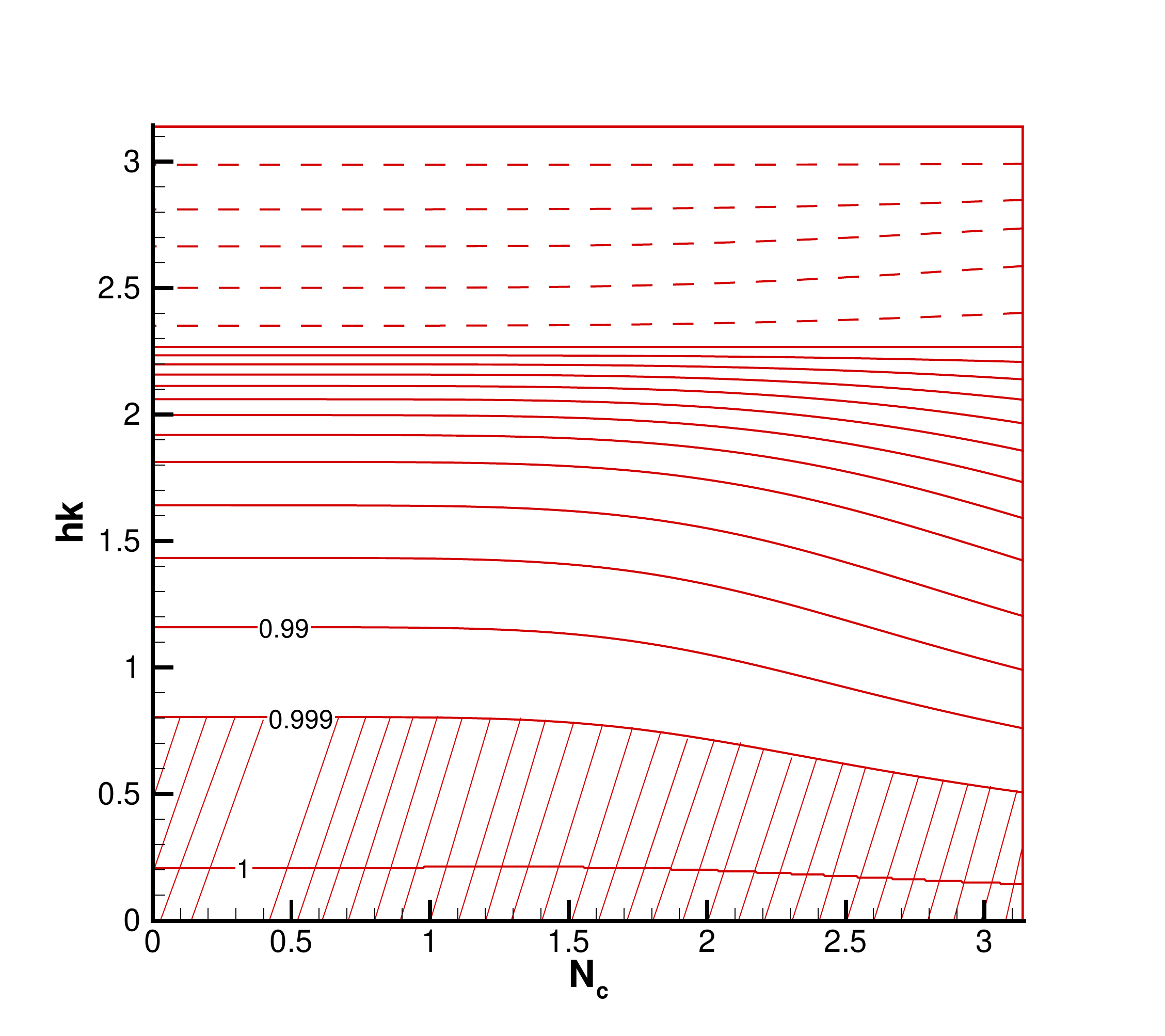,width=0.7\linewidth}
 (f)
\end{minipage}            \hspace{-2.5mm}
\begin{minipage}[b]{.6\linewidth}   \
\centering\psfig{file=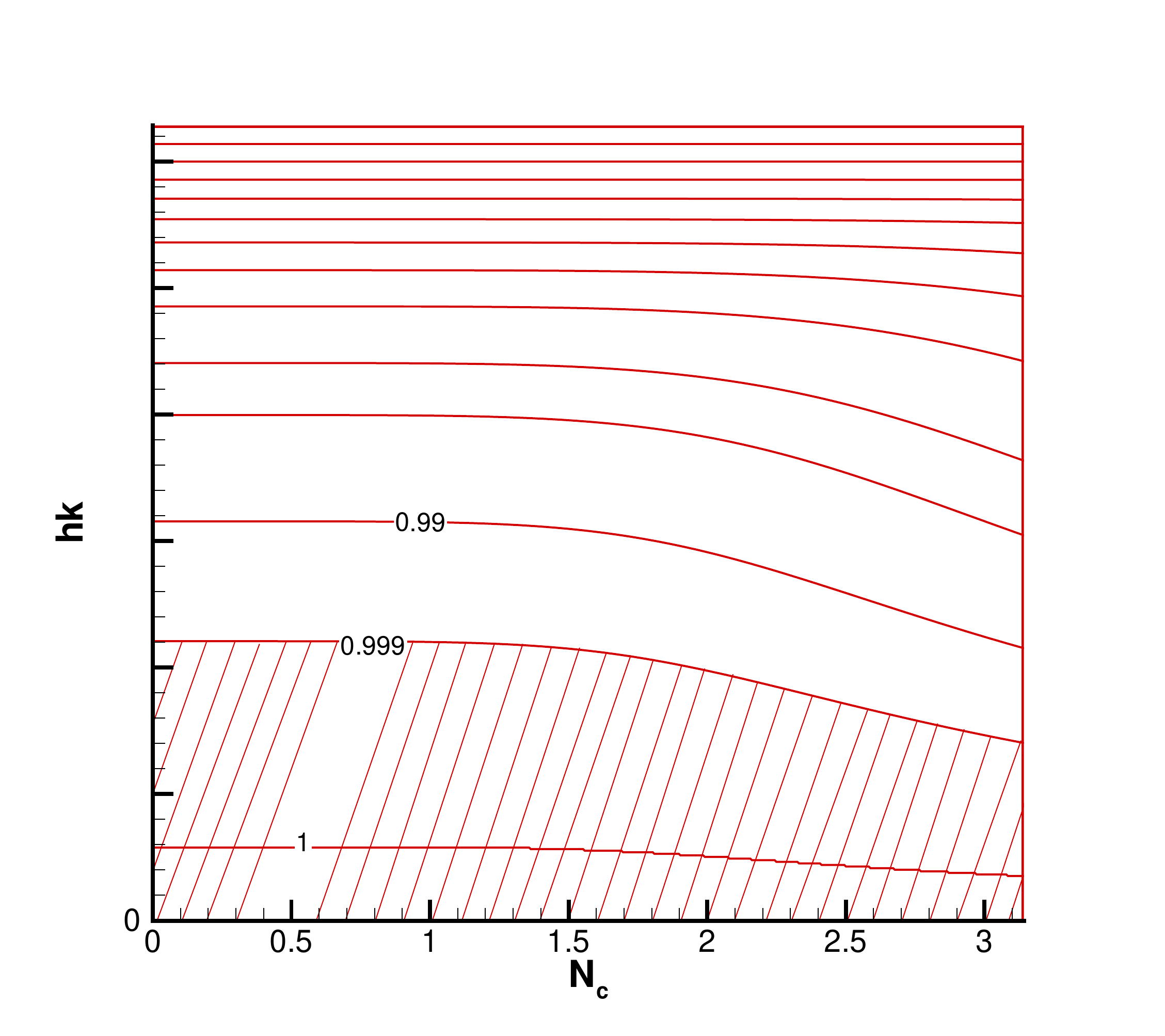,width=0.7\linewidth}
 (g)
\end{minipage}            \hspace{-2.5mm}
\begin{minipage}[b]{.6\linewidth}
\centering\psfig{file=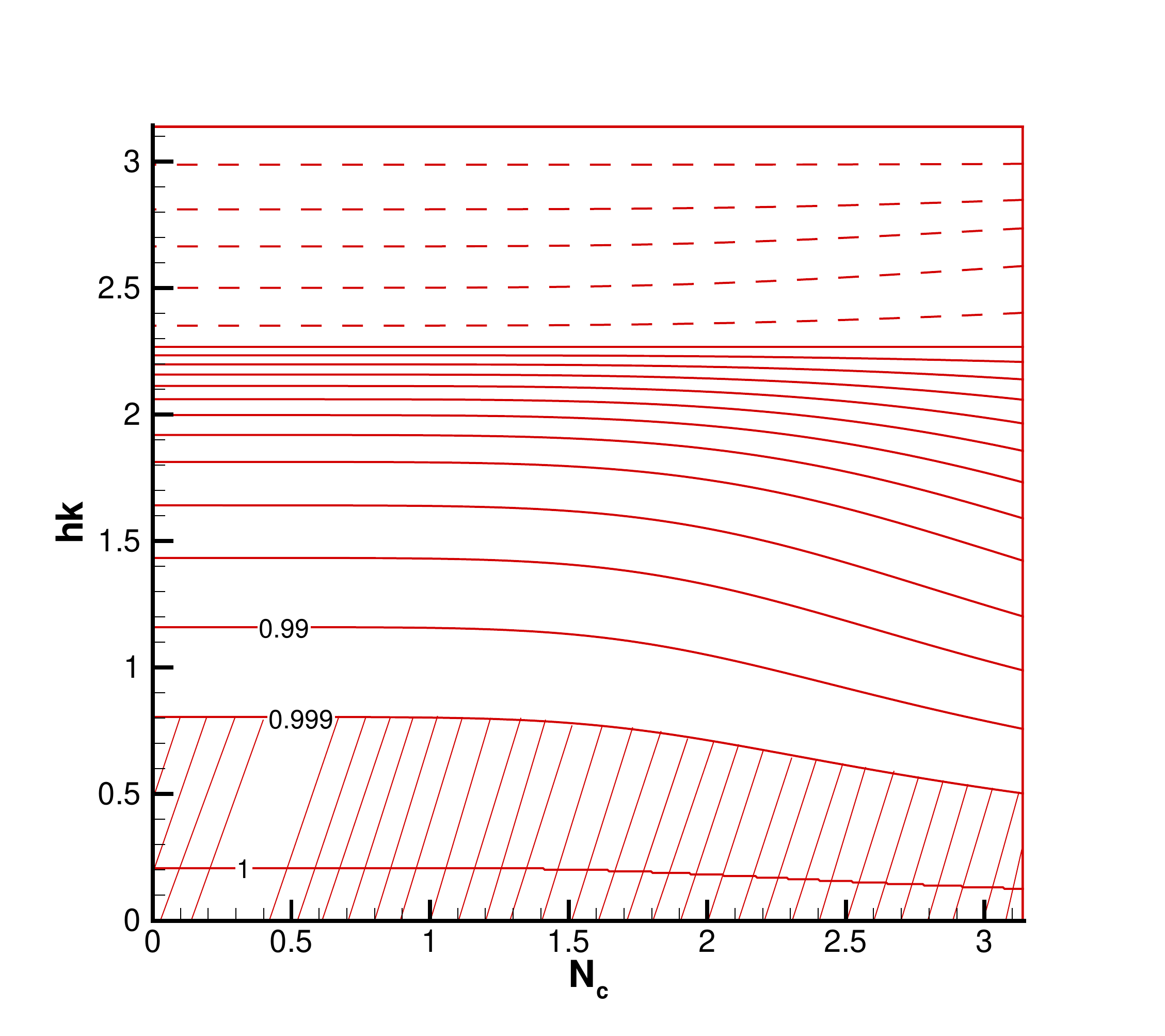,width=0.7\linewidth}
 (h)
\end{minipage}            \hspace{-2.5mm}
\begin{center}
\caption{{\sl Contours of normalized numerical phase velocity (left) and group velocity (right) for indicated schemes plotted in the $(N_c, hk)$ plane at mid-node when Lele scheme is used for spatial discretization: (a)-(b) S3A $(\alpha=0)$, (c)-(d) S3B $(\alpha=4)$, (e)-(f) S3C $(\alpha=16)$, (g)-(h) IRK36, S3D $(\alpha\rightarrow\infty)$.}}
\label{fig:3S_impl}
\end{center}
\end{figure}

\begin{figure}[p]
\vspace{-18 mm}
\begin{minipage}[b]{.6\linewidth}
\centering\psfig{file=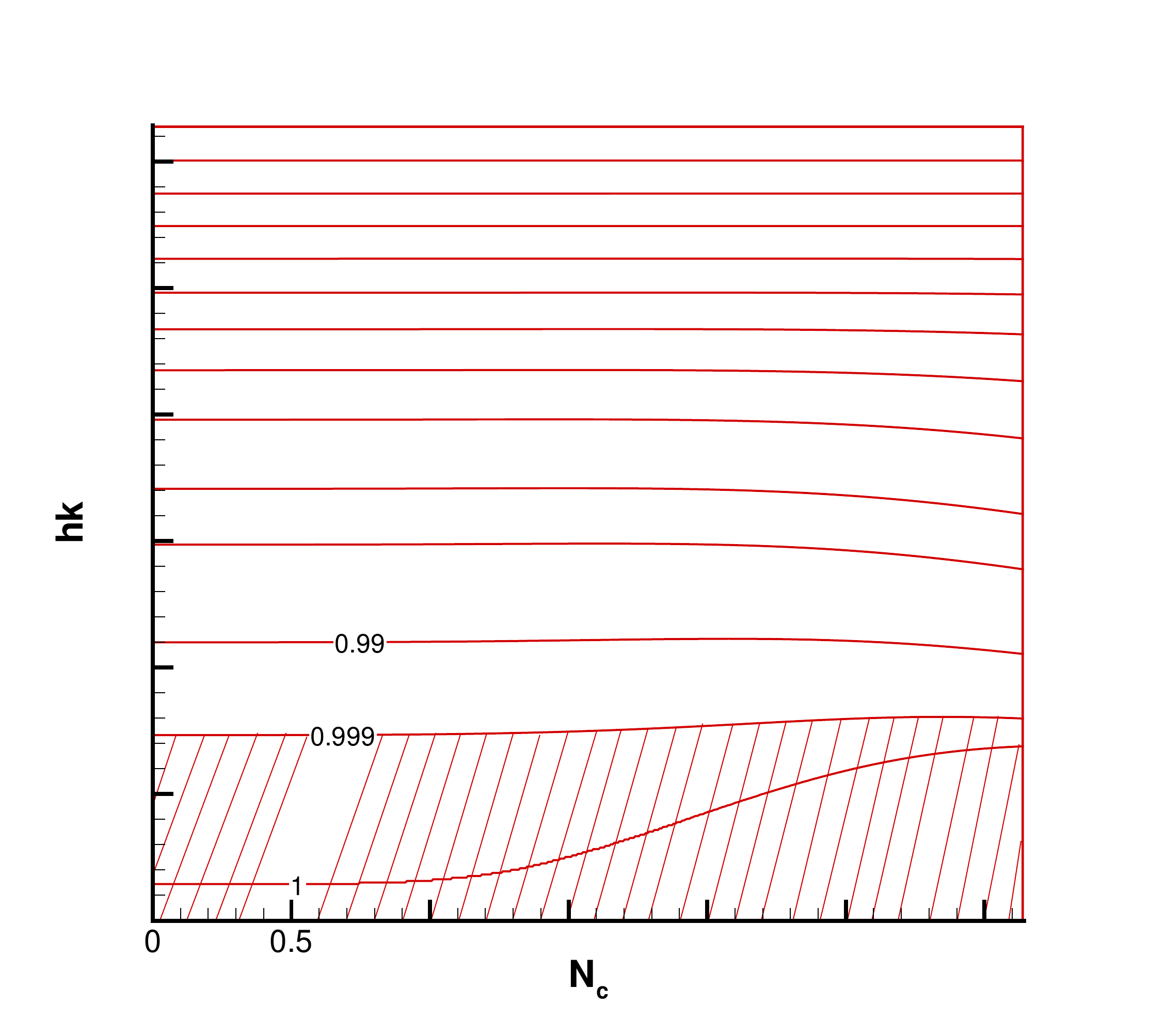,width=0.7\linewidth}
 (a)
\end{minipage}            \hspace{-2.5mm}
\begin{minipage}[b]{.6\linewidth}
\centering\psfig{file=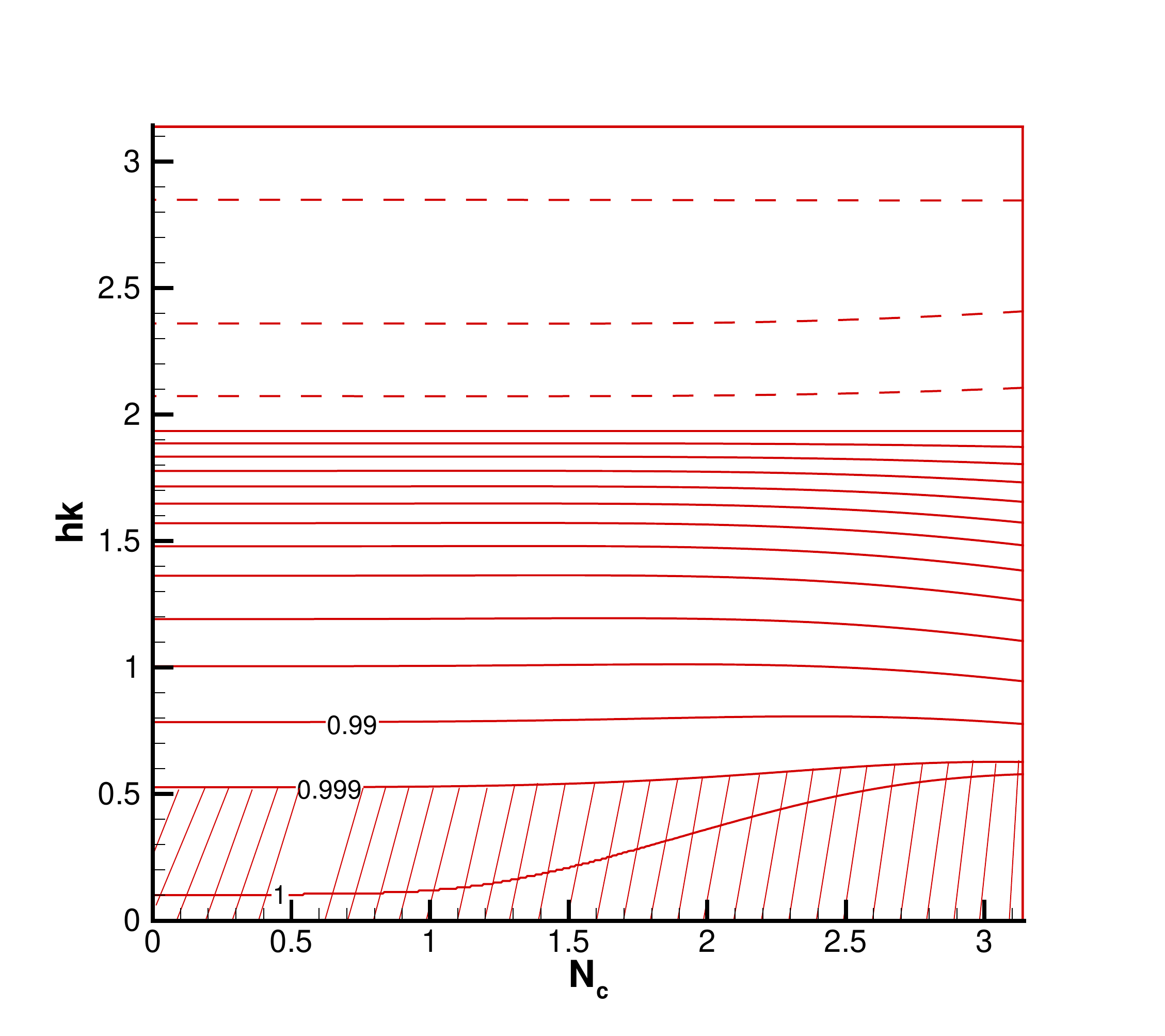,width=0.7\linewidth}
 (b)
\end{minipage}            \hspace{-2.5mm}
\begin{minipage}[b]{.6\linewidth}   \
\centering\psfig{file=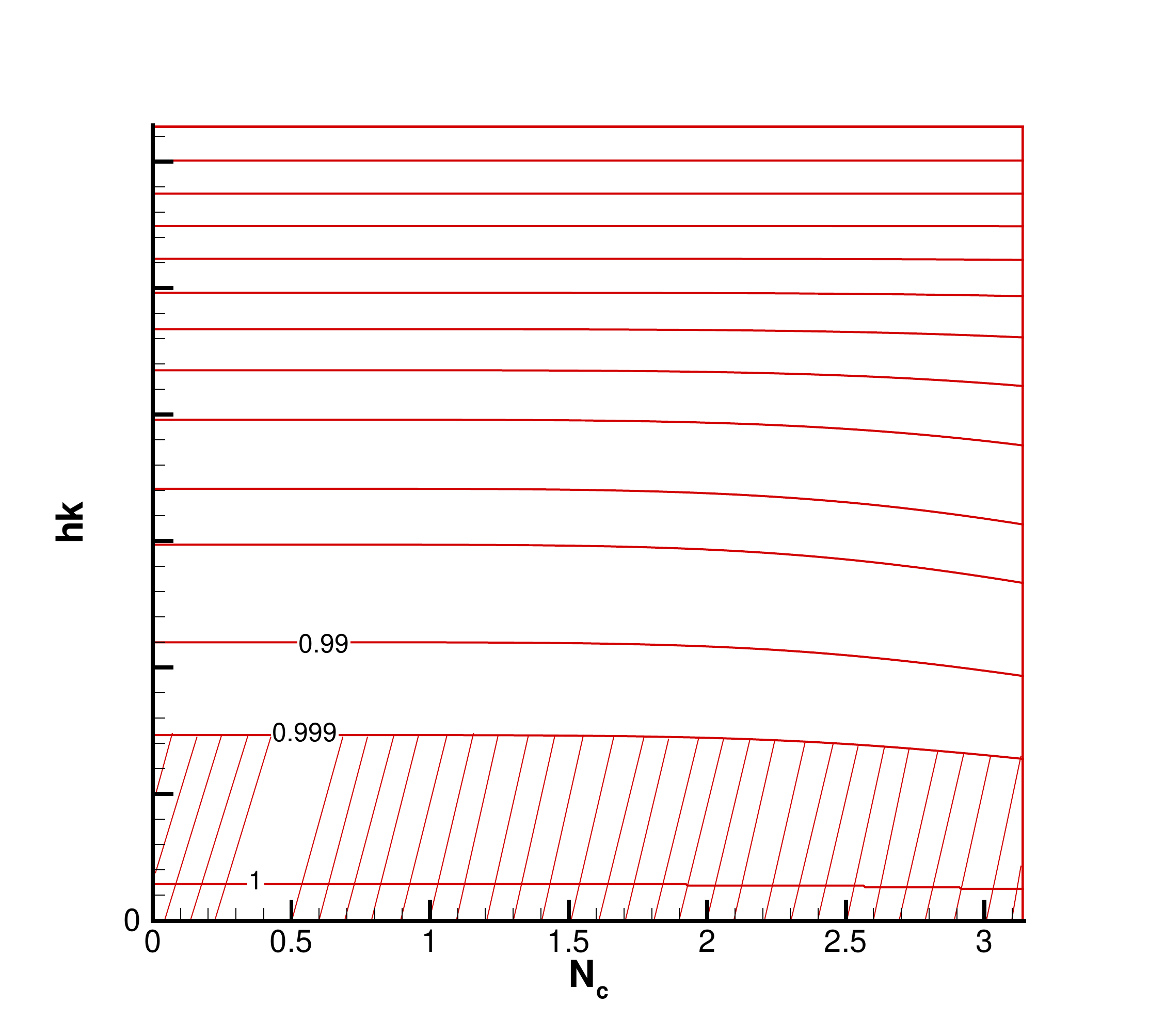,width=0.7\linewidth}
 (c)
\end{minipage}            \hspace{-2.5mm}
\begin{minipage}[b]{.6\linewidth}
\centering\psfig{file=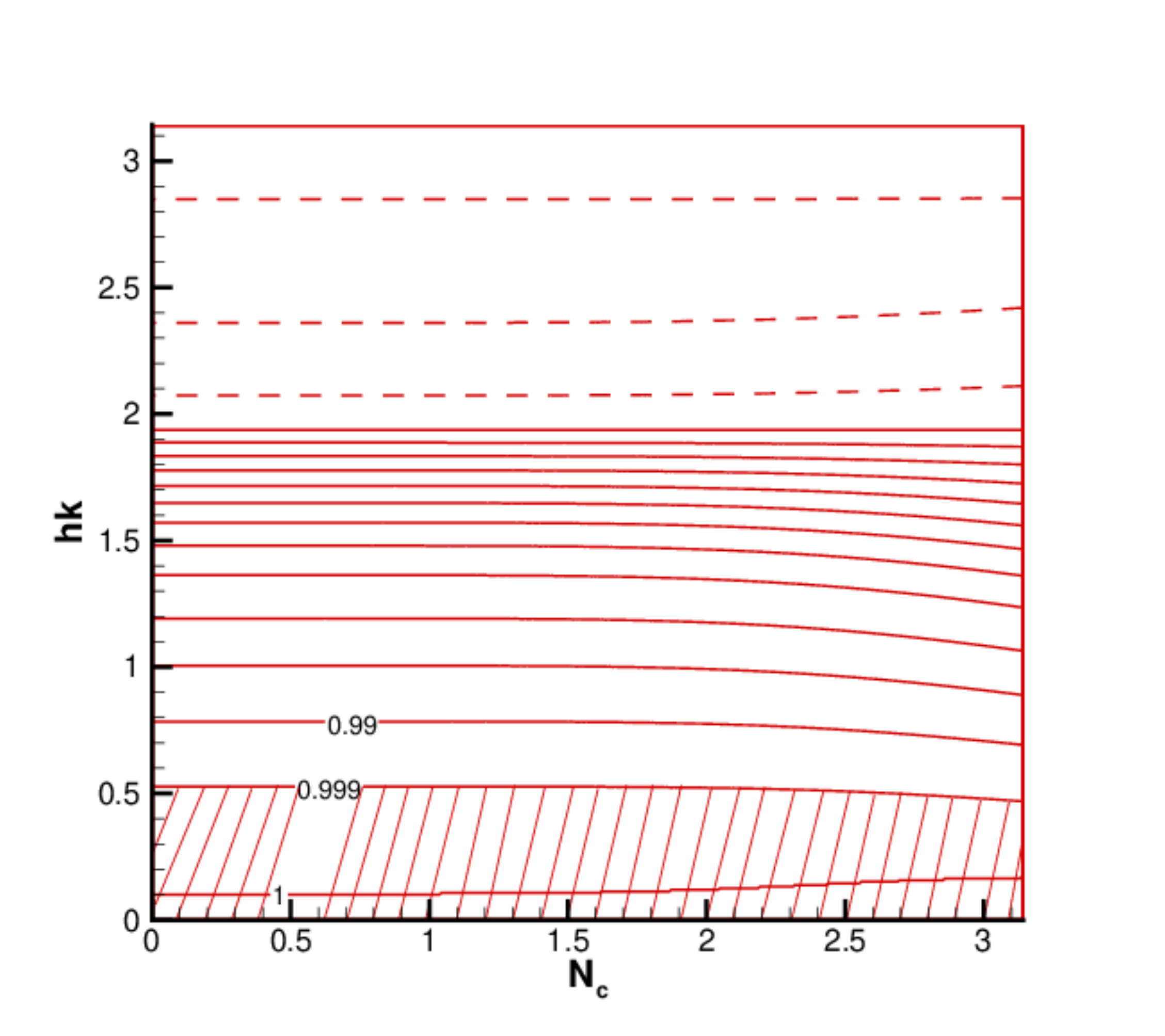,width=0.7\linewidth}
 (d)
\end{minipage}            \hspace{-2.5mm}
\begin{minipage}[b]{.6\linewidth}   \
\centering\psfig{file=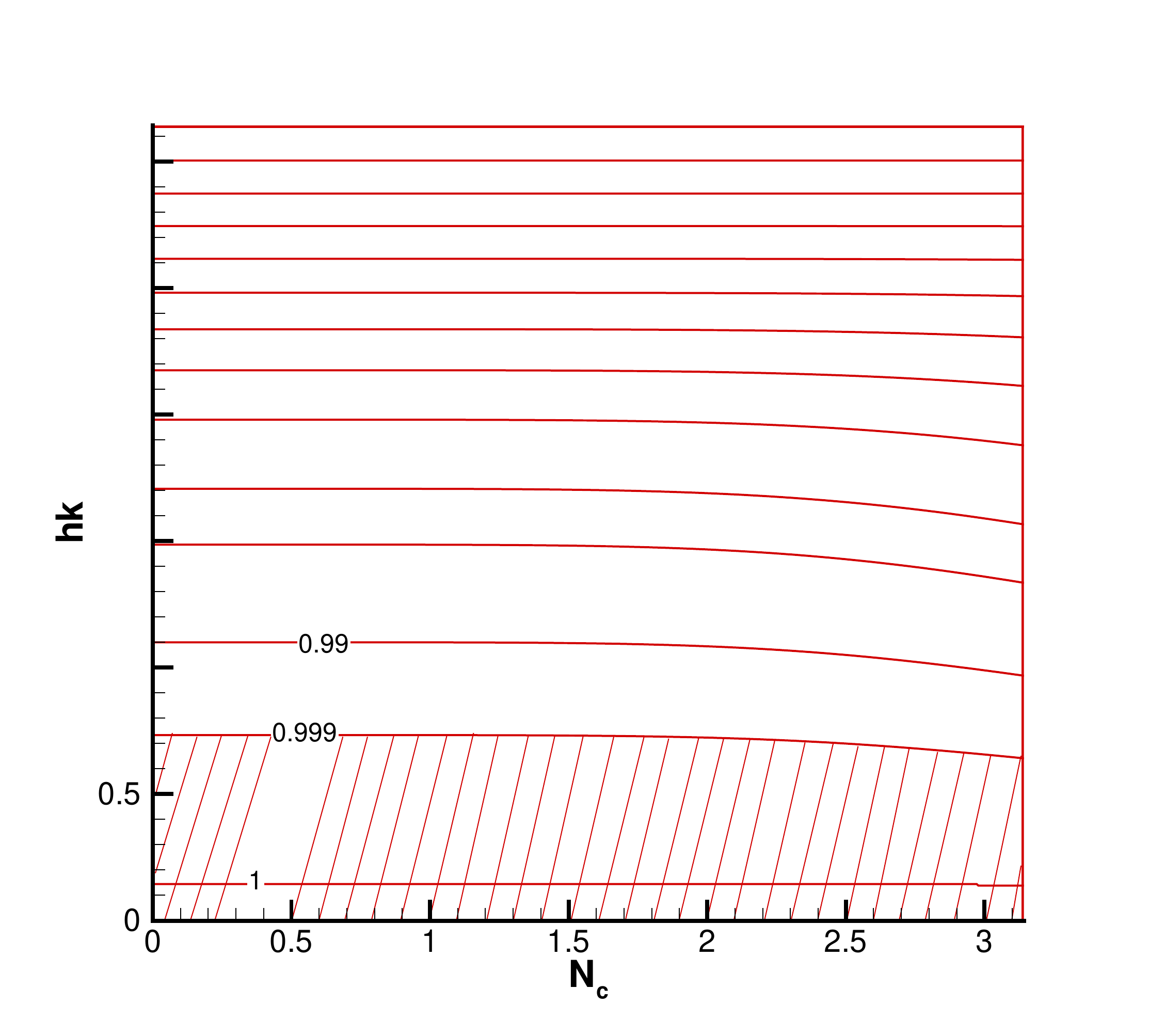,width=0.7\linewidth}
 (e)
\end{minipage}            \hspace{-2.5mm}
\begin{minipage}[b]{.6\linewidth}
\centering\psfig{file=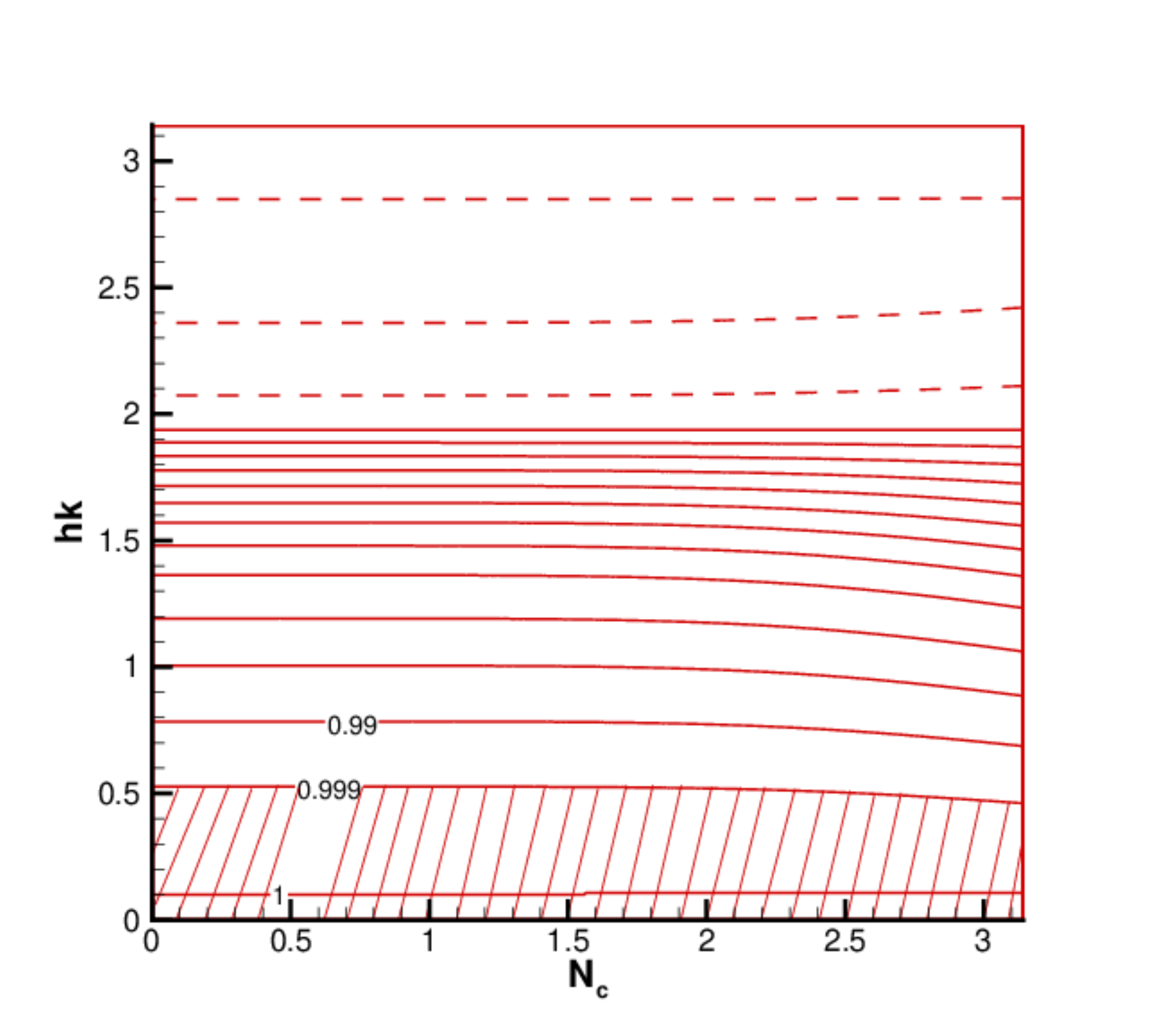,width=0.7\linewidth}
 (f)
\end{minipage}            \hspace{-2.5mm}
\begin{minipage}[b]{.6\linewidth}   \
\centering\psfig{file=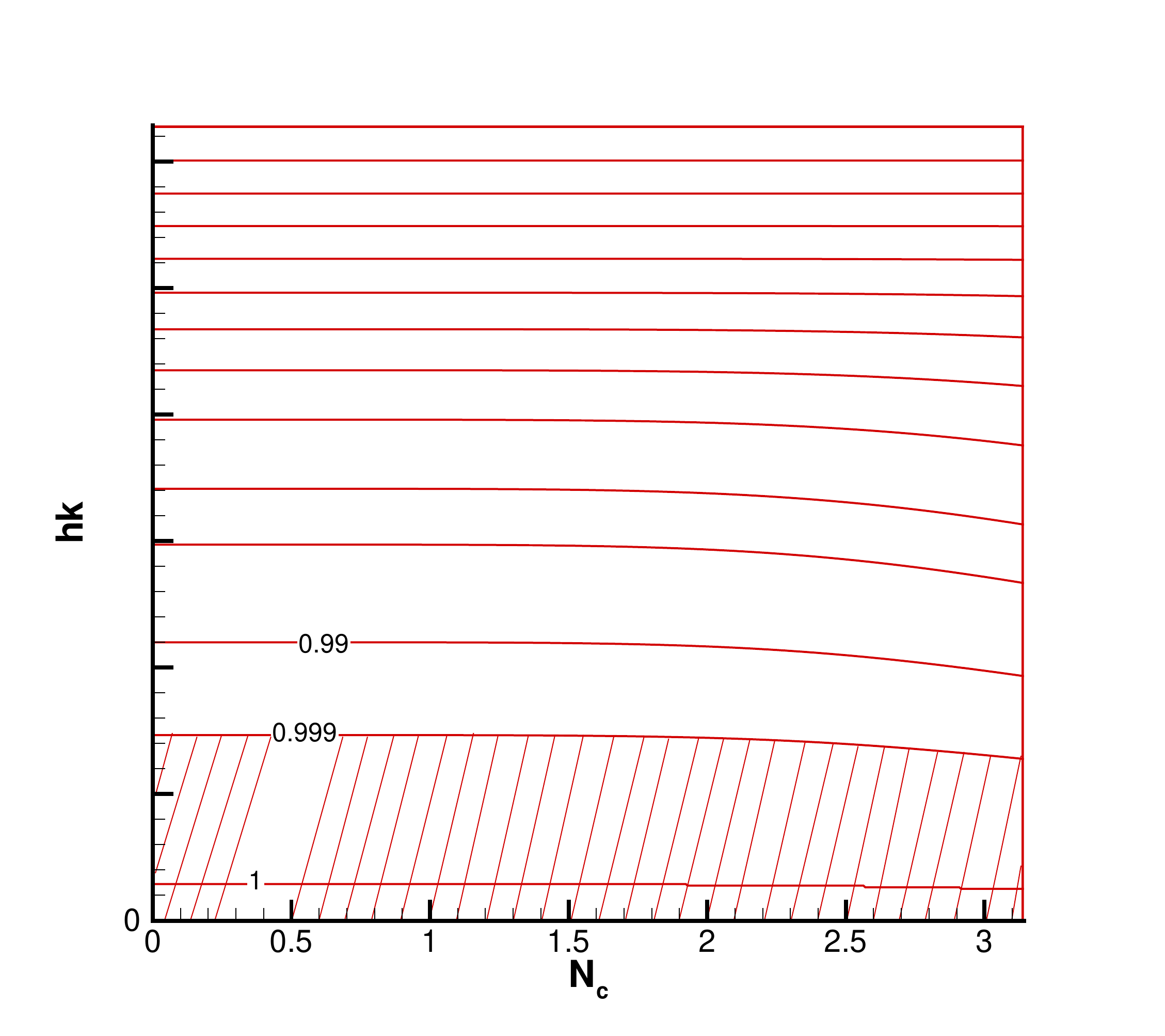,width=0.7\linewidth}
 (g)
\end{minipage}            \hspace{-2.5mm}
\begin{minipage}[b]{.6\linewidth}
\centering\psfig{file=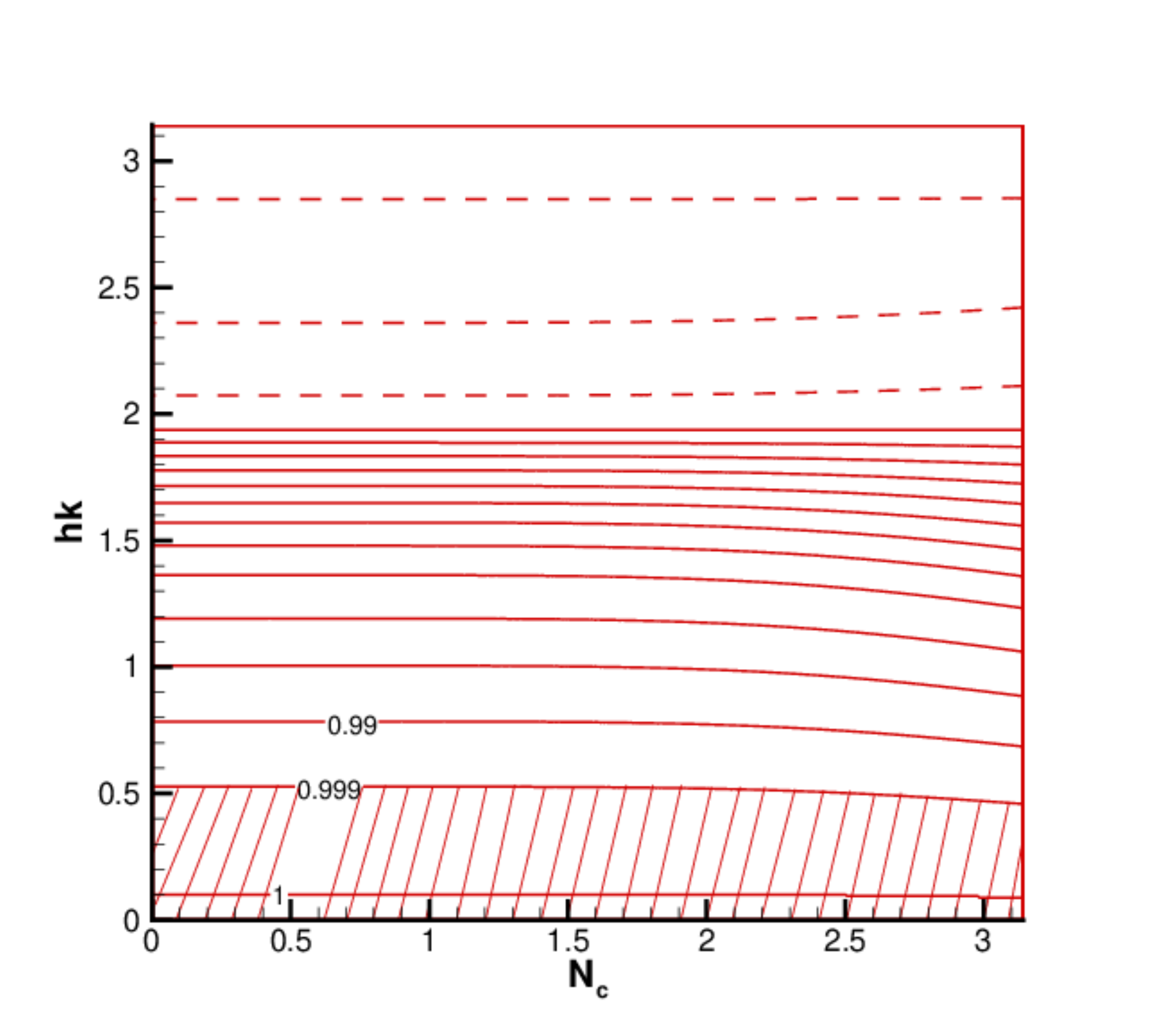,width=0.7\linewidth}
 (h)
\end{minipage}            \hspace{-2.5mm}
\begin{center}
\caption{{\sl Contours of normalized numerical phase velocity (left) and group velocity (right) for indicated schemes plotted in the $(N_c, hk)$ plane at mid-node when explicit CD6 scheme is used for spatial discretization: (a)-(b) S3A $(\alpha=0)$, (c)-(d) S3B $(\alpha=4)$, (e)-(f) S3C $(\alpha=16)$, (g)-(h) IRK36, S3D $(\alpha\rightarrow\infty)$.}}
\label{fig:3S_expl}
\end{center}
\end{figure}

Significantly analogous to each weight parameter chosen there is a narrow region of computation where dispersive error free computation is very much possible. The same was also noticed for the two stage schemes. This region being very sensitive and amenable to round-off error the same has not been attempted in this study. Similar to earlier discussion on two stage schemes, it is noticed that beyond some threshold $\sigma$ values S2B and S2C carries lesser dispersive error than IRK36.

Analysis of spectral properties at the central node of various three stage schemes in conjunction with implicit Lele scheme \cite{lel_92} and explicit central CD6 is carried out in figures \ref{fig:3S_impl} and \ref{fig:3S_expl}  respectively. Compared to two stage schemes a significant improvement in the region of $\pm0.1\%$ tolerance in the scaled values of numerical phase and group velocities can be seen for schemes derived by taking unity weight function. For all other set of schemes enlargement of efficient domain for computation can be seen at higher $N_c$ values pointing to better efficiency for three stage methods for stiff problems. Contours of $(v_{pN}/c)_j$ and $(v_{gN}/c)_j$ are found to be fairly similar for S3B, S3C and S3D. Carrying forward similarity from two stage here also S3D and IRK36 are found to be indistinguishable from each other for both implicit and explicit spatial discretization and as such have not been plotted separately.

\section{Methodology summarized}
Main steps followed in deriving various schemes can be summarized below.
\begin{enumerate}
\item Start with an appropriate stage number $R=2$ or $3$ in this work.
\item Introduce additional conditions such that $|G_{N}(\sigma)|=1$.
\item Look to impose highest possible order of accuracy such that free parameters are available.
\item Decide on appropriate weight kernel by keeping in mind
\begin{itemize}
\item wavenumber,
\item time step. 
\end{itemize}
\item Formulate phase error in $L^2$-norm over a suitable wavenumber space with the above chosen weight kernel.
\item Solve minimization problem.
\item Solve the resulting system of equations.
\end{enumerate}

\section{Numerical Examples}
\subsection{Problem 1: Solution of test equation}
To investigate proposed schemes vis-a-vis high accuracy schemes we begin by considering numerical solution of test equation which can be written as a following system of IVP
\begin{eqnarray}\label{P0_1}
\dot{\bm U}=\bm \Lambda\bm U,\;\;\; \bm U(0)=(1, 0)^T;\;\;\;\text{ with}\;\;\;\bm \Lambda=\left(
                    \begin{array}{cc}
                      0 & -\lambda \\
                      \lambda & 0 \\
                    \end{array}
                  \right).
\end{eqnarray}
Although computations were carried out using all schemes tabulated in previous sections we present, in table \ref{table P0_1}, error in solutions obtained using first two schemes from each category. 
\begin{table}[h!]
\caption{Problem 1: Comparison of absolute error at point of maxima and rate of convergence.}
\centering
\begin{footnotesize}
\begin{tabular}{c c c c c c c c c c}
\hline
Scheme&$\Delta t$=0.008&Rate&$\Delta t$=0.016&Rate&$\Delta t$=0.032&Rate&$\Delta t$=0.064&Rate&$\Delta t$=0.128 \\
\hline
S2A1 &5.8332e-4	&1.99	& 2.3249e-3 & 1.98  & 9.1646e-3 & 1.91 	& 3.4520e-2 & 1.61 	& 1.0541e-1\\
S2A2 &5.8332e-4	&1.99	& 2.3249e-3 & 1.98  & 9.1646e-3 & 1.91 	& 3.4520e-2 & 1.61 	& 1.0541e-1\\
\hline
S2B1 &2.9180e-5 &1.93	& 1.1130e-4 & 1.69  & 3.5894e-4 & --   	& \textbf{8.6288e-5} & 7.80 	& \textbf{1.9298e-2}\\
S2B2 &2.9181e-5	&1.93	& 1.1130e-4 & 1.69  & 3.5894e-4 & --   	& \textbf{8.6303e-5} & 7.80 	& \textbf{1.9298e-2}\\
\hline
S2C1 &7.0119e-6	&1.70	& 2.2767e-5 & --    & \textbf{7.0591e-6} & 7.51 	& 1.2852e-3 & 4.24 	& 2.4218e-2\\
S2C2 &7.0119e-6	&1.70	& 2.2767e-5 & --    & \textbf{7.0591e-6} & 7.51 	& 1.2852e-3 & 4.24 	& 2.4218e-2\\
\hline
IRK24&\textbf{4.3674e-7}	&4.00	& \textbf{6.9798e-6} & 3.99  & 1.1117e-4 & 3.97  & 1.7460e-3 & 3.89  & 2.5869e-2\\
S2D2 &\textbf{4.3674e-7}	&4.00	& \textbf{6.9799e-6} & 3.99  & 1.1117e-4 & 3.97  & 1.7460e-3 & 3.89  & 2.5869e-2\\
\hline
S3A1 &2.8049e-8	&4.00	& 4.4741e-7 & 3.98  & 7.0706e-6 & 3.93 	& 1.0757e-4 & 3.68 	& 1.3825e-3\\
S3A2 &2.8050e-8	&4.00	& 4.4742e-7 & 3.98  & 7.0707e-6 & 3.93 	& 1.0757e-4 & 3.68 	& 1.3825e-3\\
\hline
S3B1 &2.1162e-9	&3.96	& 3.2881e-8 & 3.82  & 4.6301e-7 & 2.89 	& 3.4260e-6 & 5.76 	& \textbf{1.8525e-4}\\
S3B2 &2.1172e-9	&3.96	& 3.2885e-8 & 3.82  & 4.6302e-7 & 2.89 	& 3.4261e-6 & 5.76 	& \textbf{1.8525e-4}\\
\hline
S3C1 &5.1955e-10&3.81	& 7.3011e-9 & 2.92  & \textbf{5.5070e-8} & 5.77 	& \textbf{3.0021e-6} & 6.55 	& 2.8187e-4\\
S3C2 &5.1653e-10&3.82	& 7.2889e-9 & 2.92  & \textbf{5.5023e-8} & 5.77 	& \textbf{3.0022e-6} & 6.55 	& 2.8188e-4\\
\hline 	
S3D1 &\textbf{1.7094e-11}&6.20	& \textbf{1.2610e-9} & 6.01  & 8.1420e-8 & 5.98 	& 5.1526e-6 & 5.93 	& 3.1420e-4\\
IRK36&\textbf{2.5843e-11}&5.62	& \textbf{1.2752e-9} & 6.00  & 8.1475e-8 & 5.98 	& 5.1528e-6 & 5.93 	& 3.1420e-4\\
\hline 			
\end{tabular}
\end{footnotesize}	
\label{table P0_1}
\end{table}
Computations are done employing as many as five different time steps $\Delta t=$ 0.008, 0.016, 0.032, 0.064 and 0.128. Least error at each time step for both two and three stage schemes have been highlighted. For this problem time period is $T=2\pi/10$. Smaller values of $\Delta t=0.008$ and 0.016 lead to $\sigma=$ 0.08 and 0.16 respectively. Thus it is seen that IRK24 (S2D1) and IRK36 (S3D2) produces least error in accordance with the analysis carried out in earlier sections. Further all schemes, by and large, conform to their theoretical rate of convergence when computed with such small time steps. Reported error shows a different trend with $\Delta t$ increasing which results in increase of $\sigma$ values. At $\Delta t=$ 0.032, S2C class of schemes report least error. A look at the figure \ref{fig:2S_phase}(b) reveals that dispersion error is close to minima at the corresponding $\sigma$ value for S2C class of schemes. As $\Delta t$ increases to 0.064, $\sigma$ value approaches region of minimum dispersion error for S2B class of schemes and is accordingly reported in table \ref{table P0_1}. For $\Delta t=$ 0.128, S2B retains its advantage in terms of dispersion error as predicted by figure \ref{fig:2S_phase}(b). A similar pattern can be noticed for three stage schemes. At $\Delta t=$ 0.032, S3C is supposed to carry least dispersion error and the numerical computation confirm least overall error. At $\Delta t=$ 0.064, S2C is marginally better than S2B in terms of dispersion preservation in figure \ref{fig:3S_phase}(b) and our numerical results amply demonstrate the same. For $\Delta t=$ 0.128, S2B is certainly the best of the lot. This problem duly demonstrate the importance of dispersion and dissipation relation preservation beyond rate of convergence for computations done using relatively bigger time step. S2A and S3A class of schemes carry better dispersion property only at high $\sigma$ values and carries certain shortcoming in terms of phase and group velocity as seen earlier. This will be further discussed in subsequent sections before commenting on their efficiency. In the table \ref{table P0_1}, it is seen that for cases ehere the error is less, separate schemes under a particular category may produce distinctive absolute errors. This may be attributed to the changes in underlying implicit system, because of alteration of coefficients, solved using open-source tool Lis \cite{kot_has_kaj_08}. 
 
\subsection{Problem 2: Periodic test }
In this test a periodical initial value problem represented by second order ODE
\begin{eqnarray}\label{P1_1}
\ddot{u}=-k^2u+(k^2-\omega^2)\sin(\omega t),\;\;\; t\geq 0,\;\;\; u(0)=u_0,\;\;\; \dot{u}(0)=\bar{u}_0
\end{eqnarray}
is numerically solved \cite{naz_moh_cha_15}. This problem admit analytical solution 
\begin{eqnarray}\label{P1_2}
u(t)=u_0\cos(kt)+\frac{(\bar{u}_0-\omega)\sin(kt)}{k}+\sin(\omega t).
\end{eqnarray}
For $u_0=0$ and $\bar{u}_0=\omega$, the same reduces to $u(t)=\sin(\omega t).$

\begin{table}[h!]
\caption{Problem 2: Absolute error at point of maxima and order of convergence.}
\centering
\begin{tabular}{c c c c c c c c}
\hline
Scheme&$\Delta t$=0.016&Order&$\Delta t$=0.032&Order&$\Delta t$=0.064&Order &$\Delta t$=0.128 \\
\hline
 					S2A & 1.6303e-3 & 1.95  & 6.3090e-3 & 1.81 & 2.2161e-2 & 1.31 & 5.4964e-2\\
\hline
					S2B & 7.5929e-5 & 1.48  & 2.1227e-4 & 1.40 & 5.5976e-4 & 5.36 & 2.2917e-2\\
\hline
 					S2C & 1.3538e-5 & 1.36  & 3.4885e-5 & 5.44 & 1.5168e-3 & 4.13 & 2.6517e-2\\
\hline
 				   IRK24& 7.4294e-6 & 3.99  & 1.1798e-4 & 3.96 & 1.8392e-3 & 3.91 & 2.7734e-2\\
\hline
					S3A & 5.0629e-7 & 3.97  & 7.9574e-6 & 3.89 & 1.1828e-4 & 3.51 & 1.3524e-3\\
\hline
					S3B & 3.6518e-8 & 3.72  & 4.8201e-7 & 1.44 & 1.3072e-6 & 8.08 & 3.5405e-4\\
\hline
					S3C & 7.6048e-9 & 1.53  & 2.1913e-8 & 8.07 & 5.8897e-6 & 6.28 & 4.5921e-4\\
\hline 	
				   IRK36& 2.0722e-9 & 5.99  & 1.3199e-7 & 5.97 & 8.2968e-6 & 5.90 & 4.9438e-4\\
\hline 						
\end{tabular}
\label{table P1_1}
\end{table}

For numerical computation two frequencies $\omega$ and $k$ are maintained at distinct values $10$ and $15$ respectively. This problem is studied to investigate efficiency of various schemes conceptualized in this work for higher order ODE in combination with two different frequencies. Four different time stepping $\Delta t$ = 0.016, 0.032, 0.064 and 0.128 are used to compute absolute error at time $t = 0.768$ near the point of maxima of $u(t)$. We compute using first scheme of each category of two stage methods and second of the each class of three stage methods. This helps us to compare results obtained using newly developed schemes with those of high accuracy Gauss-Legendre schemes. Results obtained are arranged in Table \ref{table P1_1}. It is seen that none of the newly proposed low-dissipation, low-dispersion methods from diverse categories are able to match accuracy of IRK24 and IRK36 among two and three stage methods respectively at $\Delta t =$ 0.016. But with increase in step size newly developed schemes show much better results. We see a similar trend as found earlier and the newly developed schemes show much better accuracy at their respective regions of low dispersion. This example demonstrate productivity of the newly developed schemes for varied IVP.

\subsection{Problem 3: Convection of a wave packet}
Next one dimensional linear convection equation
\begin{eqnarray}\label{P2_1}
u_t+u_x=0
\end{eqnarray}
is considered. We take initial profile as
\begin{eqnarray}\label{P2_2}
u(x,0)=e^{-\frac{(x-x_m)^2}{b}} \cos(k(x-x_m)).
\end{eqnarray}
Following Rajpoot et al. \cite{raj_sen_dut_10} the problem is treated as periodic in the domain $0\le x\le30$ with $b=2$ and constant grid spacing $h=0.01$. Spatial discretization is carried out using sixth order five point Lele scheme \cite{lel_92}. A fixed $x_m=5$ is taken for all reported computations. We use this problem to judge dispersive effects of various implicit schemes introduced in this work. In the process effect of phase error reduction on overall accuracy of the numerical schemes will be highlighted. For two stage methods various optimized schemes have been tabulated in each category but here we compute with first scheme of each category viz. S2A1, S2B1, S2C1 and S2D1. Since S2D1 represent IRK24 we take S2D2 as representative of S2D class of schemes. For three stage we compute with S3A1, S3B1, S3C1 and S3D1. S3D2 representing three stage Gauss-Legendre scheme (IRK36), is additionally used because of its own significance. Comparison of solution obtained using S2D1 (IRK24) and S2D2 as well as S3D1 and S3D2 will also reveal to some extent how error varies as schemes are switched within a category.

We choose $k=4$ and compute upto $t=20.0$ for $N_c=4.0$, $7.5$, $15.0$ and $20.0$. $L^2$-norm error between numerical and exact solutions is reported in table \ref{table P2_1}. At $N_c=4.0$, IRK24 and IRK36 are found to be most accurate among two and three stage schemes respectively. This is in consonance with our dispersive error analysis as $N_c=4.0$ correspond to $\sigma=0.16$ for the chosen wave number. From figures \ref{fig:2S_phase} and \ref{fig:3S_phase} it is clear that at $\sigma=0.16$ IRK24 and IRK36 indeed carry least dispersion error. As $N_c$ is increased $7.5$, $\sigma$ equals to $3.0$ and S2C and S3C are found to report least error in their respective classes, pointing towards a close relation between dispersion error and accuracy for the class of $A$-stable methods discussed here. As $\sigma$ increases further with $N_c =$ $15.0$ and $20.0$, S2B and S3B are seen to be most accurate. Despite our best effort with various solvers IRK36 does not converge at $N_c=20.0$. As temporal truncation error increases rapidly with bigger CFL number all methods starts diverging and hence S2A and S3A class of methods are unable to realise their potential at high $\sigma$ values. Hence it must be said that S2A and S3A may not ideal at any $N_c$ value.

\begin{table}[h!]
\caption{Problem 3: $L^2$-norm error between numerical and exact solutions with $k=4$ at $t=20$.}
\vspace{0.2 cm}
\centering
\begin{tabular}{l c c c c}
\hline
Scheme	&$N_c=4.0$&$N_c=7.5$&$N_c=15.0$&$N_c=20.0$\\
  \hline\\
  S2A  	&5.1067e-3	&1.7668e-2	&6.5196e-2	&1.0463e-1\\
  S2B 	&2.4235e-4	&6.5231e-4	&1.2859e-3	&7.3554e-3\\
  S2C 	&4.7893e-5	&8.0395e-5	&3.3978e-3	&1.1517e-2\\
  S2D	&2.1471e-5	&2.7127e-4	&4.2438e-3	&1.2967e-2\\
  IRK24 &1.9385e-5	&2.7131e-4	&4.2448e-3	&1.2968e-2\\
  \hline\\
  S3A  	&2.5353e-6	&1.7414e-5	&2.5751e-4	&7.4928e-4\\
  S3B 	&2.4130e-6	&2.9260e-6	&7.0320e-6	&4.0531e-5\\
  S3C 	&2.4122e-6	&2.6141e-6	&1.2926e-5	&8.4007e-5\\
  S3D	&2.4215e-6	&2.6424e-6	&1.9307e-5	&9.9518e-5\\
  IRK36 &2.4086e-6	&8.6216e-6	&4.7991e-5	&--\\
  \hline
\end{tabular}
\label{table P2_1}
\end{table}

\begin{figure}[!ht]
\begin{minipage}[b]{.6\linewidth}
\centering\psfig{file=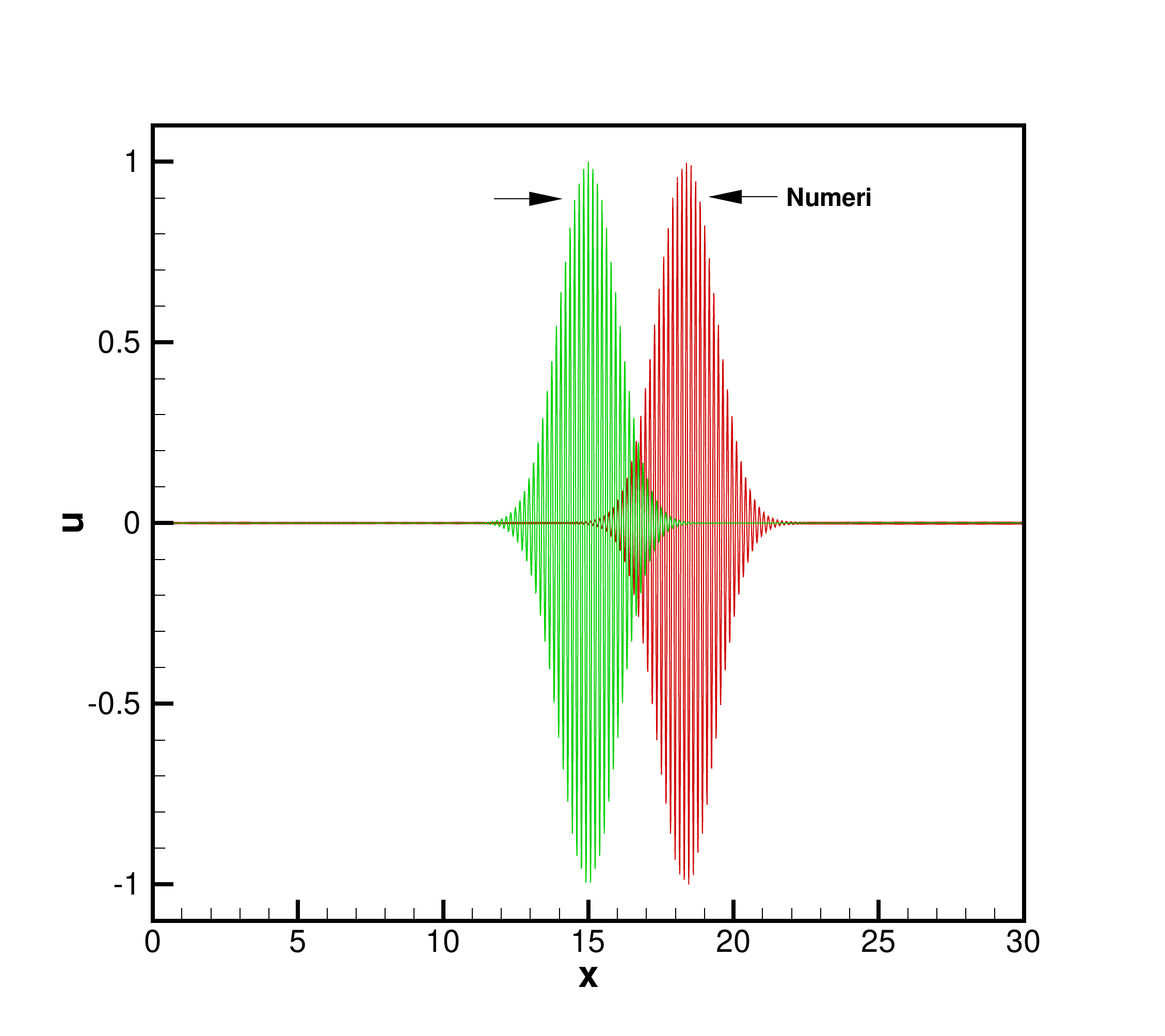,width=0.95\linewidth}
 \\(a)
\end{minipage}            \hspace{-2.5mm}
\begin{minipage}[b]{.6\linewidth}
\centering\psfig{file=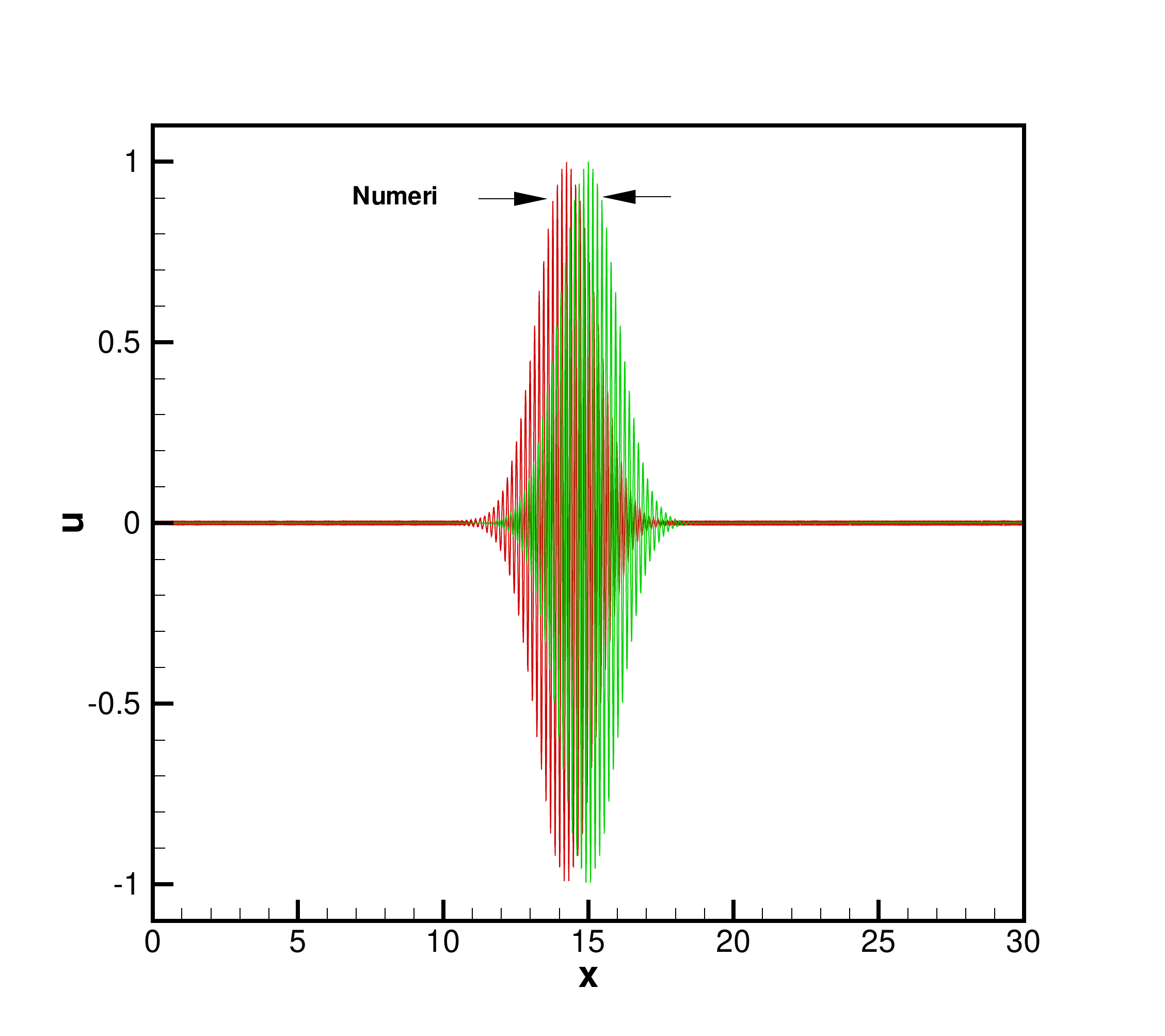,width=0.95\linewidth}
 \\(b)
\end{minipage}            \hspace{-2.5mm}
\begin{minipage}[b]{.6\linewidth}   \
\centering\psfig{file=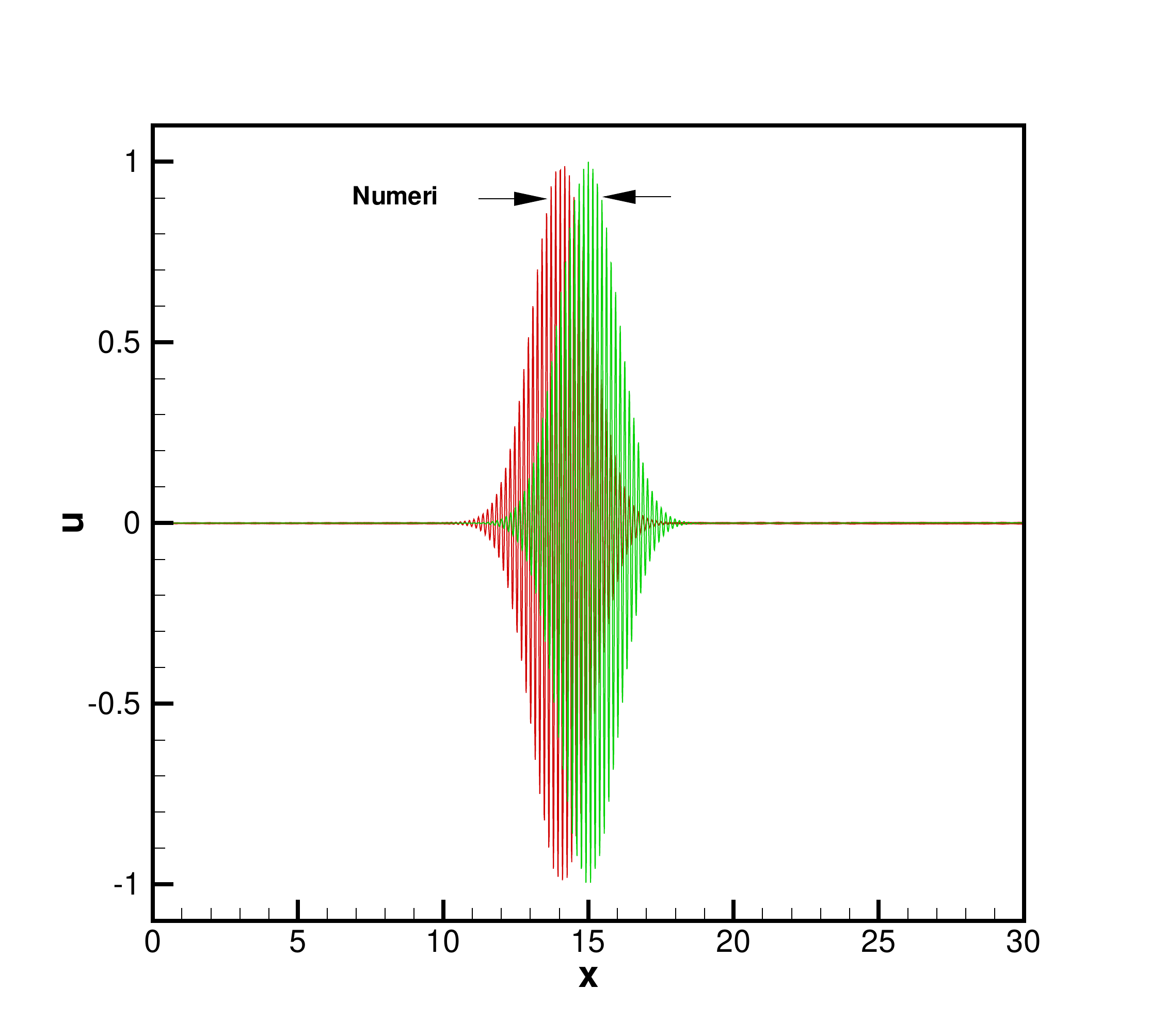,width=0.95\linewidth}
 \\(c)
\end{minipage}            \hspace{-2.5mm}
\begin{minipage}[b]{.6\linewidth}
\centering\psfig{file=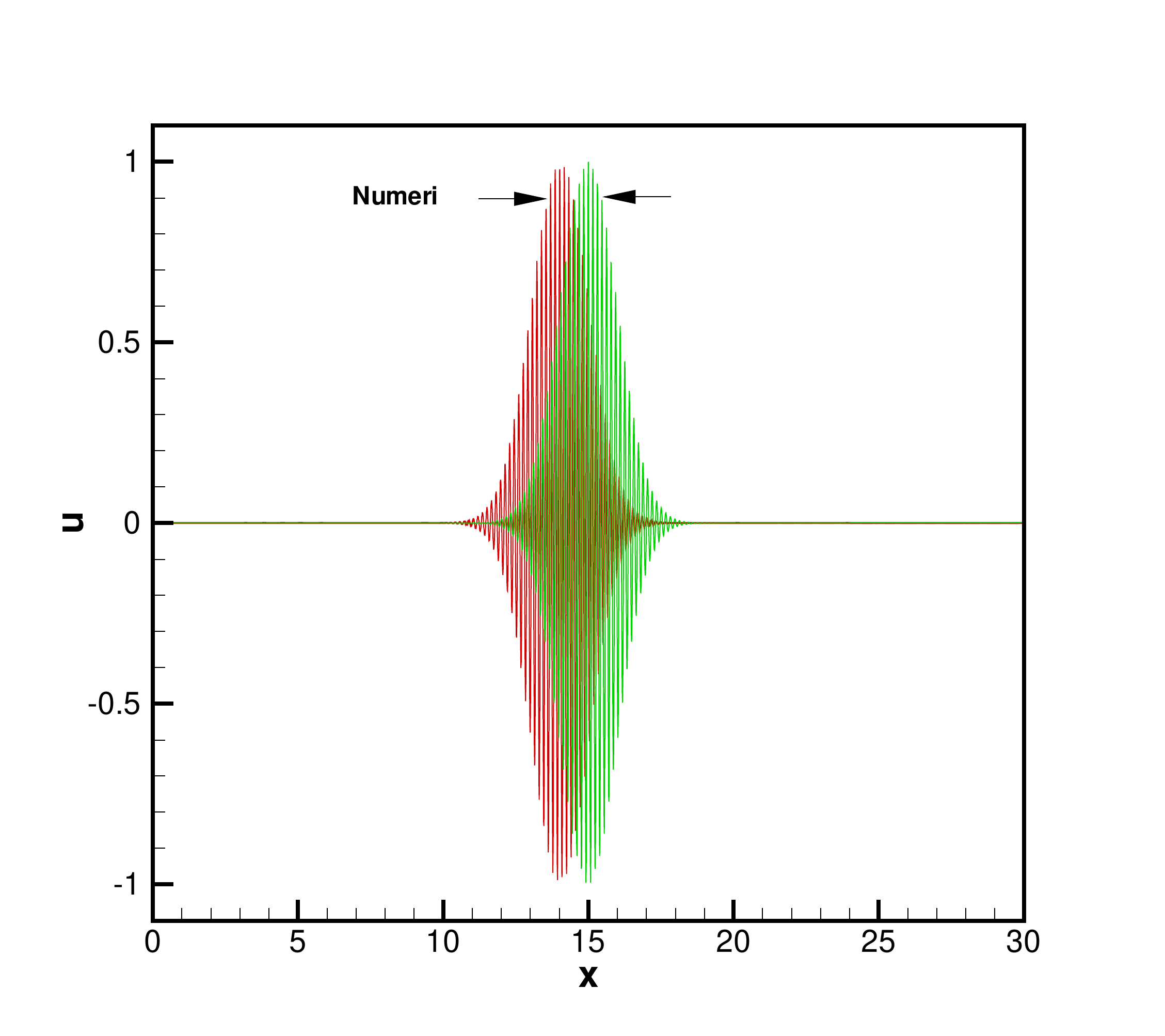,width=0.95\linewidth}
 \\(d)
\end{minipage}            \hspace{-2.5mm}
\begin{center}
\caption{{\sl Problem 3: Propagation of a computed wave-packet with $k=40$, $N_c=2.6$ at $t=130$ using various two stage methods (a) S2A, (b) S2B, (c) S2C, and (d) IRK24.} }
\label{fig:P2_1}
\end{center}
\end{figure}
\begin{figure}[!ht]
\begin{minipage}[b]{.6\linewidth}
\centering\psfig{file=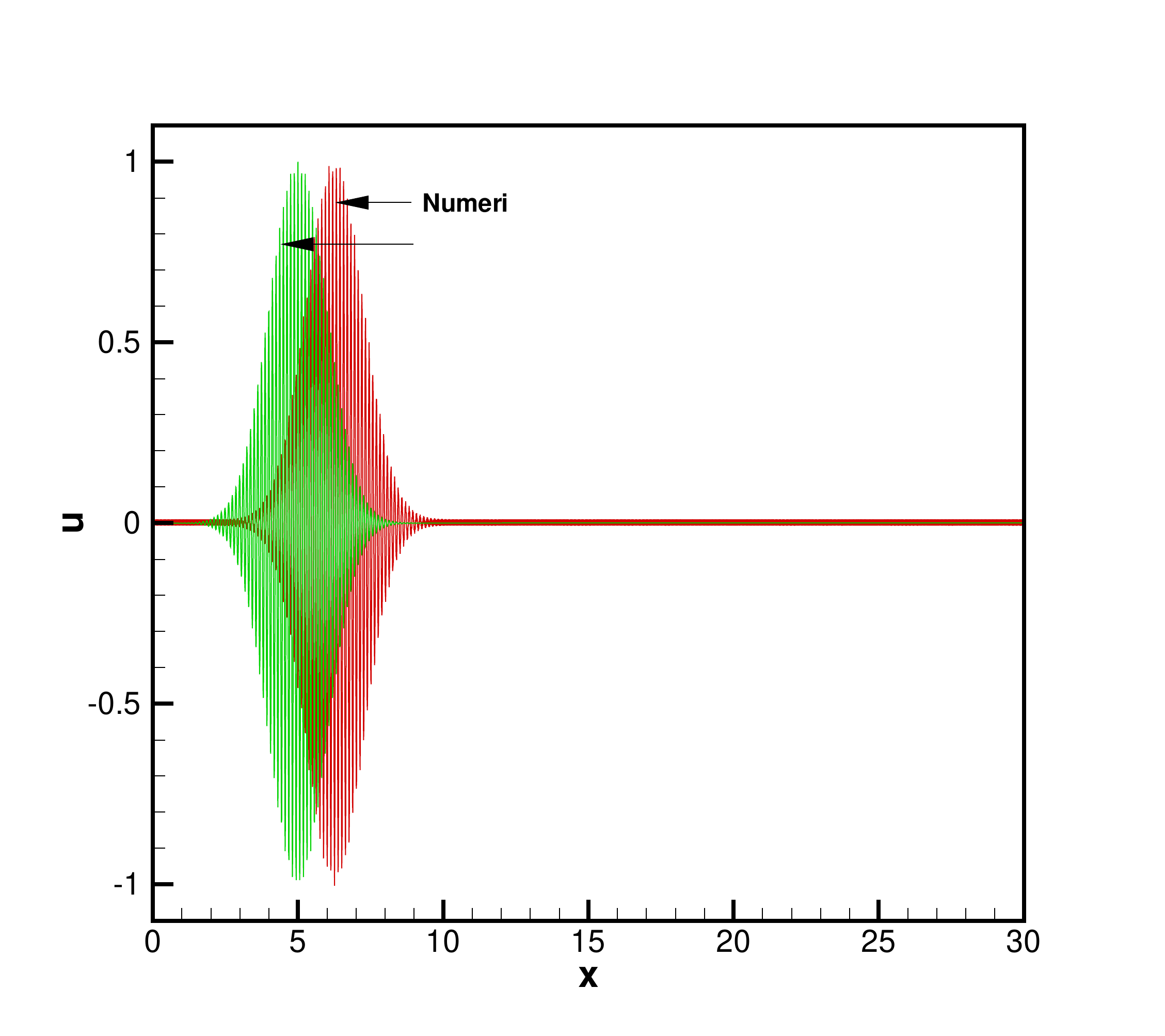,width=0.95\linewidth}
 \\(a)
\end{minipage}            \hspace{-2.5mm}
\begin{minipage}[b]{.6\linewidth}
\centering\psfig{file=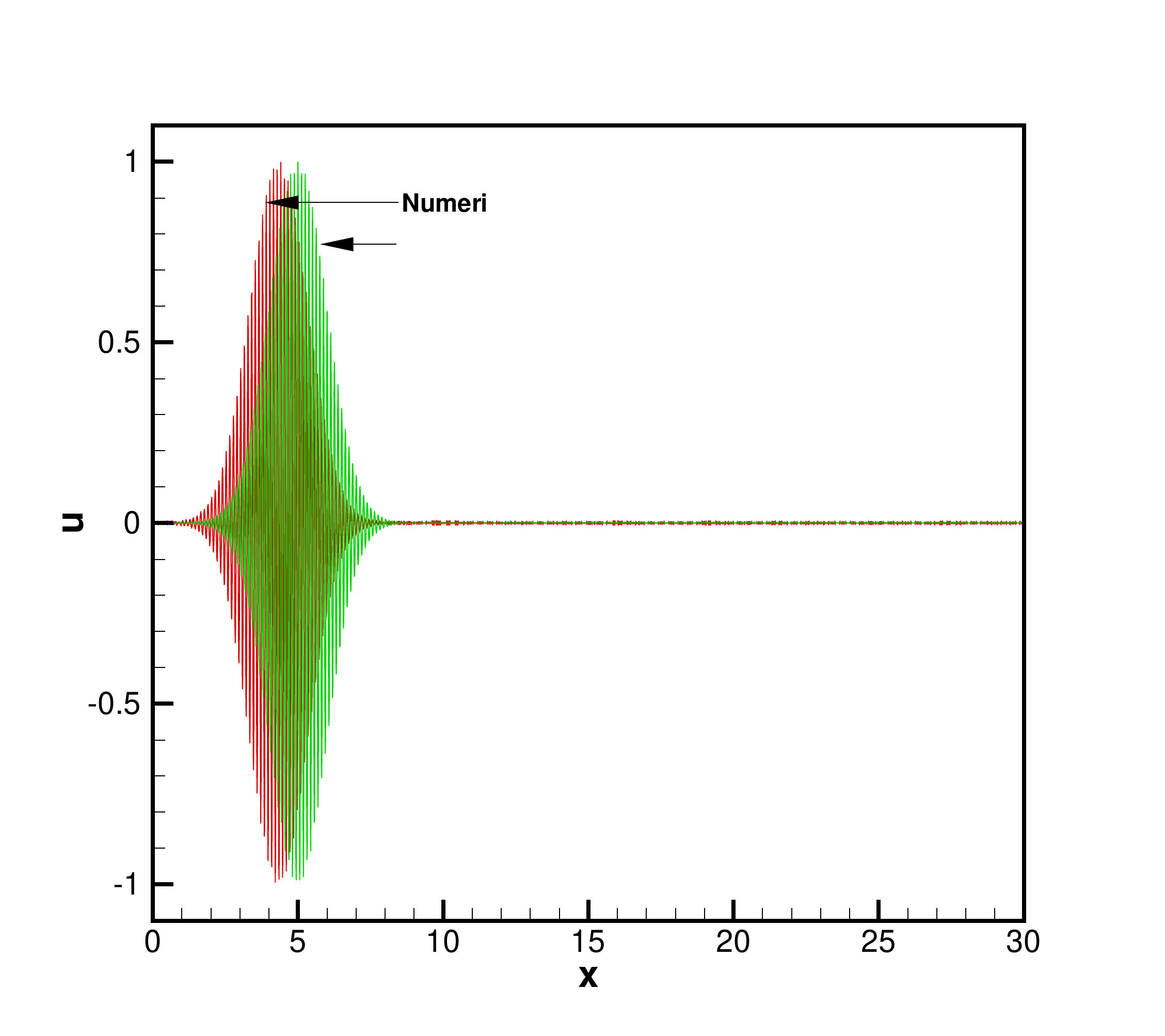,width=0.95\linewidth}
 \\(b)
\end{minipage}            \hspace{-2.5mm}
\begin{minipage}[b]{.6\linewidth}   \
\centering\psfig{file=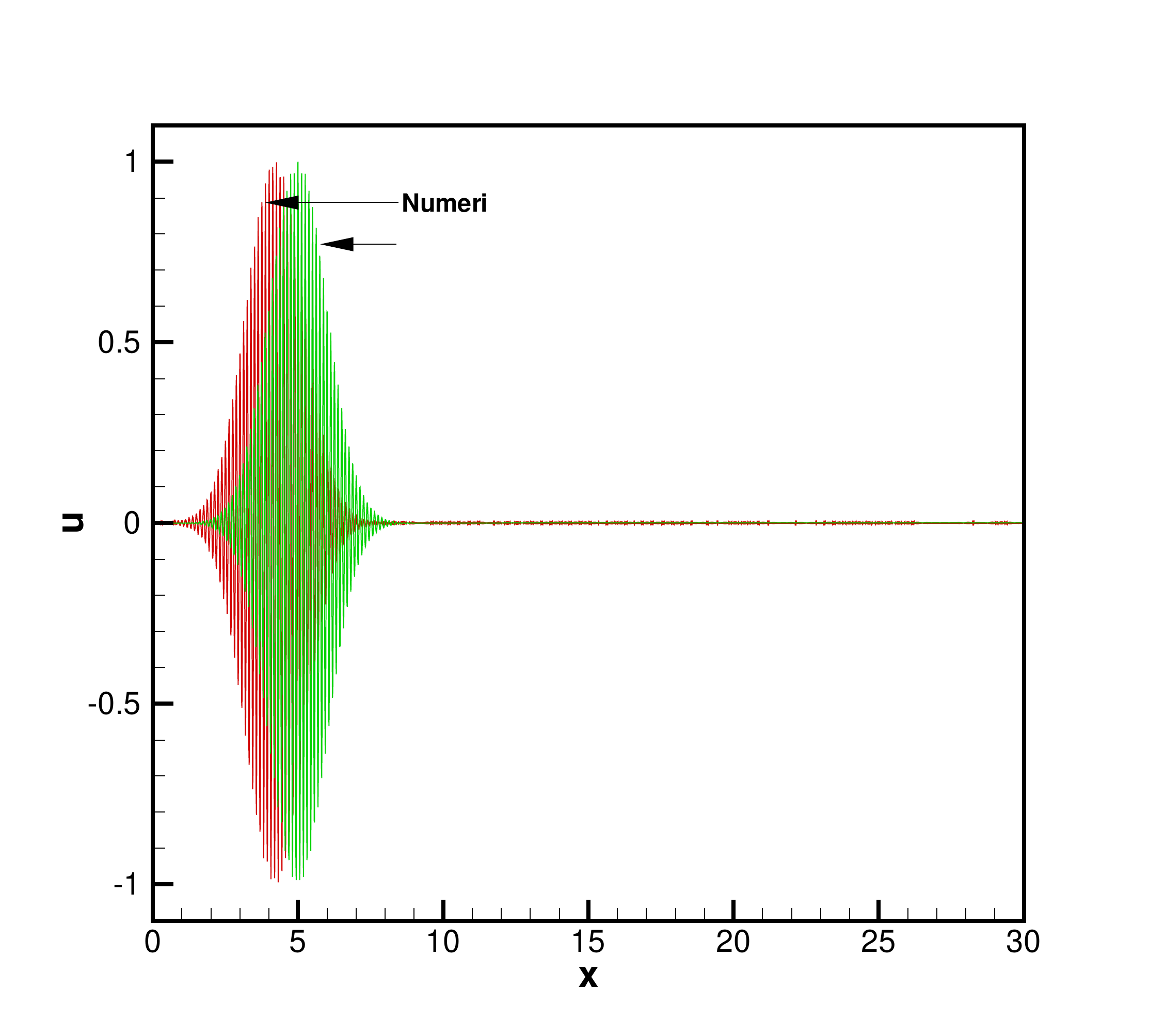,width=0.95\linewidth}
 \\(c)
\end{minipage}            \hspace{-2.5mm}
\begin{minipage}[b]{.6\linewidth}
\centering\psfig{file=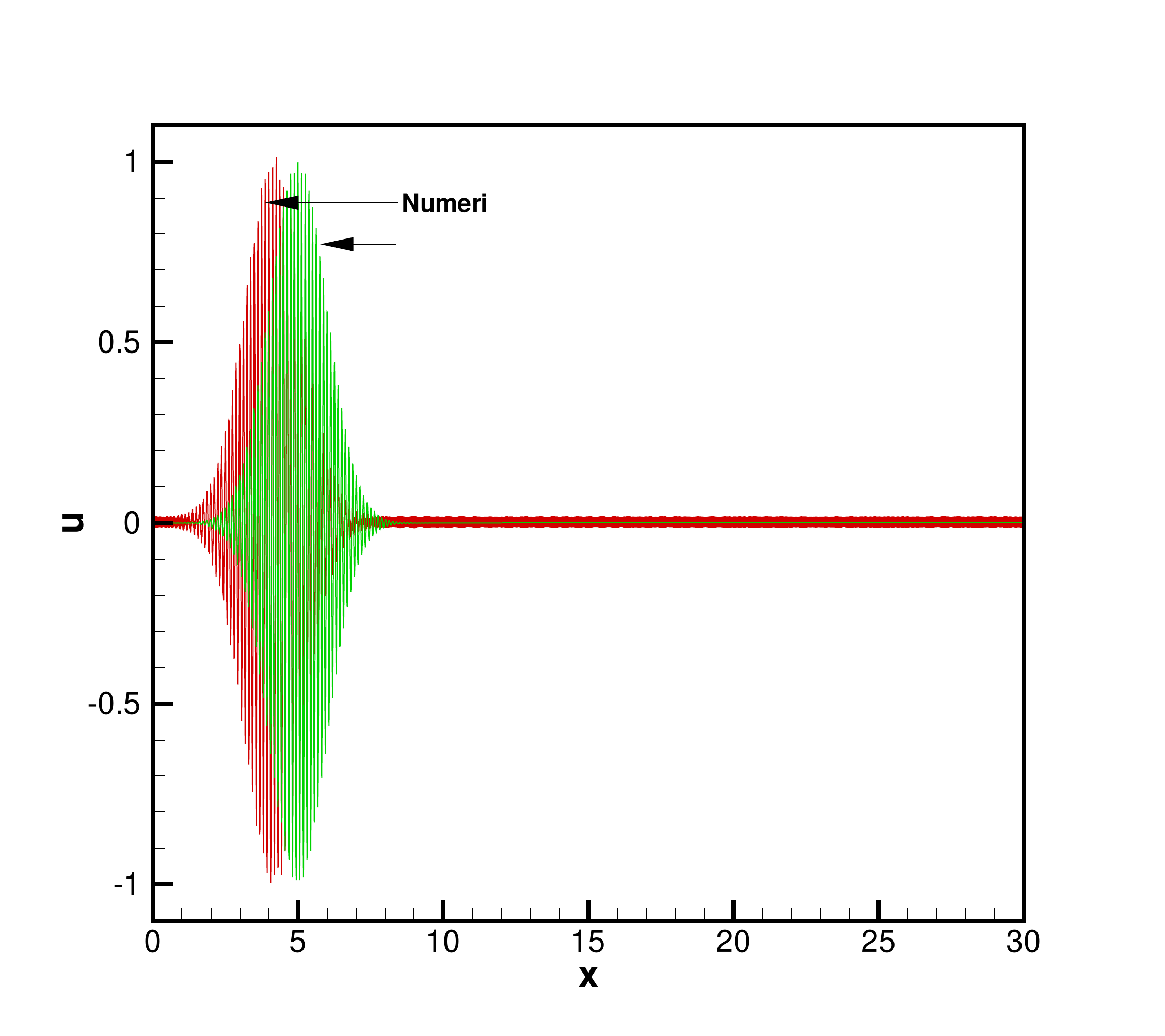,width=0.95\linewidth}
 \\(d)
\end{minipage}            \hspace{-2.5mm}
\begin{center}
\caption{{\sl Problem 3: Propagation of a computed wave-packet with $k=50$, $N_c=3.0$ at $t=1080$ using various three stage methods (a) S3A, (b) S3B, (c) S3C, and (d) IRK36.} }
\label{fig:P2_2}
\end{center}
\end{figure}
Further test are conducted with optimized schemes to reveal long time behaviour of dispersion error accumulation. For two stage scheme we chose a high wave number $k=40$ and compute with $N_c=2.6$. Such a combination takes us to the edge and beyond of the hatched region of $(N_c, hk)$ plane as indicated in figure \ref{fig:2S_impl}. Computed solution are compared with exact solution at time $t=130$ in figure
\ref{fig:P2_1}. For S2A we see that the numerical propagation speed is higher compared to actual propagation speed. This is in consonance with our analysis of figure \ref{fig:2S_impl}(a)-(b). For S2B, S2C and IRK24 there is a close resemblance between numerical and analytical group velocity. Least difference is predicted and reported in figure \ref{fig:P2_1}(b) for S2B. With three stage schemes we take a higher wave number $k=50$ and compute with a bigger $N_c=3.0$ value for longer duration of time $t=1080$. Three stage methods possesses overall better dispersion characteristics vis-a-vis two stage schemes and we notice closer match between analytical and numerical solution in figure \ref{fig:P2_2}. Nevertheless overall pattern of different types of schemes remain unaffected and marginally better efficiency of S3B set of methods for long time computation as can be seen in figure \ref{fig:P2_2}(b).

\subsection{Problem 4: Convection of a combination of waves}
Here we consider the linear convection of a combination of two waves of wavenumbers $2\pi k_1$ and $2\pi k_2$ given initially as
\begin{eqnarray}\label{P2_3}
u(x,0)=e^{-\frac{(x-x_m)^2}{b}} [\cos(2\pi k_1(x-x_m))+\cos(2\pi k_2(x-x_m))].
\end{eqnarray}
We take $x_m=90$, $b=400$, $k_1=0.125$ and $k_2=0.0625$ \cite{naj_mon_13} and compute solutions for various values of $N_c$ values ranging from 1.0 to 3.0. For spatial discretization sixth order five point Lele scheme \cite{lel_92} is used with $\Delta x=0.5$. To avoid problems of boundary reflection a large computational domain with Dirichlet boundary condition is used and solution is computed upto time $t=300.0$.

\begin{figure}[!ht]
\begin{minipage}[b]{.6\linewidth}
\centering\psfig{file=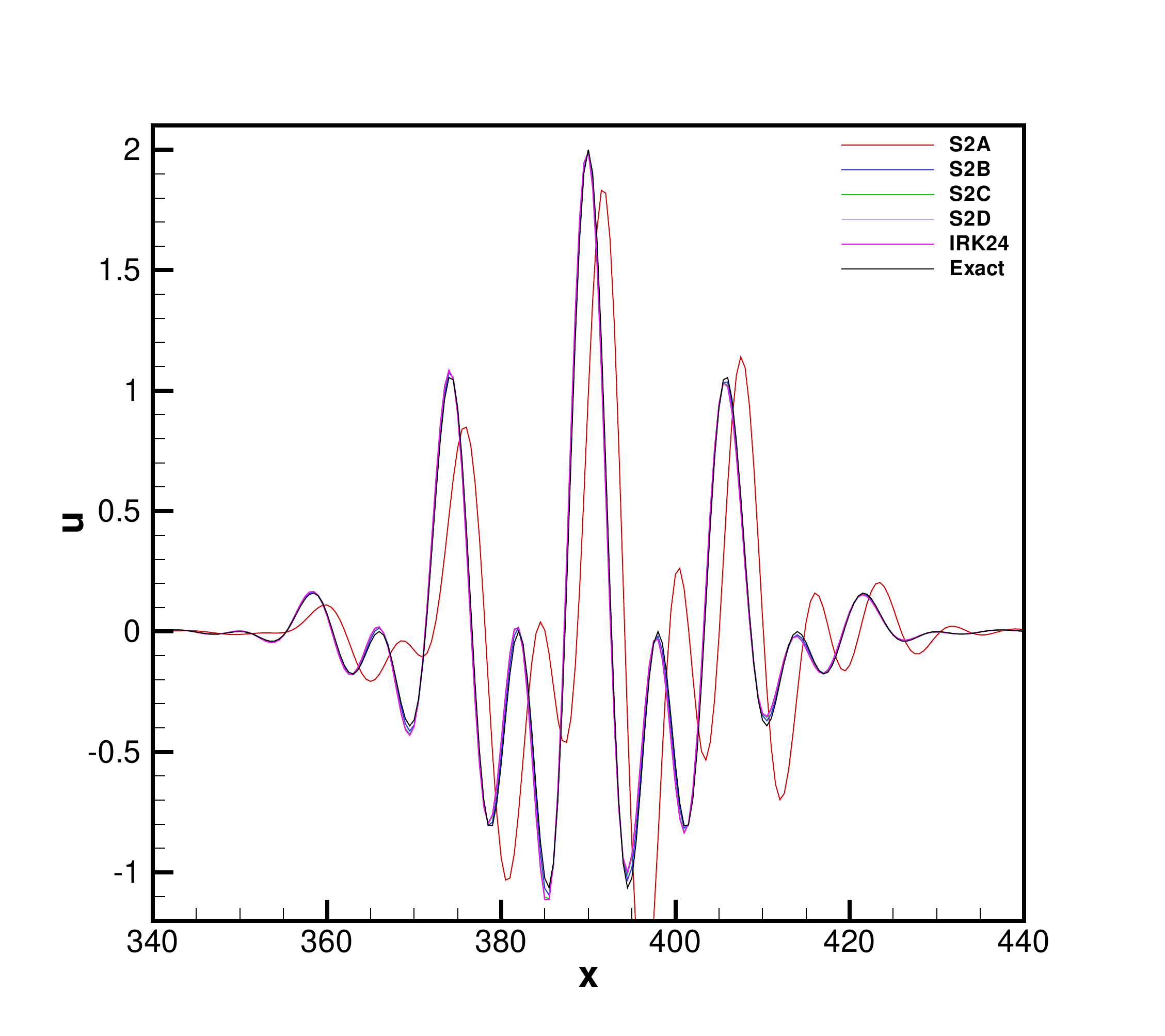,width=0.95\linewidth}
 \\(a)
\end{minipage}            \hspace{-2.5mm}
\begin{minipage}[b]{.6\linewidth}
\centering\psfig{file=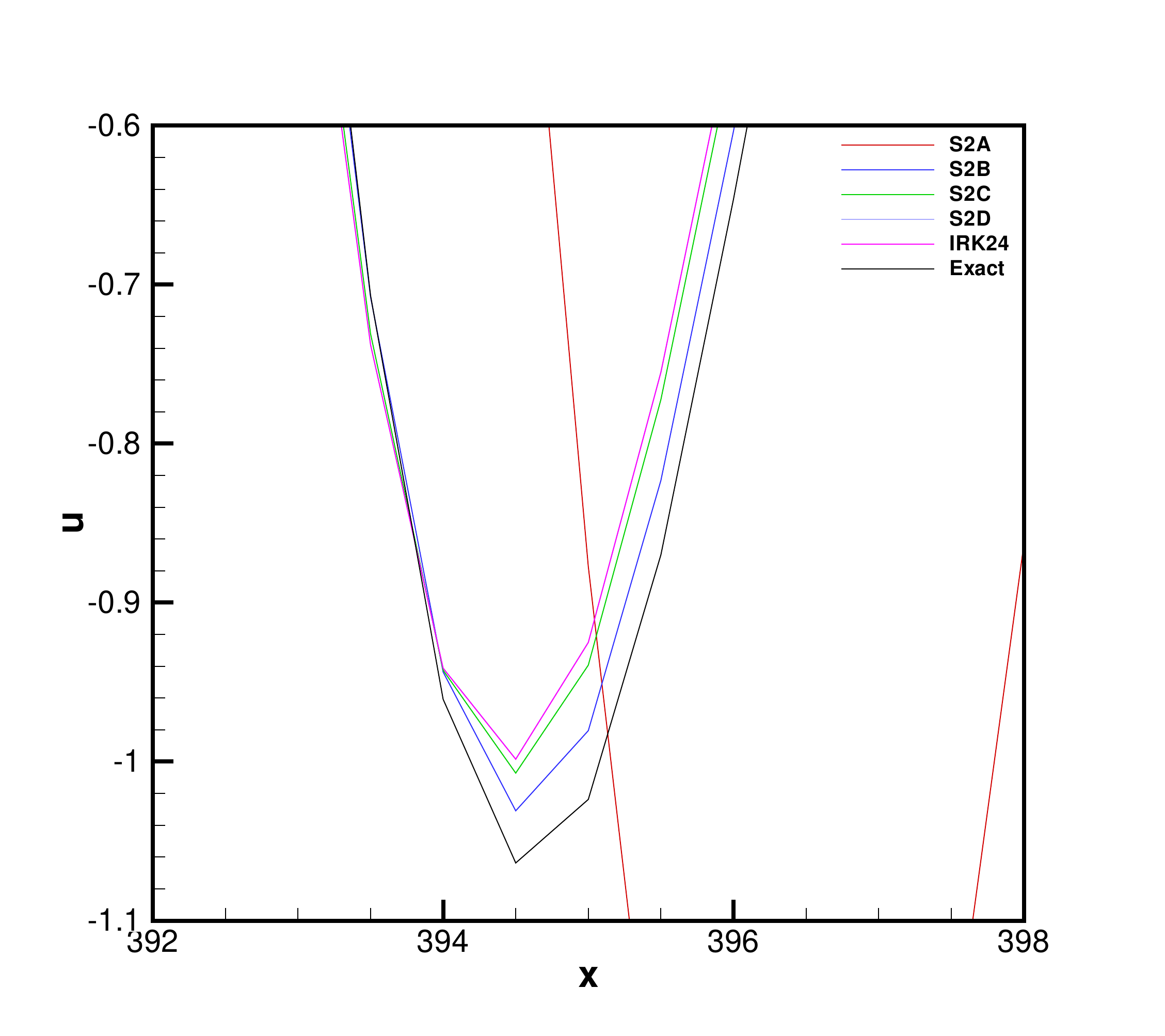,width=0.95\linewidth}
 \\(b)
\end{minipage}            \hspace{-2.5mm}
\begin{minipage}[b]{.6\linewidth}   \
\centering\psfig{file=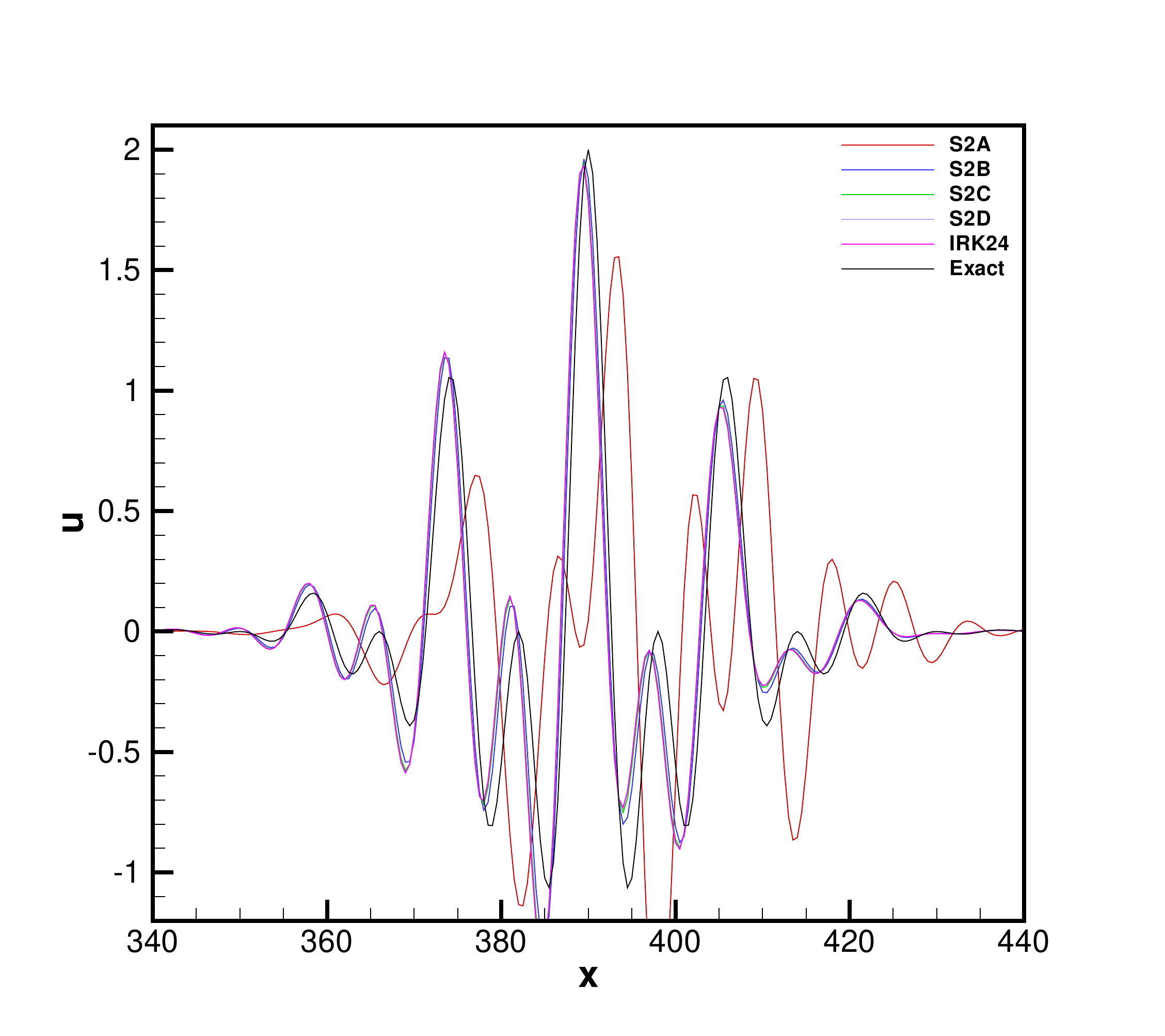,width=0.95\linewidth}
 \\(c)
\end{minipage}            \hspace{-2.5mm}
\begin{minipage}[b]{.6\linewidth}
\centering\psfig{file=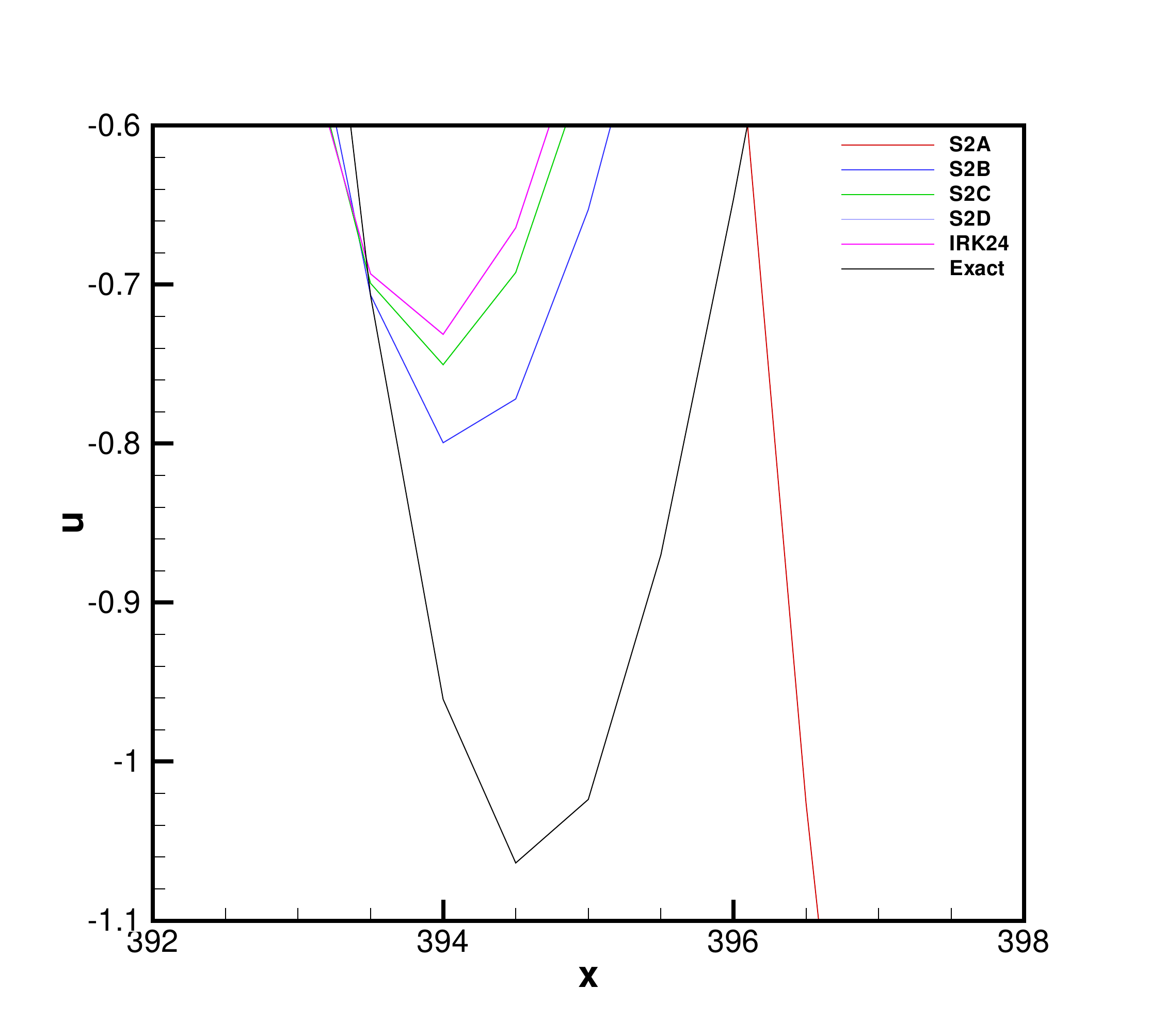,width=0.95\linewidth}
 \\(d)
\end{minipage}            \hspace{-2.5mm}
\begin{center}
\caption{{\sl Problem 4: Numerical solution at $t=300$ using two stage methods for (a) $N_c=2.0$, (b) $N_c=2.0$ close view, (c) $N_c=3.0$, and (d) $N_c=3.0$ close view.} }
\label{fig:P3_1}
\end{center}
\end{figure}

\begin{figure}[!ht]
\begin{minipage}[b]{.6\linewidth}
\centering\psfig{file=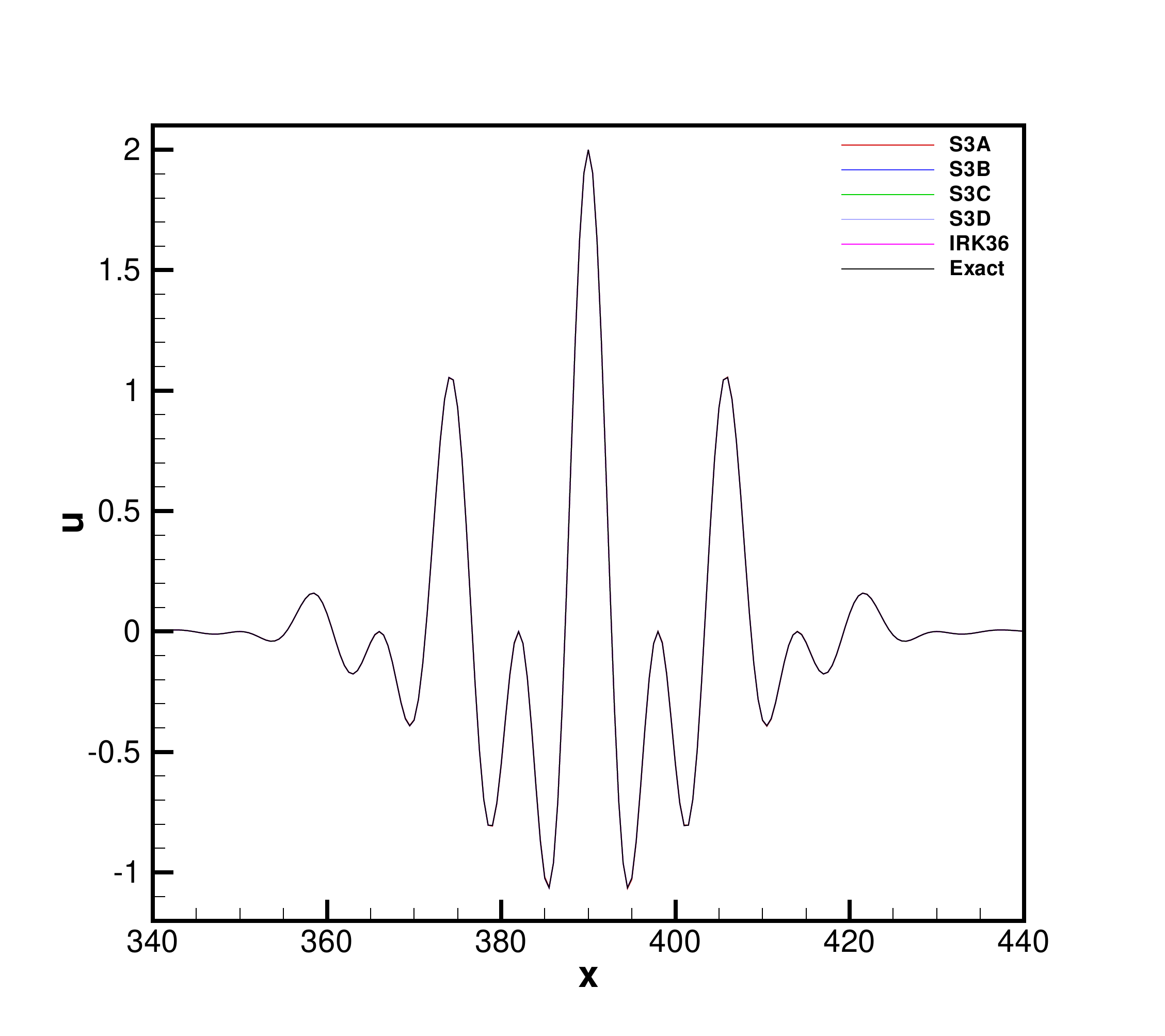,width=0.95\linewidth}
 \\(a)
\end{minipage}            \hspace{-2.5mm}
\begin{minipage}[b]{.6\linewidth}
\centering\psfig{file=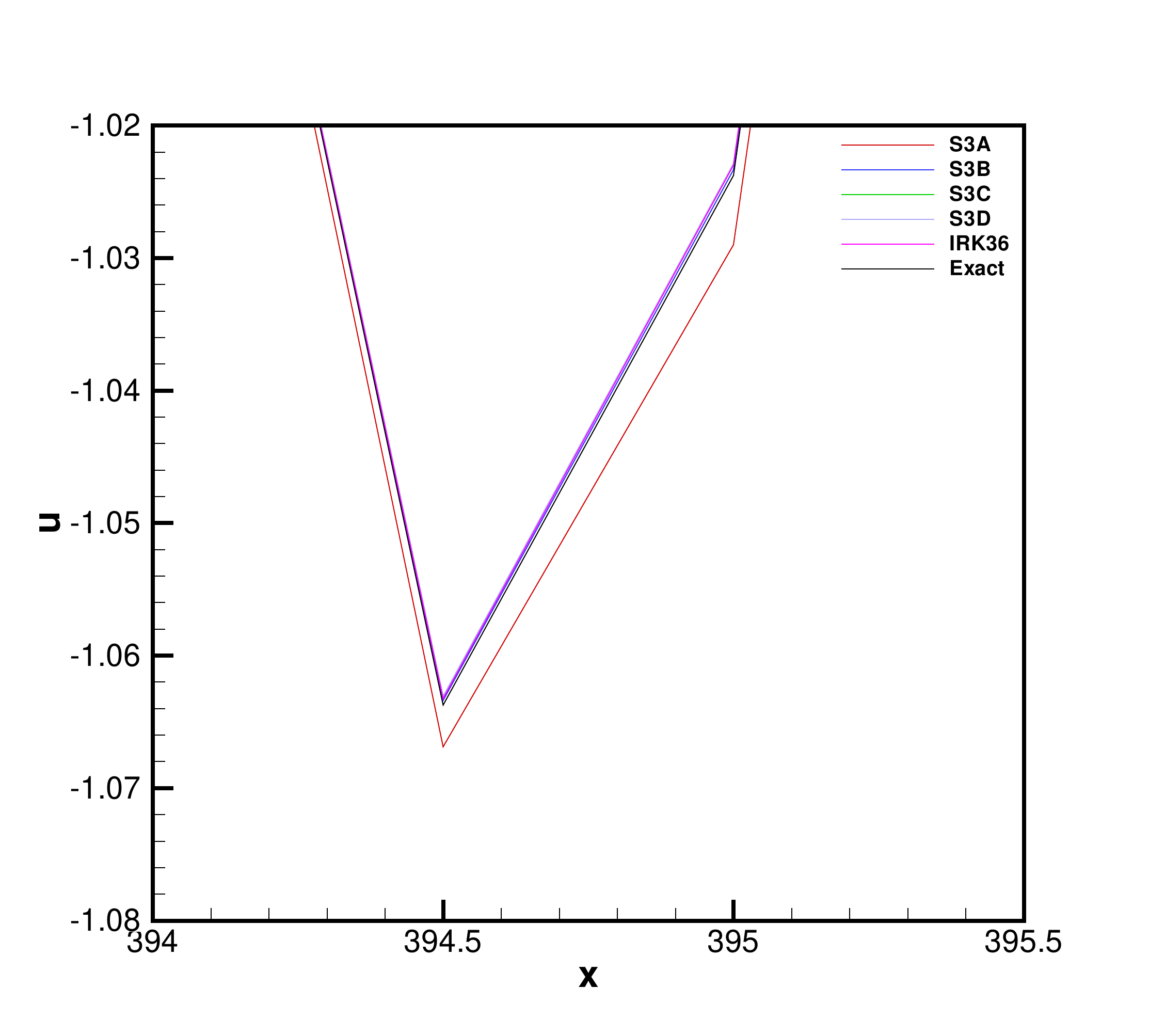,width=0.95\linewidth}
 \\(b)
\end{minipage}            \hspace{-2.5mm}
\begin{minipage}[b]{.6\linewidth}   \
\centering\psfig{file=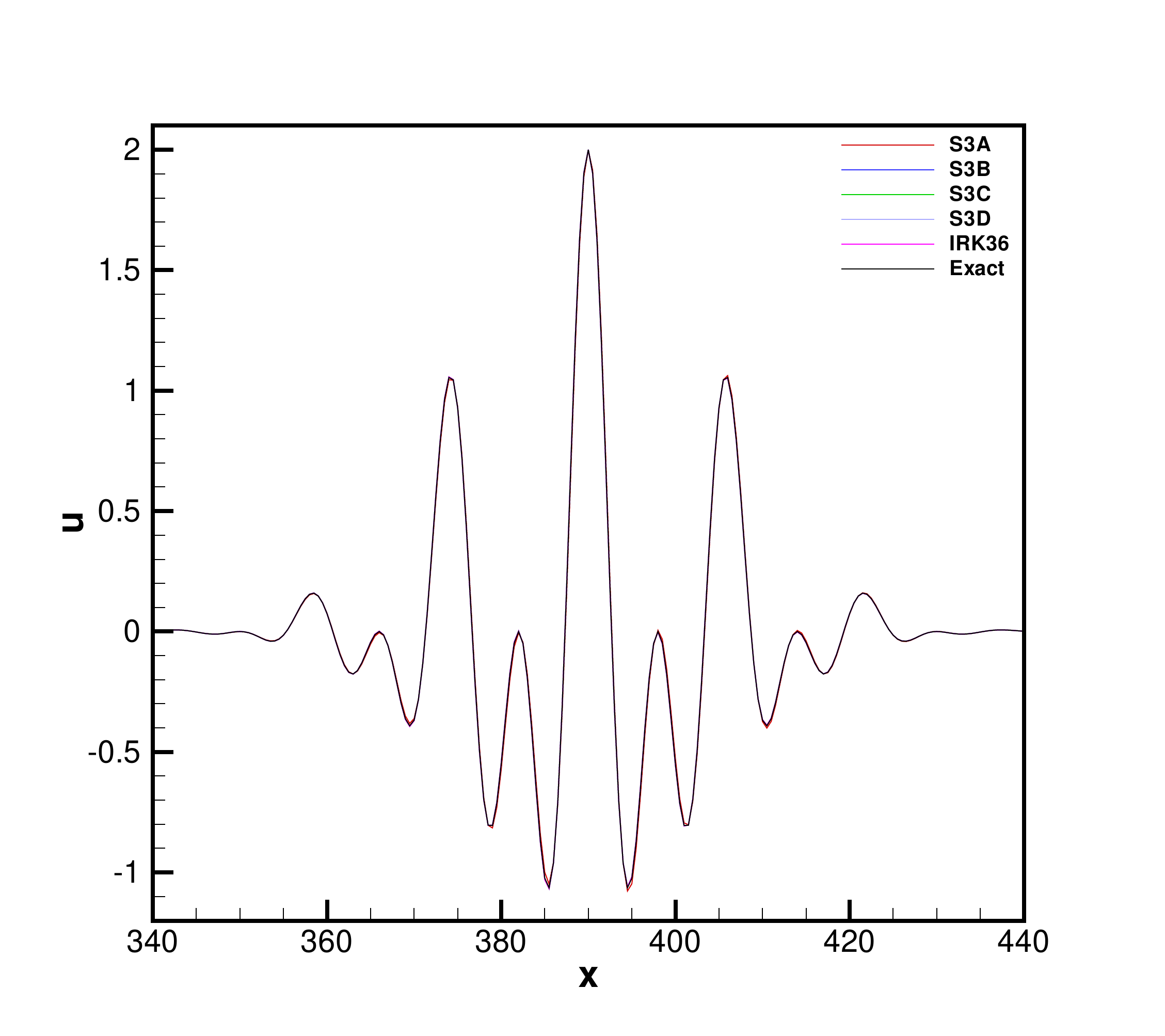,width=0.95\linewidth}
 \\(c)
\end{minipage}            \hspace{-2.5mm}
\begin{minipage}[b]{.6\linewidth}
\centering\psfig{file=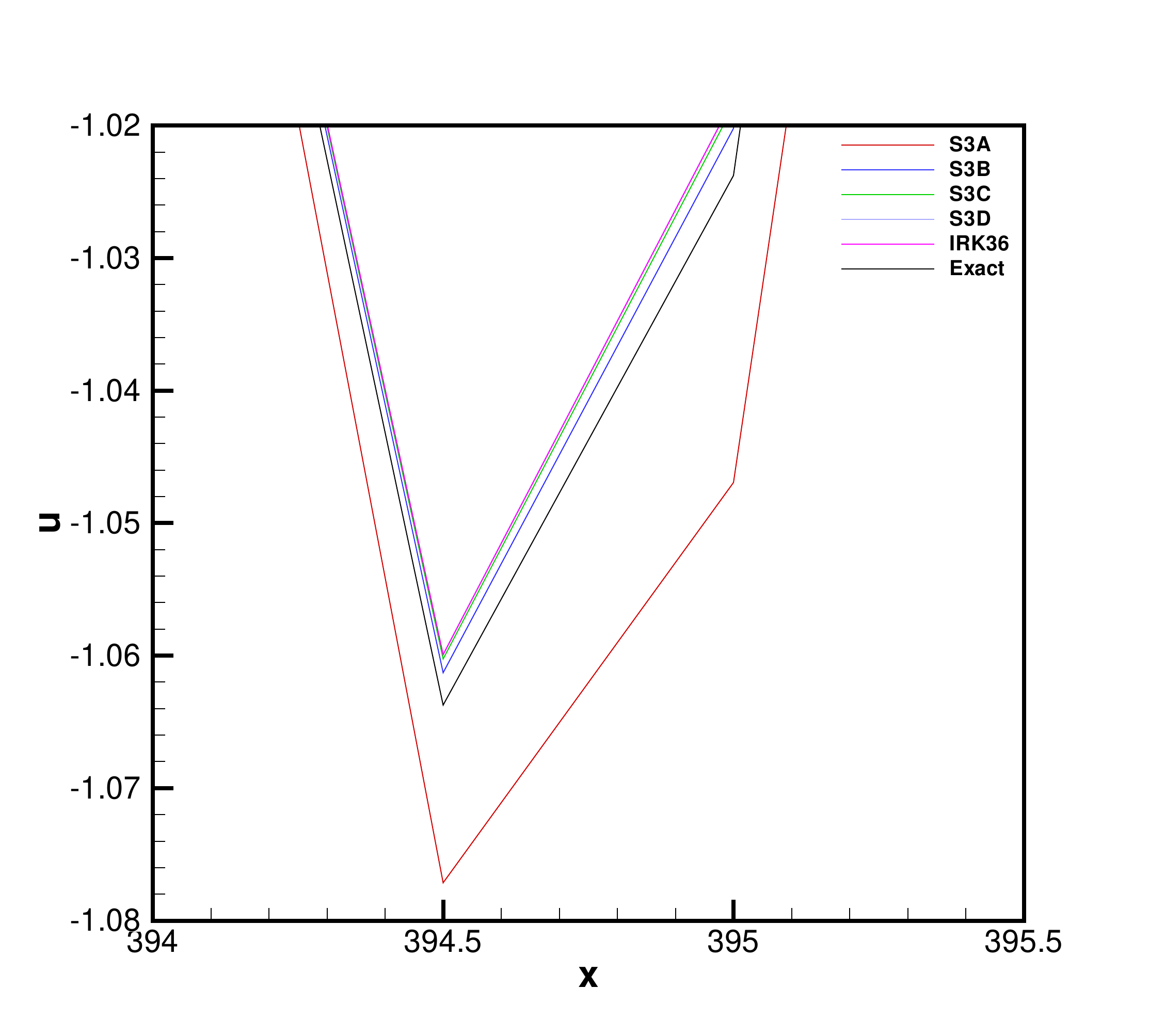,width=0.95\linewidth}
 \\(d)
\end{minipage}            \hspace{-2.5mm}
\begin{center}
\caption{{\sl Problem 4: Numerical solution at $t=300$ using three stage methods for (a) $N_c=2.0$, (b) $N_c=2.0$ close view, (c) $N_c=3.0$, and (d) $N_c=3.0$ close view.} }
\label{fig:P3_2}
\end{center}
\end{figure}

In Table \ref{table P3_1} and Figures \ref{fig:P3_1} - \ref{fig:P3_2} we present our results. From the Table \ref{table P3_1}, it is seen that the results obtained by S2D1 and IRK24, both of which carry identical dispersion characteristics but different order of accuracy, produces matching error for all CFL numbers. Similarly it can be concluded that results obtained using S3D1 and IRK36 are almost indistinguishable. This clearly exemplify importance of dispersion relation preserving schemes.

\begin{table}[h!]
\caption{Problem 4: $L^2$-norm error between numerical and exact solutions at $t=300$.}
\vspace{0.2 cm}
\centering
\begin{tabular}{l H H H c H H H c c c c}
\hline
Scheme	&$N_c=0.25$&$N_c=0.5$&$N_c=0.840$&$N_c=1.0$&$N_c=1.184$&$N_c=1.323$&$N_c=1.669$&$N_c=1.5$&$N_c=2.0$&$N_c=2.5$&$N_c=3.0$\\
  \hline
  S2A  	&3.9933e-3	&1.6013e-2	&4.4563e-2	&6.2300e-2  &8.5640e-2  &1.0489e-1 &1.5625e-1  &1.3067e-1 &2.0416e-1 &5.5814e-0 &2.8966e-1\\
  S2B 	&1.4392e-4	&6.7429e-4	&1.6356e-3	&1.9867e-3  &2.1421e-3  &1.9968e-3 &1.5372e-3  &1.4744e-3 &6.9368e-3 &2.6593e-2 &6.4304e-2\\
  S2C 	&2.5380e-5	&6.7428e-5	&1.4625e-4	&5.0056e-4  &1.3428e-3  &2.3858e-3 &7.1749e-3  &4.3809e-3 &1.5856e-2 &4.0243e-2 &8.2631e-2\\
  S2D	&8.1041e-5	&1.4033e-4	&6.8066e-4	&1.2875e-3  &2.4546e-3  &3.7717e-3 &9.3547e-3  &6.1481e-3 &1.8939e-2 &4.4865e-2 &8.8742e-2\\
  IRK24 &6.7484e-5	&1.4968e-4	&6.8013e-4	&1.2886e-3  &2.4550e-3  &3.7729e-3 &9.3609e-3  &6.1548e-3 &1.8936e-2 &4.4852e-2 &8.8735e-2\\
  \hline
  S3A  	&--	&6.7428e-5	&--	&1.0312e-4  &--  &-- &--  &3.0621e-4 &1.0610e-3 &2.5398e-3 &4.8808e-3\\
  S3B 	&--	&7.1708e-5	&--	&6.7960e-5  &--  &-- &--  &6.3170e-5 &6.4138e-5 &1.9789e-4 &6.4027e-4\\
  S3C 	&--	&7.1963e-5	&--	&7.2375e-5  &--  &-- &--  &8.5769e-5 &1.2907e-4 &3.6135e-4 &9.7353e-4\\
  S3D	&--	&7.2044e-5	&--	&7.3853e-5  &--  &-- &--  &9.3368e-5 &1.5114e-4 &4.1657e-4 &1.0857e-3\\
  IRK36 &--	&7.5025e-5	&--	&7.4967e-5  &--  &-- &--  &8.6931e-5 &1.5633e-4 &4.2084e-4 &1.0877e-3\\
  \hline
\end{tabular}
\label{table P3_1}
\end{table}

Analytical solution of this problem is a combination of two waves of time periods 8 and 16. Corresponding $\sigma$ values at $N_c=1.0$ are 0.393 and 0.196. Thus it is noteworthy to point out that among two stage schemes S2C shows least error despite IRK24 possessing fourth order of accuracy. With higher $N_c$ values least error is reported for computations carried out using S2B set of schemes. This is in conformity with our analysis in preceding sections. All results obtained using S2A shows lesser accuracy and relates to the higher dispersion error predicted in figure \ref{fig:2S_phase}. Interestingly at $N_c=2.5$, S2A reports a large error which may correlate to schemes normalized phase and group velocity deviating furthest from unity for the pair of $hk$ values. For three stage methods all class of schemes apart from S3A report almost identical error at $N_c=1.0$. Here with $N_c$ i.e. $\sigma$ increasing S3B starts producing superior results. But at $N_c=3.0$, advantage of S3B over S3C, S3D and IRK36 gets subsided which is in concurrence with variation of dispersion error presented in figure \ref{fig:3S_phase}.

In figures \ref{fig:P3_1} and \ref{fig:P3_2} numerical solutions obtained by using diverse two stage, and three stage schemes respectively have been presented for $N_c$ = 2.0 and 3.0 along with exact solution. From the figure \ref{fig:P3_1} dispersion error in numerical solution obtained using S2A is evident. For $N_c=2.0$ all other scheme are found to be quite efficient. Nevertheless zoomed view at the position of highest variation reveals better performance of S2B. At $N_c=3.0$ dispersion error is more evident. Here a close look at figure \ref{fig:P3_1}(d) establishes marginally better performance of S2B as per analytical prediction. Performance of all three stage schemes are significantly better for both these CFL numbers. At $N_c=2.0$ performance of S3B, S3C, IRK36 are indistinguishable with S3A producing some deviation from exact solution. Close view at $N_c=3.0$ display slightly better performance of S3B. From this problem it is clear that anticipated superiority of S2B and S3B at relatively higher CFL number are carried over in situations involving composition of waves. This in turn advocates for the philosophy of designing implicit R-K schemes with emphasis on wavenumber regions of interest.

\subsection{Problem 5: Inviscid Burgers' equation}
\begin{figure}[!h]
\begin{minipage}[b]{.6\linewidth}\hspace{-1cm}
\centering\psfig{file=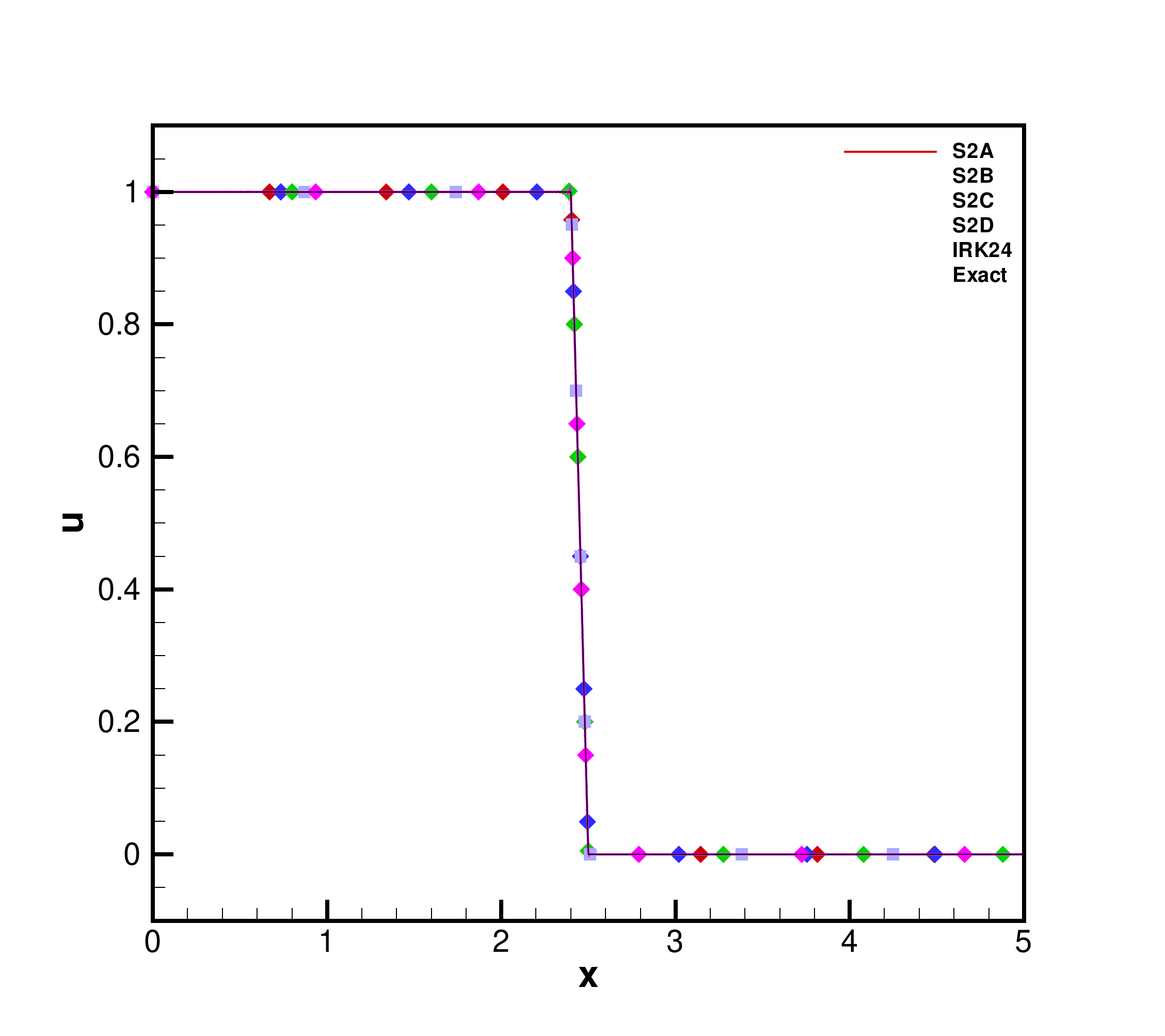,width=0.9\linewidth}\\(a)
\end{minipage}
\begin{minipage}[b]{.6\linewidth}\hspace{-1cm}
\centering\psfig{file=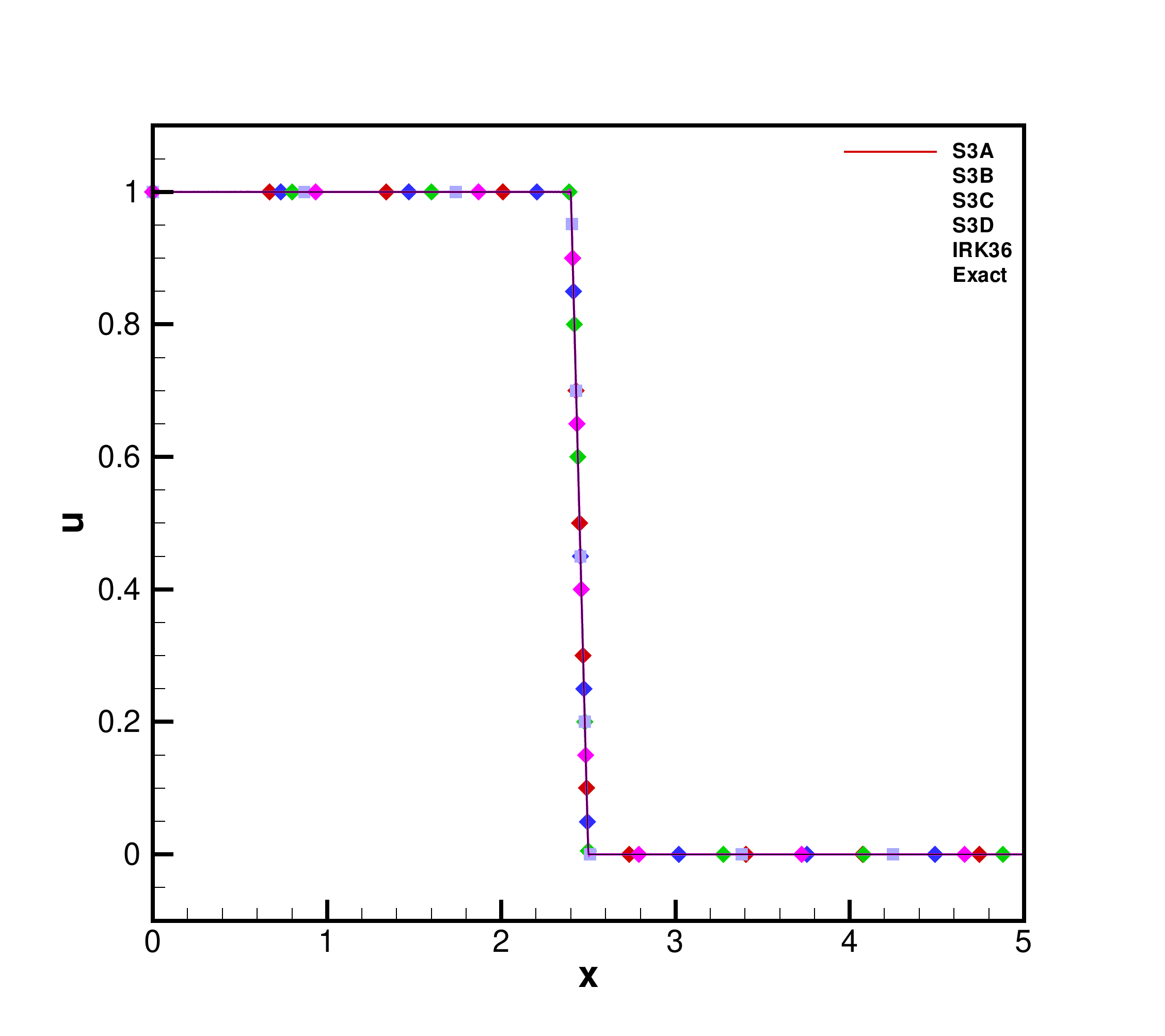,width=0.9\linewidth}\\(b)
\end{minipage}
\begin{minipage}[b]{.6\linewidth}\hspace{-1cm}
\centering\psfig{file=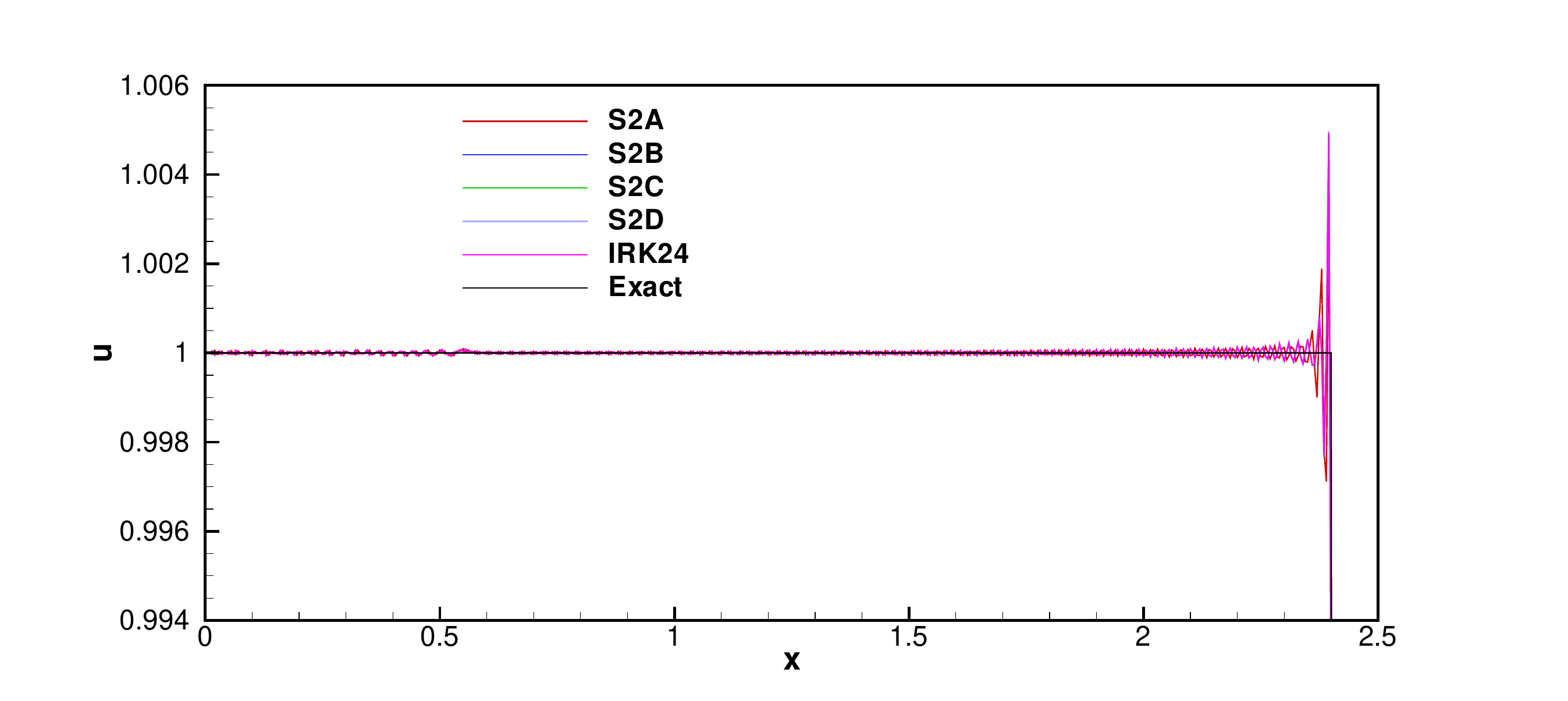,width=0.9\linewidth}\\(c)
\end{minipage}
\begin{minipage}[b]{.6\linewidth}\hspace{-1cm}
\centering\psfig{file=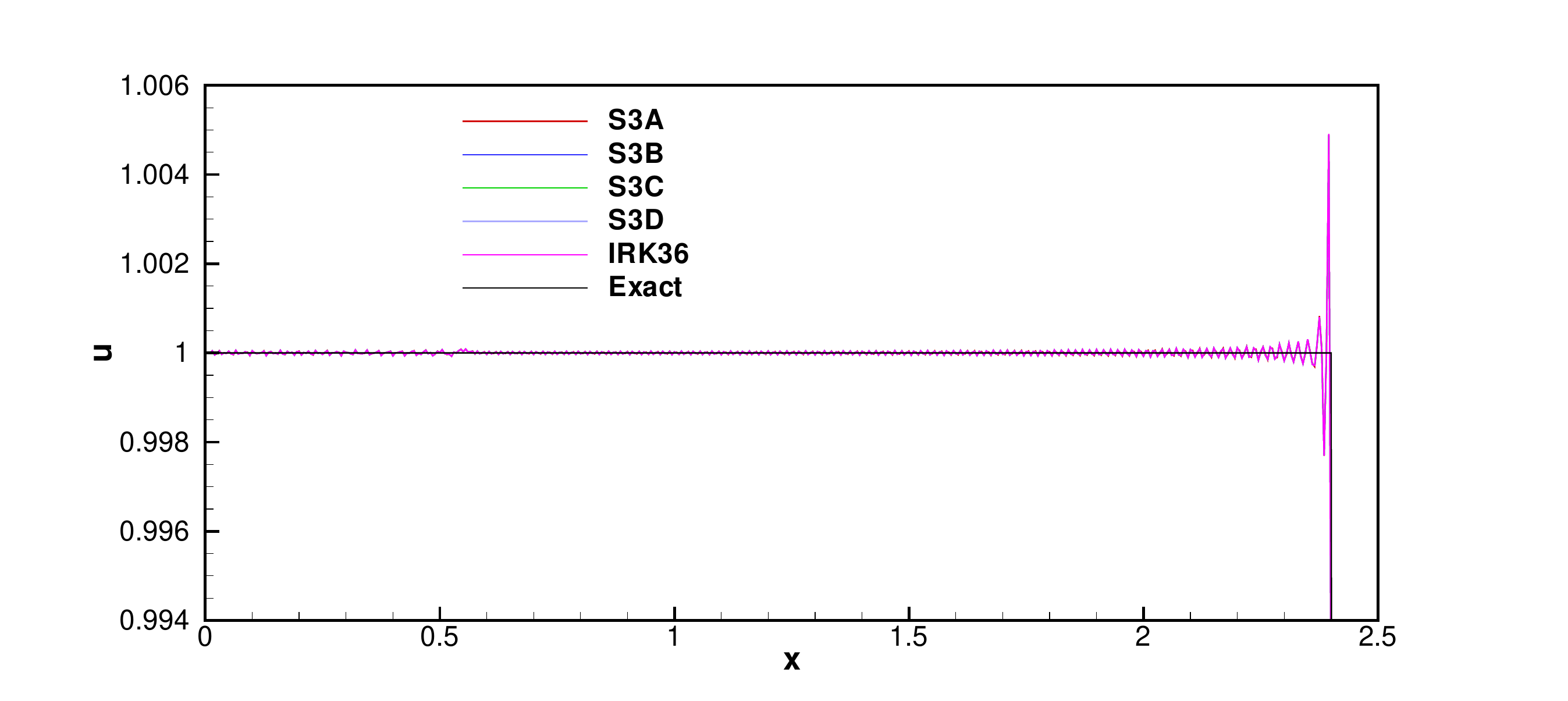,width=0.9\linewidth}\\(d)
\end{minipage}
\begin{center}
\caption{{\sl Problem 5: Numerical solution for $N_c=1.0$ at $t=0.9$ (a) two stage full view, (b) three stage full view, (c) two stage close view, (d) three stage close view.} }
\label{fig:P4_1}
\end{center}
\end{figure}

Next we consider the conservative form of inviscid non-linear Burgers' equation
\begin{eqnarray}\label{P4_1}
\frac{\partial u}{\partial t}+\frac{1}{2}\frac{\partial (u^2)}{\partial x}=0
\end{eqnarray}
with initial condition
\begin{eqnarray*}
u(x,0)=
\begin{cases}
1.0, & x\leq 1.50\\
2.5-x, & 1.50<x\leq 2.50\\
0.0, & x>2.50
\end{cases}
\end{eqnarray*}
defined over the domain $[0, 5]$. This equation, which find its application in high speed flows, displays formation of shock discontinuity with time and is traditionally difficult to capture. We use a fine mesh spacing $\Delta x=0.005$ and compute for $N_c=1.0$. Spatial discretization is carried out using sixth order compact Lele's scheme \cite{lel_92}. At $t=0.9$ exact solution exhibit a very steep shock. In figure \ref{fig:P4_1}(a) we present comparison of exact solution with numerical solutions obtained using different two stage methods. From the figure it is clear that all newly developed numerical schemes are able to correctly predict the steep shock. In fact numerical solutions match quite well with the exact solution highlighting success of newly developed schemes. Similarly in figure \ref{fig:P4_1}(b) one can note the success of all three stage schemes. Zoomed view presented in figures \ref{fig:P4_1}(c) and \ref{fig:P4_1}(d) for two and three stage computations respectively help us closely analyse the error reported by two and three stage methods. Here difference between the exact and computed solutions is evident. Almost identical pattern of high wavenumber grid scale oscillations restricted to a small region near the solution discontinuity is visible for both two and three stage methods in these figures. This may be attributed to numerical characteristics of spatial discretization used. Oscillation are found to be little bit more spread out for S2A and S3A in their respective classes.

\subsection{Problem 6: Convection in 2D}
A two dimensional wave propagation problem given by
\begin{eqnarray}\label{P5_1}
u_t+c u_x+d u_y=0
\end{eqnarray}
with initial condition
\begin{eqnarray}\label{P5_2}
u(x,y,0)=e^{-\frac{(x-x_m)^2+(y-y_m)^2}{b}}\sin (k_x x+k_y y)
\end{eqnarray}
is studied in a square domain $\{(x,y): 0\leq x,y\leq 60\}$ with a uniform grid spacing $\Delta x = 0.1 = \Delta y$. This problem serve as a test case for multi dimensional validation. Discretization of first order space derivatives are are carried out using wide stencil optomized FDo13p \cite{bog_bai_04} in combination with two stage implicit R-K schemes and twelfth order standard finite difference scheme FDs13p \cite{bog_bai_04} in conjunction with three stage schemes. Such choices are made to accentuate portability of the newly developed schemes and to reduce effect of spatial truncation error. We take an identical wave number $2\pi$ and velocity $c=0.5=d$ as also $x_m=30=y_m$. For two stage scheme we compute with a stiff initial wave with $b=0.2$ upto time $t=3.0$ for three different $N_c=N_{c_x}=N_{c_y}$ values. For three stage schemes $b=20.0$ is chosen and solution is computed till $t=7.2$. Temporal accuracy of various schemes discussed in this work is shown using $L^2$-norm error between numerical and exact solutions in table \ref{table P5_1} for two stage and in table \ref{table P5_2} for three stage methods.
\begin{table}[h!]
\caption{Problem 6: $L^2$-norm error between numerical and exact solution at $t=3.0$ for two stage schemes, $b=0.2$.}
\vspace{0.2 cm}
\centering
\begin{tabular}{l c H c H c}
  \hline
  Scheme 	& $N_c=0.4$  & Order 	& $N_c=0.5$  	&Order 	& $N_c=0.6$ \\
  \hline
  S2A	 	&5.0352e-04  &1.72   	&7.3966e-04  	&1.66   &1.0006e-03  \\
  S2B		&8.4053e-06  &6.12  	&3.2916e-05  	&4.72   &7.7888e-05  \\
  S2C	 	&2.5234e-05  &3.94  	&6.0835e-05  	&4.05   &1.2735e-04  \\
  S2D		&3.1786e-05	 &3.58		&7.0696e-05		&3.79	&1.4104e-04	 \\
  IRK24 	&3.1787e-05  &3.58  	&7.0696e-05 	&3.79   &1.4104e-04  \\
  \hline
\end{tabular}
\label{table P5_1}
\end{table}

\begin{table}[h!]
\caption{Problem 6: $L^2$-norm error between numerical and exact solution at $t=7.2$ for three stage schemes, $b=20.0$.}
\vspace{0.2 cm}
\centering
\begin{tabular}{l c H c H c}
  \hline
  Scheme 	& $N_c=0.3$  & Order 	& $N_c=0.6$  	&Order 	& $N_c=0.9$ \\
  \hline
  S3A	 	&1.7805e-07  &4.02   	&2.8985e-06  	&3.87   &1.3927e-05  \\
  S3B		&1.1849e-08  &1.51  	&3.3638e-08  	&8.76   &1.1715e-06  \\
  S3C	 	&6.2670e-09  &4.56  	&1.4773e-07  	&6.55   &2.1006e-06  \\
  S3D		&7.3068e-09	 &4.82		&2.0673e-07		&6.06	&2.4117e-06	 \\
  IRK36 	&7.3060e-09  &4.82  	&2.0673e-07 	&6.06   &2.4117e-06  \\
  \hline
\end{tabular}
\label{table P5_2}
\end{table}

\begin{figure}[!ht]
\begin{minipage}[b]{.6\linewidth}
\centering\psfig{file=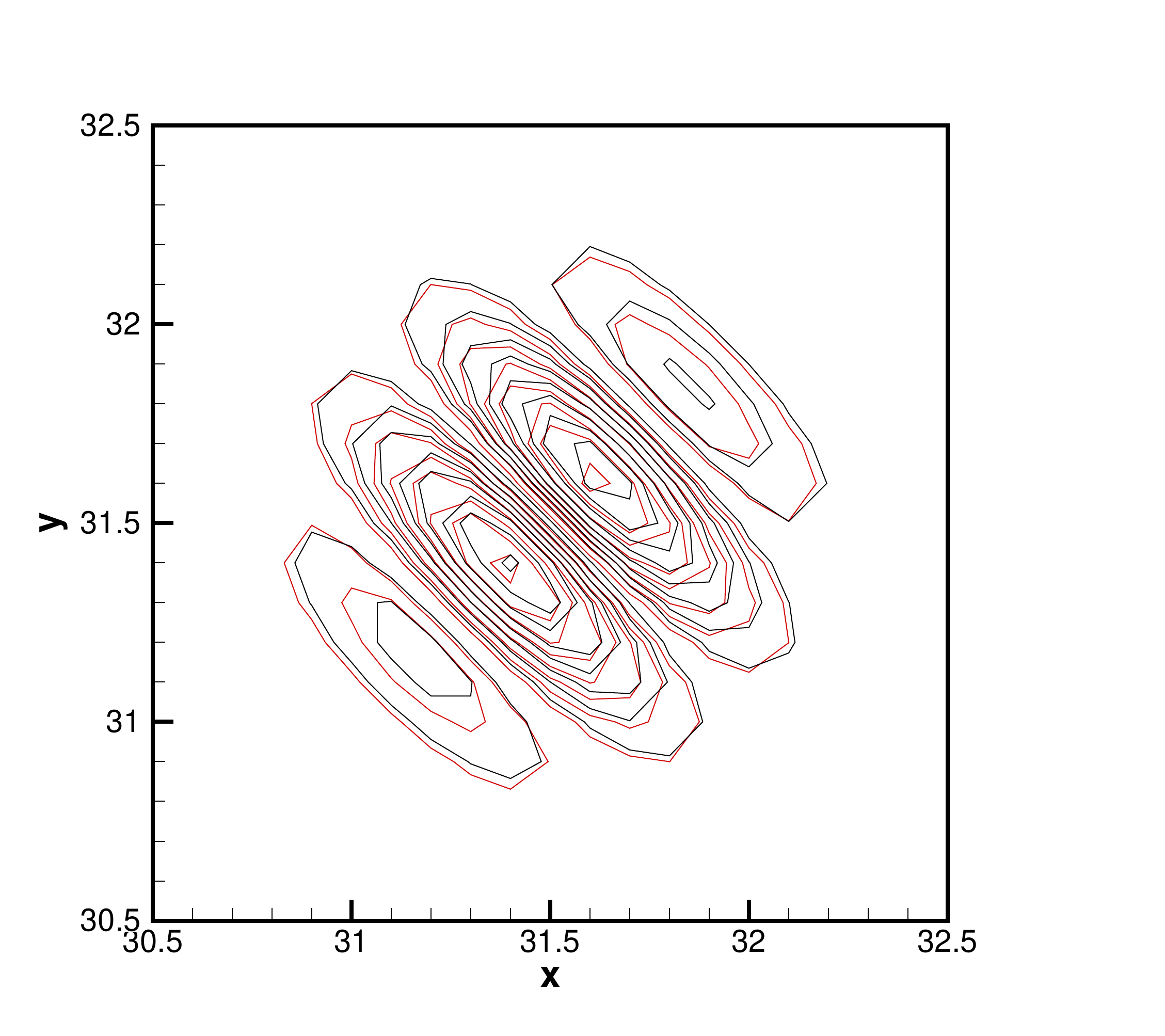,width=0.95\linewidth}
 \\(a)
\end{minipage}            \hspace{-2.5mm}
\begin{minipage}[b]{.6\linewidth}
\centering\psfig{file=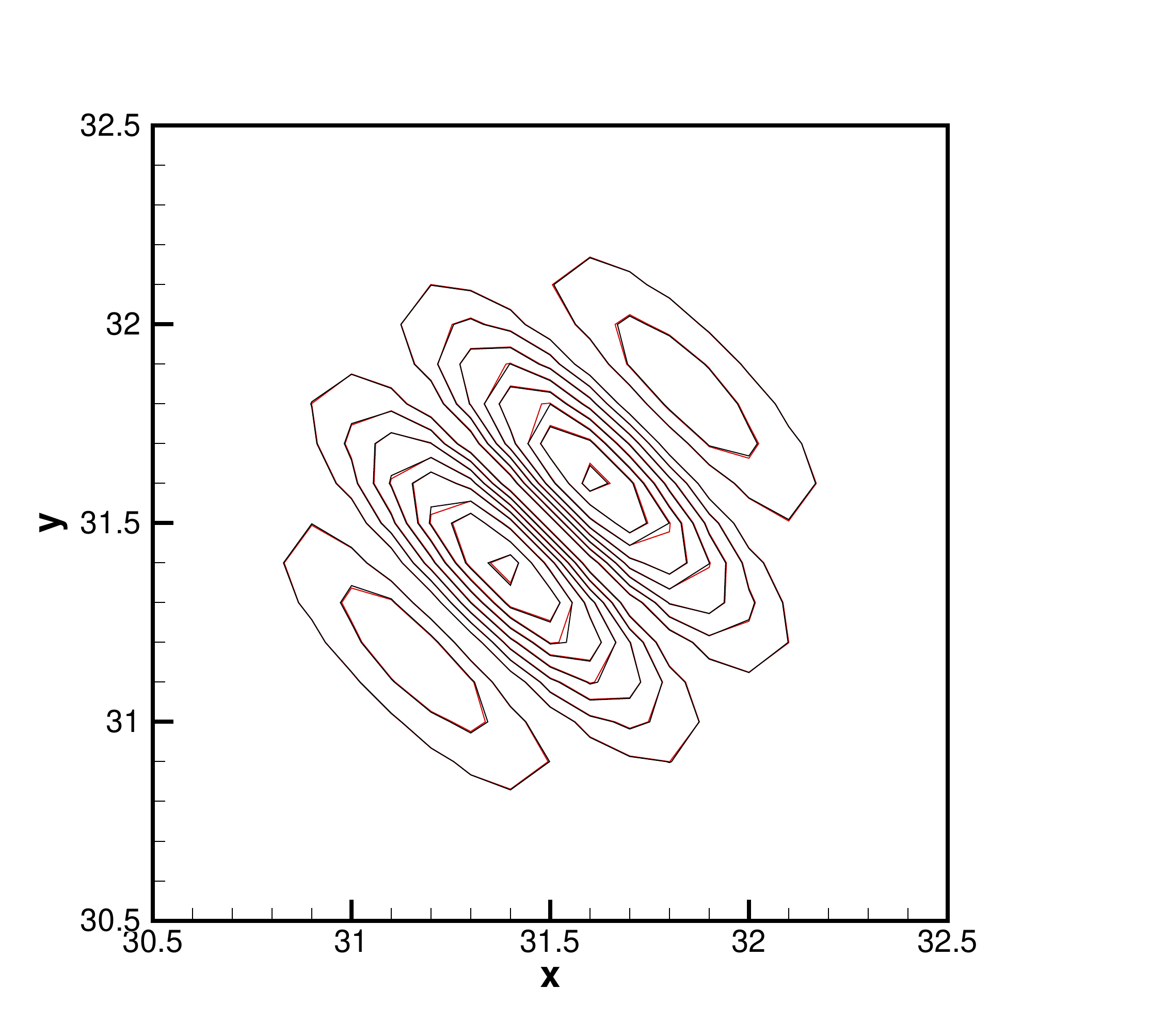,width=0.95\linewidth}
 \\(b)
\end{minipage}            \hspace{-2.5mm}
\begin{minipage}[b]{.6\linewidth}   \
\centering\psfig{file=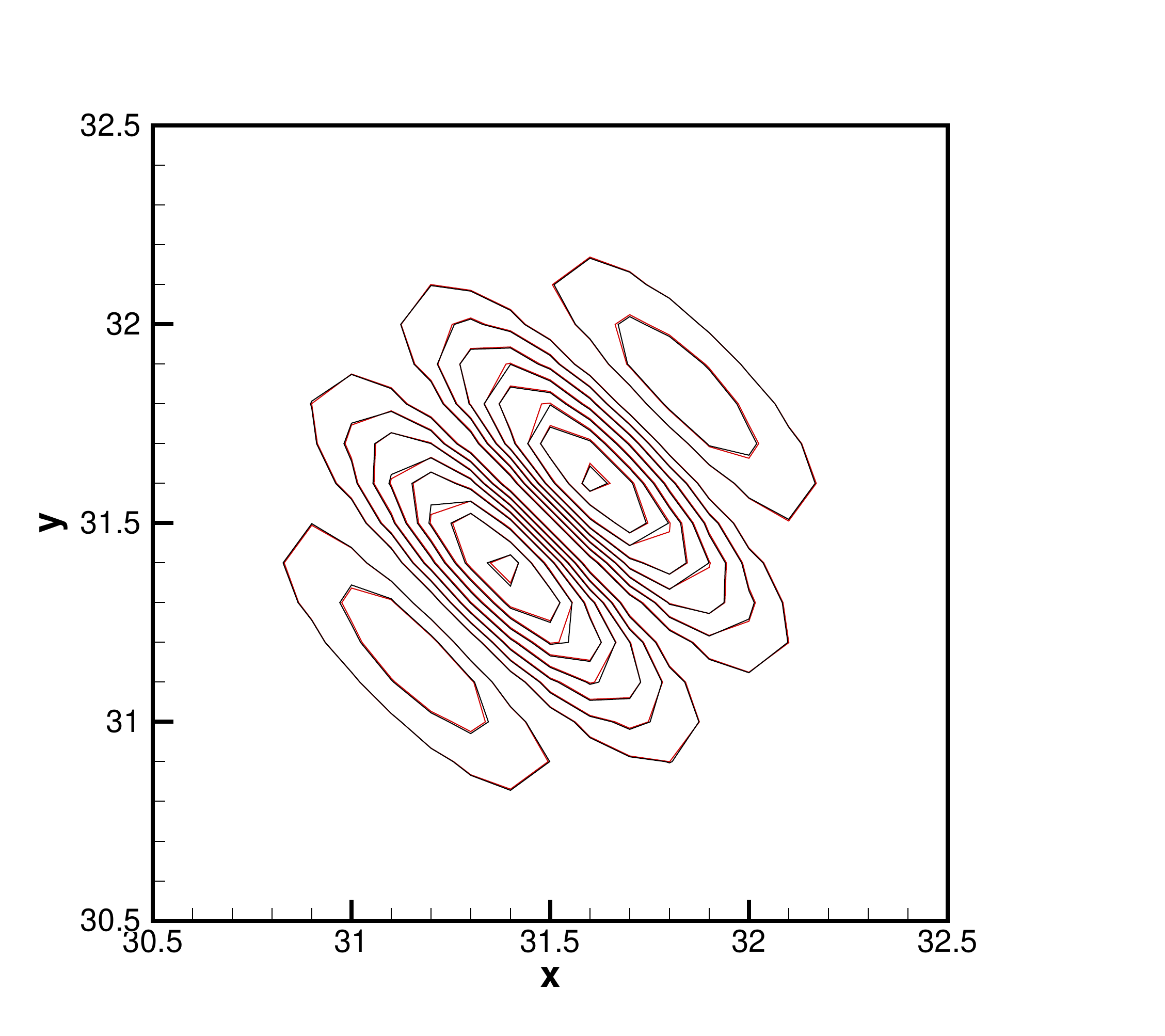,width=0.95\linewidth}
 \\(c)
\end{minipage}            \hspace{-2.5mm}
\begin{minipage}[b]{.6\linewidth}
\centering\psfig{file=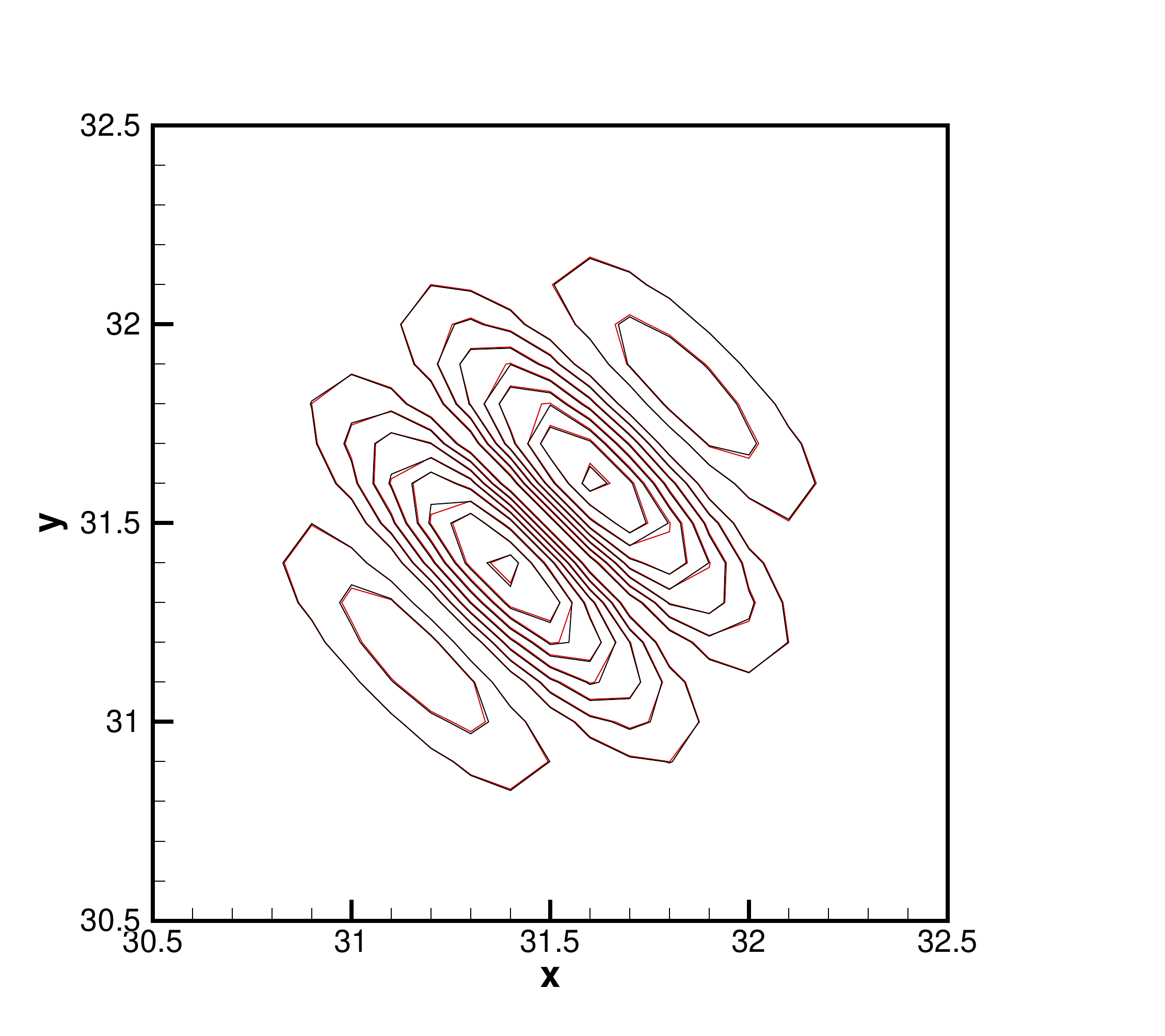,width=0.95\linewidth}
 \\(d)
\end{minipage}            \hspace{-2.5mm}
\begin{center}
\caption{{\sl Problem 6: Comparison of numerical (black) and exact (red) solutions at $t=3.0$ using two stage methods for $b=0.2$ and $N_c=0.6$ (a) S2A, (b) S2B, (c) S2C, and (d) IRK24.} }
\label{fig:P5_1}
\end{center}
\end{figure}

\begin{figure}[!ht]
\begin{minipage}[b]{.6\linewidth}
\centering\includegraphics[width=75mm]{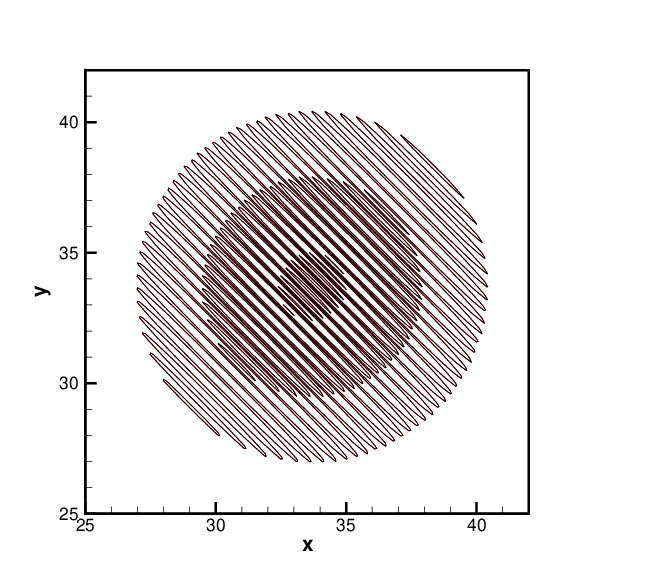}\\(a)
\end{minipage}            \hspace{-2.5mm}
\begin{minipage}[b]{.6\linewidth}
\centering\includegraphics[width=75mm]{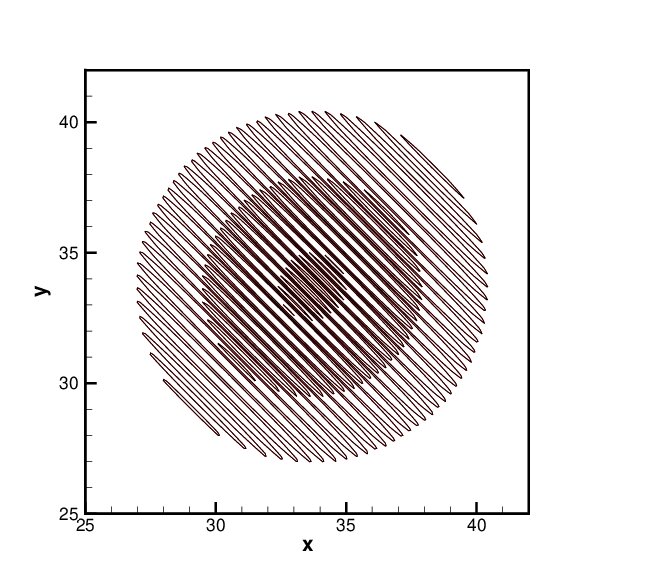}\\(b)
\end{minipage}            \hspace{-2.5mm}
\begin{minipage}[b]{.6\linewidth}   \
\centering\includegraphics[width=75mm]{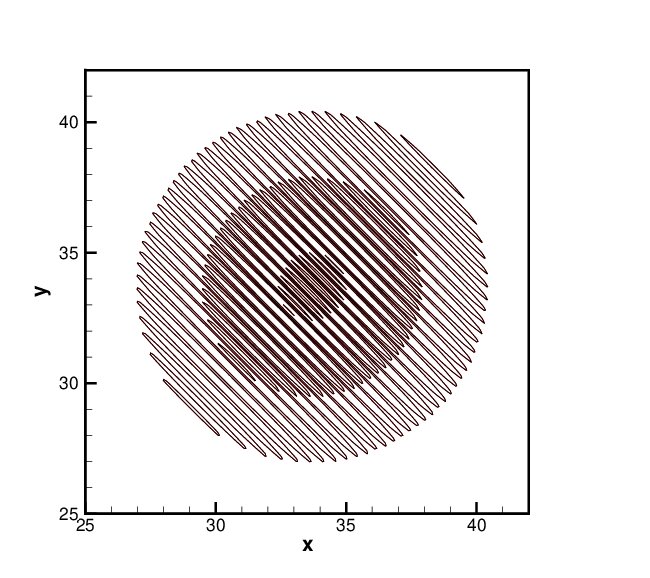}\\(c)
\end{minipage}            \hspace{-2.5mm}
\begin{minipage}[b]{.6\linewidth}
\centering\includegraphics[width=75mm]{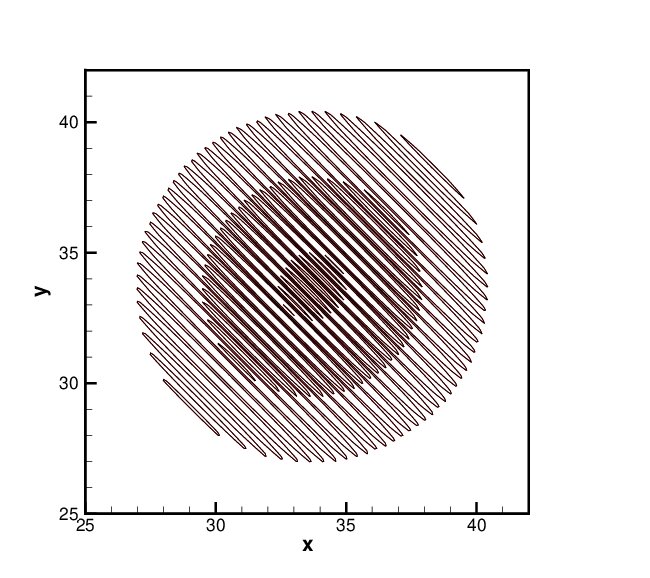}\\(d)
\end{minipage}            \hspace{-2.5mm}
\begin{center}
\caption{{\sl Problem 6: Comparison of numerical (black) and exact (red) solutions at $t=7.2$ using three stage methods for $b=20.0$ and $N_c=0.9$ (a) S3A, (b) S3B, (c) S3C, and (d) IRK36.} }
\label{fig:P5_2}
\end{center}
\end{figure}
For $0.4\le N_c\le 0.6$ we have $0.503\le\sigma\le0.754$ and hence similar to one dimensional cases we notice that S2B holds advantage over all other group of two stage schemes in table \ref{table P5_1}. At $N_c=0.3$, $\sigma=0.377$ and as analytically expected S3C produces least error among all three stage schemes. With increasing $N_c$ values S3B outperform all other schemes in its categories as seen in table \ref{table P5_2}. This is consistent with the schemes analytical emphasis of providing relative weightage to higher $\sigma$ values.

Comparison of numerical and exact solutions obtained using diverse two stage methods are presented in figure \ref{fig:P5_1} for $N_c=0.6$. Marginal superiority of S2B in preserving dispersion character can be seen in this figure. As the error reported are quite small for all three stage methods even at $N_c=0.9$, it is hard to distinguish between numerical and exact solutions for all schemes in figure \ref{fig:P5_2}. Nevertheless it documents efficiency of the newly developed strategy for problems in 2D.

\section{Conclusion}
General two and three stage R-K methods have been thoroughly investigated with a view to minimize dissipation and dispersion error. With enough free parameters amplitude error reduction can be done quite efficiently. Additionally weighted phase error reduction with targeted wavenumber space is found to be possible. Schemes are $A$-stable and even while computation is done with relatively bigger time step they are found to be quite accurate. Subsequently, the algorithm is used to come up with distinct classes of two and three stage implicit R-K methods. It is seen that optimized reduction of dispersion error with emphasis spread equally over the domain $[0, \pi]$ is futile. Schemes developed with comparatively more emphasis for smaller wave number is suggested and are found to be effective for computations done at relatively higher CFL number. Comparison carried out with optimal order schemes available in the literature amply demonstrate efficiency of the newly developed methods. Probably for the first time, numerical characteristics of two and three stage Gauss-Legendre methods, viz. IRK24 and IRK36 are earnestly documented to highlight their inherent potential beyond order of accuracy. These methods, known for their immense accuracy, are found to exhibit very low phase error at relatively lower angular frequency. We found that in two and three stage Gauss-Legendre methods emphasis is restricted to a small neighbourhood of wave number space in the vicinity of zero and relates to the limiting case of the strategy advocated in this study. That phase error could be quantified in terms of integral formula with exponential kernel indeed helped in this course. It can be safely concluded that for error minimization suitable choice of wavenumber domain with appropriate weightage significantly influence characteristics of the schemes and in this context our proposed algorithm can be used to design problem specific time integration formula.

Six different numerical tests, covering one and two dimensional cases including a non-linear problem, are envisaged to illustrate strength of the schemes derived. They also demonstrate efficiency and accuracy of the schemes proposed. Methods with identical characteristics produce comparable results irrespective of order of accuracy. For schemes with negligible dissipation error we see that accuracy decreases with increase in dispersion error and is not always dependent on the order of convergence. It is found that, given a particular wavenumber, error reduction is inherently linked to preservation of dissipation and dispersion characteristics.

\section*{Acknowledgement}
Authors acknowledge use of facilities at the High Performance Computing Centre, Tezpur University sponsored by DeitY, India in collaboration with C-DAC, India.

First author is supported by University Grant Commission, India under Rajiv Gandhi National Fellowship (F1-17.1/2015-16/RGNF-2015-17-SC-WES-12451). Second author is thankful to Science \& Engineering Research Board, India for assistance under Mathematical Research Impact Centric Support  (File Number: MTR/2017/000038).

\section*{References}
\bibliographystyle{plain}

\begin{thebibliography}{}

\bibitem{but_08}
J. Butcher, Numerical Methods for Ordinary Differential Equations, second edition, John Wiley \& Sons Ltd., 2008.

\bibitem{ale_77}
R. Alexander, Diagonally implicit Runge-Kutta methods for stiff O.D.E.'s, SIAM Journal on Numerical Analysis 14 (1977) 1006-1021.

\bibitem{tam_web_93}
C.K.W. Tam, J.C. Webb, Dispersion-Relation-Preserving finite difference schemes for computational acoustics, Journal of Computational Physics 107 (1993) 262-281.

\bibitem{sim_93}
T.E. Simos, Runge-Kutta Interpolants with Minimal Phase-Lag, Computers and Mathematics with Applications 8 (1993) 43-49.

\bibitem{hu_hus_man_96}
F.Q. Hu, M.Y. Hussaini, J.L. Manthey, Low-dissipation and low-dispersion Runge-Kutta schemes for computational acoustics, Journal of Computational Physics 124 (1996) 177-191.

\bibitem{cal_fra_ran_03}
M. Calvo, J.M. Franco, L. R\'{a}ndez, Minimum Storage Runge-Kutta Schemes for Computational Acoustics,  Computers and Mathematics with Applications 45 (2003) 535-545.

\bibitem{bog_bai_04}
C. Bogey, C. Bailly, A family of low dispersive and low dissipative explicit schemes for flow and noise computations, Journal of Computational Physics 194 (2004) 194-214.

\bibitem{ber_bog_bai_06}
J. Berland, C. Bogey, C. Bailly, Low-dissipation and low-dispersion fourth-order Runge-Kutta algorithm, Computers \& Fluids 35 (2006) 1459-1463.

\bibitem{ana_sim_05}
Z.A. Anastassi, T.E. Simos, An optimized Runge-Kutta method for the solution of orbital problems, Journal of Computational and Applied Mathematics 175 (2005) 1-9.

\bibitem{tse_sim_05}
K. Tselios, T.E. Simos, Runge-Kutta methods with minimal dispersion and dissipation for problems arising from computational acoustics, Journal of Computational and Applied Mathematics 175 (2005) 173-181.

\bibitem{fra_gom_ran_97}
J.M. Franco, I. G\'{o}mez, L. R\'{a}ndez, SDIRK methods for stiff ODEs with oscillating solutions, Journal of Computational and Applied Mathematics 81 (1997) 197-209.

\bibitem{naj_mon_13}
A. Najafi-Yazdi, L. Mongeau, A low-dispersion and low-dissipation implicit Runge-Kutta scheme, Journal of Computational Physics 233 (2013) 315-323. 	 	

\bibitem{naz_moh_cha_14}
F. Nazari, A. Mohammadian, M. Charron, A. Zadra, Optimal high-order diagonally-implicit Runge-Kutta schemes for nonlinear diffusive systems on atmospheric boundary layer, Journal of Computational Physics 271 (2014) 118-130.

\bibitem{naz_moh_cha_15}
F. Nazari, A. Mohammadian, M. Charron, High-order low-dissipation low-dispersion diagonally implicit Runge-Kutta schemes, Journal of Computational Physics 286 (2015) 38-48.

\bibitem{fer_spi_08}
L. Ferracina, M.N. Spijker, Strong stability of singly-diagonally-implicit Runge-Kutta methods, Applied Numerical Mathematics 58 (2008) 1675-1686.

\bibitem{bha_sen_sen_13}
S. Bhaumik, S. Sengupta, A. Sengupta, Wave properties of fourth-order fully implicit Runge-Kutta time integration schemes, Computers \& Fluids 81 (2013) 110-121.

\bibitem{har_taa_94}
Z. Haras, S. Ta'asan, Finite difference scheme for long-time integration, Journal of Computational Physics 114 (1994) 265-279.

\bibitem{sen_dip_sau_07}
T.K. Sengupta, A. Dipankar, P. Sagaut, Error dynamics: Beyond von Neumann analysis, Journal of Computational Physics 226 (2007) 1211-1218.

\bibitem{raj_sen_dut_10}
M.K. Rajpoot, T.K. Sengupta, P.K. Dutt, Optimal time advancing dispersion relation preserving schemes, Journal of Computational Physics 229 (2010) 3623-3651.

\bibitem{lel_92}
S.K. Lele, Compact finite difference schemes with spectral-like resolution, Journal of Computational Physics 103 (1992) 16-42.

\bibitem{kot_has_kaj_08}
Lis: Library of Iterative Solvers for Linear Systems, (2017, June 23). Retrieved  from https://www.ssisc.org/lis/

%

\end{thebibliography}

\end{document}